\documentclass[10pt]{article}

\usepackage{amsfonts}
\usepackage{latexsym}
\usepackage{amsthm}
\usepackage{amsmath}
\usepackage{amssymb}
\usepackage{amsopn}
\usepackage{mathrsfs}  
\usepackage[english]{babel}

\newtheorem{theorem}{Theorem}

\newtheorem{coro}{Corollary}
\newtheorem{defi}{Definition}
\newtheorem{hp}{Hypothesis}
\newtheorem{lemma}{Lemma}
\newtheorem{prop}{Proposition}
\newtheorem{rmk}{Remark}

\newcommand{\zerarcounters}{\setcounter{equation}{0}}

\newcommand{\ZZZ}{\mathbb{Z}}
\newcommand{\CCC}{\mathbb{C}}
\newcommand{\NNN}{\mathbb{N}}

\newcommand{\RRR}{\mathbb{R}}
\newcommand{\TTT}{\mathbb{T}}
\newcommand{\SSS}{\mathbb{S}}
\newcommand{\DD}{\mathbb{D}}

\newcommand{\BB}{{\mathcal B}}
\newcommand{\CC}{{\mathcal C}}

\newcommand{\GG}{{\mathcal G}}

\newcommand{\II}{{\mathcal I}}

\newcommand{\MM}{{\mathcal M}}
\newcommand{\NN}{{\mathcal N}}

\newcommand{\PPP}{{\mathcal P}}

\newcommand{\RR}{{\mathcal R}}
\newcommand{\SSSS}{{\mathcal S}}

\newcommand{\np}{\noindent}
\newcommand{\gota}{{\mathfrak a}}

\newcommand{\gotA}{{\mathfrak A}}

\newcommand{\DDD}{{\mathfrak D}}
\newcommand{\gotE}{{\mathfrak E}}

\newcommand{\gotG}{{\mathfrak G}}
\newcommand{\gotH}{{\mathfrak H}}
\newcommand{\gotI}{{\mathfrak I}}

\newcommand{\gotK}{{\mathfrak K}}

\newcommand{\gotM}{{\mathfrak M}}

\newcommand{\gotO}{{\mathfrak O}}
\newcommand{\gotP}{{\mathfrak P}}
\newcommand{\gotQ}{{\mathfrak Q}}
\newcommand{\gotR}{{\mathfrak R}}
\newcommand{\gotS}{{\mathfrak S}}

\newcommand{\gotU}{{\mathfrak U}}

\newcommand{\Val}{{\rm Val}}

\newcommand{\Fullbox}{{\rule{2.0mm}{2.0mm}}}

\newcommand{\EP}{\hfill\Fullbox\vspace{0.2cm}}
\newcommand{\prova}{\noindent{\it Proof. }}
\newcommand{\io}{\infty}
\newcommand{\eps}{\varepsilon}
\newcommand{\al}{\alpha}
\newcommand{\be}{\beta}
\newcommand{\h}{\eta}
\newcommand{\de}{\delta}

\newcommand{\p}{\pi}

\newcommand{\g}{\gamma}
\newcommand{\om}{\omega}
\newcommand{\la}{\lambda}

\newcommand{\s}{\sigma}

\newcommand{\pr}{\partial}
\newcommand{\laa}{\langle}
\newcommand{\raa}{\rangle}

\newcommand{\nn}{{\boldsymbol{\nu}}}

\newcommand{\ii}{{\rm i}}
\newcommand{\xx}{{\bf x}}

\def\ins#1#2#3{\vbox to0pt{\kern-#2 \hbox{\kern#1 #3}\vss}\nointerlineskip}

\begin{document}

\title{{\bf Periodic solutions for a class of nonlinear
partial differential equations in higher dimension}}
\author
{\bf Guido Gentile$^{1}$ and Michela Procesi$^{2}$
\vspace{2mm} \\ \small 
$^{1}$ Dipartimento di Matematica, Universit\`a di Roma Tre, Roma,
I-00146, Italy \\ \small
$^{2}$ Dipartimento di Matematica, Universit\`a di
Napoli ``Federico II'', Napoli, I-80126, Italy \\ \small 
E-mail: gentile@mat.uniroma3.it, procesi@mat.uniroma3.it}

\date{}
\maketitle

\begin{abstract}
We prove the existence of periodic solutions in a class
of nonlinear partial differential equations, including
the nonlinear Schr\"odinger equation, the nonlinear wave equation,
and the nonlinear beam equation, in higher dimension.
Our result covers cases where the bifurcation equation
is infinite-dimensional, such as the nonlinear Schr\"odinger
equation with zero mass, for which solutions which at leading order
are wave packets are shown to exist.
\end{abstract}

%%%%%%%%%%%%%%%%%%%%%%%%%%%%%%%%%%%%%%%%%%%%%%%%%%%%%%%%%%%%%%%%%%%%%%%%%
%%%%%%%%%%%%%%%%%%%%%%%%%%%%%%%%%%%%%%%%%%%%%%%%%%%%%%%%%%%%%%%%%%%%%%%%%
\zerarcounters
\section{Introduction and main results}
\label{sec:1}
%%%%%%%%%%%%%%%%%%%%%%%%%%%%%%%%%%%%%%%%%%%%%%%%%%%%%%%%%%%%%%%%%%%%%%%%%
%%%%%%%%%%%%%%%%%%%%%%%%%%%%%%%%%%%%%%%%%%%%%%%%%%%%%%%%%%%%%%%%%%%%%%%%%

The problem of the existence of finite-dimensional tori for
infinite-dimensional systems, such as nonlinear PDE equations,
has been extensively studied in the literature. Up to very
recent times, the only available results were confined to the
case of one space dimension ($D=1$). In this context the first results
were obtained by Wayne, Kuksin, and P\"oschel \cite{W,Ku,KP,P}, 
for the nonlinear Schr\"odinger equation (NLS) and the
nonlinear wave equation (NLW) with Dirichlet boundary conditions,
by using KAM techniques. Later on,
Craig and Wayne proved similar results, for both Dirichlet and
periodic boundary conditions \cite{CW}, with a rather different
method based on the Lyapunov-Schmidt decomposition.
The case of periodic boundary condition within the framework of
KAM theory was then obtained by Chierchia and You \cite{CY}.
The case of completely resonant systems, i.e. systems where
all eigenvalues of the linear operator are commensurate
with each other, was discussed by several authors,
and theorems on the existence of periodic solutions
for a large measure set of frequencies were obtained
by Bourgain \cite{Bo3} for the NLW with periodic boundary conditions,
by Gentile, Mastropietro and Procesi \cite{GMP}, and Berti and
Bolle \cite{BB1} for the NLW with Dirichlet boundary conditions,
and by Gentile and Procesi \cite{GP1} for the NLS with
Dirichlet boundary conditions. The existence of quasi-periodic
solutions for the completely resonant NLW with periodic boundary
conditions has been proved by Procesi \cite{Pr} 
for a zero-measure set of two-dimensional rotation vectors,
by Baldi and Berti \cite{BB} for a large measure set of
two-dimensional rotation vectors,
and by Yuan \cite{Y1} for a large measure set of
-- at least three-dimensional -- rotation vectors.

Extending the results to higher space dimensions ($D>1$) introduces
a lot of difficulties, mainly due to the high degeneracy of the
eigenvalues of the linear operator. The first achievements
in this direction were due to Bourgain,
and concerned the existence of periodic solutions for NLW \cite{Bo1}
and of periodic solutions (also quasi-periodic in $D=2$) for the
NLS \cite{Bo2}. The case of quasi-periodic solutions
in arbitrary dimension was solved by Bourgain \cite{Bo4} for
the NLS and the NLW. Bourgain's method is based on a Nash-Moser
algorithm, which does not imply the linear stability.

A proof of existence and stability of quasi-periodic
solutions in high dimension was given by 
Geng and You using KAM theory. Their result holds for a class of PDE's,
which includes the nonlinear beam equation (NLB) \cite{GY1}
and the NLS with a smoothing nonlinearity \cite{GY2},
with periodic boundary conditions and with nonlinearities
which do not depend on the space variable. Both conditions are required
in order to ensure a symmetry for the Hamiltonian which simplifies
the problem in a remarkable way. Their approach  does not extend to
the NLS with local nonlinearities -- mainly because it requires a
``second Melnikov condition'' at each iterative KAM step,
and such a condition does not appear
to be satisfied by the local NLS.

Successively, Eliasson and Kuksin \cite{EK}, by using KAM techniques,
proved the existence and stability of quasi-periodic solutions for the
NLS with local nonlinearities. In their paper the main point is
indeed to prove that one may impose a second Melnikov condition
at each iterative KAM step. However, given a PDE equation,
in general (see for instance the case of the NLW in $D>1$),
it can be too hard to impose a second Melnikov condition -- even
on the unperturbed eigenvalues. Very recently, Yuan \cite{Y2}
proposed a KAM-like approach which does not require the
second Melnikov condition, and hence allows to extend the proof
of existence to other kinds of equations, including the NLW:
with respect to Eliasson and Kuksin's approach the linear
stability of the solutions does not follow from the construction.

In both Eliasson and Kuksin's and Yuan's papers Sobolev norms
are used to control the regularity of the solutions in the space
variables, so that only finite smoothness is found
even if the nonlinearity is assumed to be analytic.
This is a drawback which does not arise in Bourgain's approach \cite{Bo4},
where an exponential decay of the Fourier coefficients is obtained.

Again very recently, Berti and Bolle \cite{BB3} proved the
existence of periodic solutions for PDE systems with eigenvalues
of the linear part satisfying rather general separation properties
-- weaker than those considered in this paper.
They use a Nash-Moser algorithm suited for finitely
differentiable nonlinearities, already employed in the
one-dimensional case \cite{BB2}, and they find solutions
belonging to suitable Sobolev classes.
By construction, their method looks for a Sobolev regularity,
and hence it produces only a finite smoothness
even when applied to systems with analytic nonlinearities and
with stronger separation properties,
as in the cases discussed in this paper.
It is very likely that, if we considered analytic nonlinearities
and the same weaker separation properties as in \cite{BB3},
we would obtain solutions with only a finite smoothness.

In this paper we revisit the case of periodic solutions with
a different method, based on renormalisation group ideas
and originally introduced in \cite{GM}.
We consider analytic nonlinearities, and formulate
a general theorem on the existence of periodic solutions
in Gevrey class, which emphasises the main assumptions
that we need in the proof. From a technical point of view,
besides the more abstract formulation -- and hence the wider
range of application --, the present paper represents an
improvement of the renormalisation group method of \cite{GP2},
and allows to considerably simplify the technical aspects of the proof.

For the NLS, with respect to \cite{GP2} and \cite{GY2},
here we remove the condition for the nonlinearity to be smoothed
by a convolution function, so recovering the case of local
nonlinearities, as in \cite{Bo2}.
Moreover, we obtain results for other equations,
including the NLW and the NLB.
Finally -- and this represents the main novelty of this paper --
we discuss cases in which the bifurcation equation
is infinite-dimensional, such as the zero-mass NLS and NLB, where the
other methods have not been applied so far. 
In the resonant case the linearised equation has an
infinite-dimensional space of periodic solutions with the same period,
so that in principle we have at our disposal infinitely many
linear solutions with the same period which can be extended to
solutions of the nonlinear equation. Indeed we find a denumerable
infinity of solutions with the same minimal period even in the
presence of the nonlinearity. More precisely, we prove the
existence of periodic solutions
which at leading order involve an arbitrary finite number of harmonics,
and which therefore can be described as distorted wave packets.
Solutions of this kind are very natural in the case of completely
resonant PDE, where all harmonics are commensurate in the absence
of the nonlinearity. An essential ingredient for the
existence of such solutions is the particular form of
the bifurcation equation: the proof strongly relies on the fact
that the leading order of the nonlinearity is cubic and gauge-invariant.
Moreover, in order to prove the non-degeneracy of the solutions
of the bifurcation equation we need some condition on the higher
orders of the nonlinearity. A sufficient condition is that the
nonlinearity does not depend explicitly on the space variables.

The problem of existence of periodic and quasi-periodic
solutions in completely resonant systems in higher dimension
was already considered by Bourgain in \cite{Bo2},
where he constructed quasi-periodic solutions with two frequencies,
in $D=2$, for the NLS with periodic boundary conditions.
In the case of Dirichlet boundary conditions, proving the non-degeneracy
of the solutions becomes rather involved.
We use a combinatorial lemma, proved in \cite{GP2},
and some results in algebraic number theory.
With respect to the nonlocal NLS considered in \cite{GP2},
the proof we give here is much simpler, however it has the drawback
that a stronger assumption on the nonlinearity is required.

In the remaining part of this section, we give a rigorous
description of the PDE systems we shall consider, and a formal
statement of the results that we shall prove in the paper.
Throughout the paper we shall call a function $F(x,t)$,
with $x=(x_{1},\ldots,x_{D})\in\RRR^{D}$ and $t\in\RRR$,
even [resp. odd] in $x$ -- or even [resp. odd] {\it tout court} --
if it even [resp. odd] in each of its arguments $x_{i}$.

Let $\SSS$ be the $D$ dimensional square $[0,\p]^D$, and let
$\partial \SSS$ be its boundary. We consider for instance
the following class of equations
\begin{equation}
\begin{cases}
(\ii\partial_t+P(-\Delta)+\mu) \, v =  f(x,v,\bar v),
\quad & (x,t)\in \SSS\times \RRR , \\
v(x,t)=0 \quad &  (x,t)\in \partial \SSS\times \RRR , \end{cases}
\label{eq:1.1} \end{equation}
where $\Delta$ is the Laplacian operator, $P(x)$ is a strictly
increasing convex $C^{\io}$ function with $P(0)=0$,
$\mu$ is a real parameter which -- we can assume -- belongs
to some finite interval $(0,\mu_{0})$, with $\mu_{0}>0$, and
$x \to f(x,v(x,t),\bar v(x,t)) $ is an analytic
function which is super-linear in $v,\bar v$ and
odd (in $x$) for odd $v(x,t)$:
\begin{equation}
f(x,v,\bar v) = \!\!\!\!\!\!
\sum_{r,s\in \NNN : r+s\ge N+1} \!\!\!\!\!\!
a_{r,s}(x) \, v^r \bar v^s , \qquad N \ge 1 , 
\label{eq:1.2} \end{equation}
with $a_{r,s}(x)$ even for odd $r+s$ and odd otherwise.
We shall look for odd $2\pi$-periodic solutions
with periodic boundary conditions in $[-\pi,\pi]^{D}$.

We require for $f$ in (\ref{eq:1.2}) to be of the form
\begin{equation}
f(x,v,\bar v)=\frac{\pr}{\pr \bar v}
H(x,v,\bar v) + g(x,\bar v) , \qquad
\overline {H(x,v,\bar v)} = H(x,v,\bar v) .
\label{eq:1.3} \end{equation}

We also consider the class of equations
\begin{equation}
\begin{cases}
\left( \partial_{tt}+(P(-\Delta)+\mu)^{2} \right)
v = f(x,v) , \quad & (x,t)\in \SSS\times \RRR , \\
v(x,t)=0 \quad & (x,t)\in \partial \SSS\times \RRR , 
\end{cases}
\label{eq:1.4} \end{equation}
and finally the wave equation
\begin{equation}
\begin{cases} (\partial_{tt}-\Delta + \mu) \, v = 
f(x,v), \quad & (x,t)\in \SSS\times \RRR , \\
v(x,t)=0 ,\quad & \forall (x,t)\in \partial \SSS\times \RRR ,
\end{cases}
\label{eq:1.5} \end{equation}
where $f(x,v)$ is of the form (\ref{eq:1.2}) with $s$ identically zero
and $a_{r}(x):=a_{r,0}(x)$ real (by parity $a_{r}(x)$ is even for
odd $r$ and odd for even $r$).

We shall consider also (\ref{eq:1.1}), (\ref{eq:1.4})
and (\ref{eq:1.5}) with periodic boundary conditions: in that case,
we shall drop the condition for $f$ to be odd.

For all these classes of equations we prove the existence of
small periodic solutions with frequency $\om$ close to
the linear frequency  $\om_{0}=P(D)+\mu$ for (\ref{eq:1.1})
and (\ref{eq:1.4}) and $\om_{0}= \sqrt{P(D) + \mu}$ for (\ref{eq:1.5}),
with $\om$ in an appropriate Cantor set of positive measure.
We introduce a smallness parameter by rescaling
\begin{equation}
v(x,t)= \eps^{1/N} u(x,\om t) , \qquad
\qquad \eps > 0 ,
\label{eq:1.6}
\end{equation}
with $\om = P(D) + \mu - \eps$ for (\ref{eq:1.1}) and (\ref{eq:1.4})
and $\om^{2} = P(D) + \mu - \eps$ for (\ref{eq:1.5}).

We shall formulate our results in a more abstract context, by
considering the following classes of equations with Dirichlet
boundary conditions:
\begin{subequations}
\begin{align}
{\rm (I)} & \qquad \begin{cases} \DD(\eps) \, u = \eps f(x,u,
\bar u,\eps^{1/N}) , \qquad (x,t) \in \SSS\times \TTT , \\ 
u(x,t)=0 , \qquad (x,t)\in \pr \SSS\times \TTT , \end{cases}
\label{eq:1.7a} \\
{\rm (II)} & \qquad \begin{cases}  \DD(\eps) \, u = \eps
f(x,u,\eps^{1/N}) , \qquad (x,t) \in \SSS\times \TTT , \\ 
u(x,t)=0 , \qquad (x,t)\in \pr \SSS\times \TTT , \end{cases}
\label{eq:1.7b}
\end{align}
\label{eq:1.7}
\end{subequations}
\vskip-.3truecm
\noindent where $\TTT:= \RRR/2\p\ZZZ$ and $\DD(\eps)$ is
a linear (possibly integro-)differential wave-like operator
with constant coefficients depending on a (fixed once and
for all) real parameter $\om_{0}$ and on the parameter $\eps$.

We can treat the case of periodic boundary conditions in the same way:
\begin{subequations}
\begin{align}
{\rm (I)} & \qquad \DD(\eps) \, u = \eps f(x,u,\bar u,\eps^{1/N}) ,
\qquad (x,t)\in \TTT^{D} \times \TTT ,
\label{eq:1.8a }\\
{\rm (II)} & \qquad \DD(\eps) \, u = \eps f(x,u,\eps^{1/N}) ,
\qquad (x,t)\in \TTT^{D} \times \TTT ,
\label{eq:1.8b}
\end{align}
\label{eq:1.8}
\end{subequations}
\vskip-.5truecm
\noindent with the same meaning of the symbols as in (\ref{eq:1.7}).

In Case (I) we assume that $f(x,u,\bar u,\eps^{1/N})$ is a rescaling
of a function $f(x,u,\bar u)$ defined as in (\ref{eq:1.2}) and
satisfying (\ref{eq:1.3}).
In Case (II) we suppose  $\DD(\eps)$ real  and  $f$  real for real
$u$, so that it is natural to look for real solutions $u=\bar u$.

For $\nn\in\ZZZ^{D+1}$ set $\nn=(\nu_{0},m)$, with
$\nu_{0} \in \ZZZ$ and $m=(\nu_{1},\ldots,\nu_{D})\in\ZZZ^{D}$
and $|\nn|=|\nu_{0}|+|m|=|\nu_{0}|+|\nu_{1}|+\ldots+|\nu_{D}|$.
For $\xx=(t,x)=(t,x_{1},\ldots,x_{D})\in\RRR^{D+1}$ set
$\nn\cdot\xx=\nu_{0}t+m \cdot x = \nu_{0} t + \nu_{1} x_{1} +
\ldots + \nu_{D}x_{D}$. Set also $\ZZZ_{+}=\{0\}\cup\NNN$ and
$\ZZZ^{D+1}_{*} = \ZZZ^{D+1} \setminus \{{\bf 0}\}$.
Finally denote by $\delta(i,j)$ the Kronecker delta, i.e.
$\delta(i,j)=1$ if $i=j$ and $\delta(i,j)=0$ otherwise.
Given a finite set $\gotA$ we denote by $|\gotA|$
the cardinality of the set. Throughout the paper, for $z\in\CCC$
we denote by $\overline z$ the complex conjugate of $z$.

Since all the results of the paper are local (that is, they concern
small amplitude solutions), we shall always assume that
the hypotheses below are satisfied for all $\eps$ sufficiently small.

%%%%%%%%%%%%%%%%%%%%%%%%%%%%%%%%%%%%%%%%%%%%%%%%%%%%%%%%%%%%%%%%%%%%%%%%%
\begin{hp} \label{hp:1}
\textbf{ (Conditions on the linear part).}
\begin{enumerate}
\item $\DD(\eps)$ is diagonal in the Fourier basis
$\{{\rm e}^{\ii\nn\cdot \xx}\}_{\nn\in \ZZZ^{D+1}}$ with real
eigenvalues $\de_{\nn}(\eps)$ which are $C^{\io}$
in both $\nn$ and $\eps$.
\item For all $\nn\in\ZZZ^{D+1}_{*}$ one has either
$\de_{\nn}(0)=0$ or $|\de_{\nn}(0)| \ge \g_{0}|\nn|^{-\tau_{0}}$,
for suitable constants $\g_{0},\tau_{0}>0$.
\item For all $\nn\in\ZZZ^{D+1}_{*}$ one has
$\left| \pr_{\eps} \de_{\nn}(\eps) \right| <
c_{2} |\nn|^{c_{0}}$ and, if $|\de_{\nn}(\eps)|<1/2$,
one has $\left| \pr_{\eps} \de_{\nn}(\eps) \right| >
c_{1} |\nn|^{c_{0}}$ as well, for suitable
$\eps$-independent constants $c_{0},c_{1},c_{2}>0$.
\item For all $\nn\in\ZZZ^{D+1}_{*}$
such that $|\de_\nn(\eps)|< 1/2$ one has 
$ \left| \pr_{\eps} \pr_{\nn} \de_{\nn}(\eps) \right|
\le c_{3} |\nn|^{c_{0}-1}$,
for a suitable $\eps$-in\-de\-pen\-dent constant $c_{3}>0$.
\item In case (I) we require that if for some $\eps$ and for some
$\nn_{1},\nn_{2}\in \ZZZ^{D+1}$ one  has
$|\de_{\nn_{1}}(\eps)|,|\de_{\nn_2}(\eps)|<1/2$
then $|\nn_{1}-\nn_{2}|\le |\nn_{1}+\nn_{2}|$.
\end{enumerate}
\end{hp}
%%%%%%%%%%%%%%%%%%%%%%%%%%%%%%%%%%%%%%%%%%%%%%%%%%%%%%%%%%%%%%%%%%%%%%%%%

We now pass to the equation for the Fourier coefficients. We write
\begin{equation}
u(x,t) = \sum_{\nn \in \ZZZ^{D+1}}
u_{\nn} \, {\rm e}^{\ii\nn \cdot \xx } ,
\label{eq:1.9}
\end{equation}
and introduce the coefficients $u^{\pm}_{\nn}$ by setting
$u^{+}_{\nn}:=u_{\nn}$ and $u^{-}_{\nn}:=\overline{u_{\nn}}$.
Analogously we define
\begin{equation}
f_{\nn}(\{u\},\eta) :=
\left[ f(x,u,\bar u,\eta) \right]_{\nn} =
\!\!\!\!\!\!\!\!\! \sum_{r,s\in \NNN : r+s = N+1} \!\!\!\!\!\!\!\!\!\!\!
[a_{r,s}(x)u^r \bar u^s]_{\nn} \;\;+ \!\!\!\!\!\!\!\!\!
\sum_{r,s\in \NNN : r+s > N+1} \!\!\!\!\!\!\!\!\!\!\!\!\! \eta^{r+s-N-1}
[a_{r,s}(x)u^r \bar u^s]_{\nn}
\nonumber
\end{equation}
where $\{u\} = \{u^{\s}_{\nn}\}^{\s=\pm}_{\nn\in\ZZZ^{D+1}}$,
$[\cdot]_{\nn}$ denotes the Fourier coefficient with label $\nn$,
and we set $f^{+}_{\nn}:=f_{\nn}$ and $f^{-}_{\nn}:=\overline{f_{\nn}}$.
Naturally $f_{\nn}$ depends also on the Fourier coefficients of the
functions $a_{r,s}(x)$, which we denote by $a_{r,s,m}$,
with $m\in \ZZZ^{D}$; we set $a^{+}_{r,s,m} := a_{r,s,m}$ and
$a^{-}_{r,s,m} := \overline{a_{r,s,m}}$.

Then in Fourier space the equations (\ref{eq:1.7})
and (\ref{eq:1.8}) give
\begin{equation}
\de_{\nn}(\eps) \, u^{\s}_{\nn} =
\eps f^{\s}_{\nn}(\{u\},\eps^{1/N}),
\qquad \nn \in \ZZZ^{D+1} , \qquad \s = \pm ,
\label{eq:1.10} \end{equation}
and in the case of Dirichlet boundary conditions we shall require
$u_{\nn} = - u_{S_i(\nn)}$ for all $i=1,\dots,D$, where $S_{i}(\nn)$
is the linear operator that changes the sign
of the $i$-th component of $\nn$.

%%%%%%%%%%%%%%%%%%%%%%%%%%%%%%%%%%%%%%%%%%%%%%%%%%%%%%%%%%%%%%%%%%%%%%%%%%
\begin{rmk} \label{rmk:1}
The reality condition on $H$ in (\ref{eq:1.3}) spells
\begin{equation}
\left( s+1 \right) a^{-}_{s+1,r-1,m}= r \, a^{+}_{r,s,-m} .
\label{eq:1.11}
\end{equation}
Moreover, by the analyticity assumption on the nonlinearity, one has
$|a_{r,s,m}| \le A_{1}^{r+s} {\rm e}^{- A_{2}|m|}$ for suitable
positive constants $A_{1}$ and $A_{2}$ independent of $r$ and $s$.
\end{rmk}
%%%%%%%%%%%%%%%%%%%%%%%%%%%%%%%%%%%%%%%%%%%%%%%%%%%%%%%%%%%%%%%%%%%%%%%%%%

%%%%%%%%%%%%%%%%%%%%%%%%%%%%%%%%%%%%%%%%%%%%%%%%%%%%%%%%%%%%%%%%%%%%%%%%%%  
\begin{rmk} \label{rmk:2}
We have doubled our equations by considering separately the
equations for $u_{\nn}^{+}$ and $u_{\nn}^{-}$ -- which clearly must
satisfy a compatibility condition. In Case (II) one can work
only on $u_{\nn}^{+}$, since $u^{-}_{\nn}=u^{+}_{-\nn}$. In
other examples it may be possible to reduce to solutions with
$u_{\nn}$ real for all $\nn\in\ZZZ^{D+1}$, but we found more
convenient to introduce the doubled equations in order to
deal with the general case.
\end{rmk}
%%%%%%%%%%%%%%%%%%%%%%%%%%%%%%%%%%%%%%%%%%%%%%%%%%%%%%%%%%%%%%%%%%%%%%%%%%  

Following the standard Lyapunov-Schmidt decomposition scheme we split
$\ZZZ^{D+1}$ into two subsets called $\gotP$ and $\gotQ$ and treat
the equations separately. By definition we call $\gotQ$ the set of
those $\nn\in \ZZZ^{D+1}$ such that $ \de_{\nn}(0)=0$; then we define
$\gotP=\ZZZ^{D+1}\setminus \gotQ$.
The equations (\ref{eq:1.10}) restricted to the $\gotP$ and $\gotQ$
subset are called respectively the $P$ and $Q$ equations.

%%%%%%%%%%%%%%%%%%%%%%%%%%%%%%%%%%%%%%%%%%%%%%%%%%%%%%%%%%%%%%%%%%%%%%%%%
\begin{hp} \label{hp:2}
\textbf{ (Conditions on the $\boldsymbol Q$ equation).}
\begin{enumerate}
\item For all $\nn\in\gotQ$ one has $\la_{\nn}(\eps):=
\eps^{-1}\de_{\nn}(\eps) \ge c>0$, where $c$ is $\eps$-independent. 
%\item  The set $\gotQ$ is finite dimensional.
\item  The $Q$ equation at $\eps=0$, 
$$ \la_{\nn}(0) \, u^{\s}_{\nn}=f^{\s}_{\nn}(\{u\},0) ,
\qquad \nn\in\gotQ , \qquad \s=\pm 1, $$
has a non-trivial non-degenerate solution
$$ q^{(0)}(x,t) = \sum_{\nn\in\gotQ} u^{(0)}_{\nn}
{\rm e}^{\ii\nn\cdot\xx} , $$
where non-degenerate means that the matrix
$$ J^{\s,\s'}_{\nn,\nn'} = \la_{\nn}(0) \, \de(\nn,\nn')\,\de(\s,\s') -
\frac{\partial f^{\s}_{\nn}}{\partial u^{\s'}_{\nn'}}(\{q^{(0)}\},0) $$
is invertible. Moreover one has $|u^{(0)}_{\nn}| \le \Lambda_{0} {\rm e}^{-
\lambda_{0}|\nn|}$ and $\left| (J^{-1})^{\s,\s'}_{\nn,\nn'} \right|
\le \Lambda_{0} {\rm e}^{-\lambda_{0}|\nn-\nn'|}$,
for suitable constants $\Lambda_{0}$ and $\lambda_{0}$.
\end{enumerate}
\end{hp}
%%%%%%%%%%%%%%%%%%%%%%%%%%%%%%%%%%%%%%%%%%%%%%%%%%%%%%%%%%%%%%%%%%%%%%%%%

%%%%%%%%%%%%%%%%%%%%%%%%%%%%%%%%%%%%%%%%%%%%%%%%%%%%%%%%%%%%%%%%%%%%%%%%%%
\begin{rmk} \label{rmk:3}
The solution of the {\rm bifurcation equation}, i.e. of the $Q$ equation
at $\eps=0$, could be assumed to be only Gevrey-smooth.
Note also that, even when $\gotQ$
is infinite-dimensional, the number of non-zero Fourier components
of $q^{(0)}(x,t)$ can be finite.
\end{rmk}
%%%%%%%%%%%%%%%%%%%%%%%%%%%%%%%%%%%%%%%%%%%%%%%%%%%%%%%%%%%%%%%%%%%%%%%%%%

%%%%%%%%%%%%%%%%%%%%%%%%%%%%%%%%%%%%%%%%%%%%%%%%%%%%%%%%%%%%%%%%%%%%%%%%%
\begin{defi} \label{def:1}
\textbf{ (The sets $\boldsymbol \gotE_{\boldsymbol 0}$, $\boldsymbol
\gotO \boldsymbol( \boldsymbol \eps \boldsymbol)$
and $\boldsymbol \gotO$).} 
Given $\eps\in \gotE_{0}:=[0,\eps_0]$  we set
$\gotO(\eps):=\{ \nn \in \gotP :  \left| \de_{\nn}(\eps) \right| <
1/2 \}$ and $\gotO =\cup_{\eps\in\gotE_{0}} \gotO(\eps)$.
Finally we call $\gotR$ the subset $\gotP\setminus\gotO$.
\end{defi}
%%%%%%%%%%%%%%%%%%%%%%%%%%%%%%%%%%%%%%%%%%%%%%%%%%%%%%%%%%%%%%%%%%%%%%%%%

%%%%%%%%%%%%%%%%%%%%%%%%%%%%%%%%%%%%%%%%%%%%%%%%%%%%%%%%%%%%%%%%%%%%%%%%%%
\begin{rmk} \label{rmk:4}
Note that $\nn\in\gotR$ means that $|\de_{\nn}(\eps)| \ge 1/2$
for all $\eps\in\gotE_{0}$.
\end{rmk}
%%%%%%%%%%%%%%%%%%%%%%%%%%%%%%%%%%%%%%%%%%%%%%%%%%%%%%%%%%%%%%%%%%%%%%%%%%

The following definitions appear (in a slightly different form)
in the papers by Bourgain. The notations which we use are
those proposed by Berti  and Bolle in \cite{BB3}. 

%%%%%%%%%%%%%%%%%%%%%%%%%%%%%%%%%%%%%%%%%%%%%%%%%%%%%%%%%%%%%%%%%%%%%%%%%
\begin{defi} \label{def:2}
\textbf{ (The equivalence relation $\boldsymbol\sim$).}
We say that two vectors $\nn,\nn'\in\gotO(\eps)$
are equivalent, and we write $\nn\sim\nn'$, if for $\beta$
small enough the following happens: one has $|\delta_{\nn}(\eps)|,
|\delta_{\nn'}(\eps)|< 1/2$ and there exists a sequence
$\{\nn_{1},\ldots,\nn_{N}\}$ in $\gotO(\eps)$, with
$\nn_{1}=\nn$ and $\nn_{N}=\nn'$, such that
\begin{equation}
\left| \delta_{\nn_{k}}(\eps) \right|< \frac{1}{2} , \qquad
\left| \nn_{k}-\nn_{k+1} \right| \le \frac{C_{2}}{2}
\left(|\nn_{k}|+|\nn_{k+1}| \right)^{\beta} , \qquad
k=1,\ldots,N-1, \nonumber]
\end{equation}
where $C_{2}$ is a universal constant. Denote by $\Delta_{j}(\eps)$,
$j\in\NNN$, the equivalence classes with respect to $\sim$.
\end{defi}
%%%%%%%%%%%%%%%%%%%%%%%%%%%%%%%%%%%%%%%%%%%%%%%%%%%%%%%%%%%%%%%%%%%%%%%%%%

%%%%%%%%%%%%%%%%%%%%%%%%%%%%%%%%%%%%%%%%%%%%%%%%%%%%%%%%%%%%%%%%%%%%%%%%%%
\begin{rmk} \label{rmk:5}
The equivalence relation $\sim$ induces a partition of $\gotO(\eps)$
into disjoint sets $\{\Delta_{j}(\eps)\}_{j\in \NNN}$.
Note also that, if $\nn,\nn'\in \Delta_{j}(\eps)$, then it is not
possible that for some $\eps'$ one has $\nn\in\Delta_{j_{1}}(\eps')$
and $\nn'\in \Delta_{j_{2}}(\eps')$ with $j_{1}\neq j_{2}$.
\end{rmk}
%%%%%%%%%%%%%%%%%%%%%%%%%%%%%%%%%%%%%%%%%%%%%%%%%%%%%%%%%%%%%%%%%%%%%%%%%

%%%%%%%%%%%%%%%%%%%%%%%%%%%%%%%%%%%%%%%%%%%%%%%%%%%%%%%%%%%%%%%%%%%%%%%%%
\begin{hp} \label{hp:3}
\textbf{ (Conditions on the set $\boldsymbol \gotO \boldsymbol(
\boldsymbol\eps \boldsymbol)$: separation properties).} 
There exist three $\eps$-independent positive constants
$\al,\beta,C_{1}$, with $\al$ small enough and $\beta<\al$,
such that $|\Delta_{j}(\eps)| \le C_{1} p_{j}^{\al}(\eps)$, where
$p_{j}(\eps)= \min_{\nn\in\Delta_{j}(\eps)}|\nn|$,
for all $j\in\NNN$.
\end{hp}
%%%%%%%%%%%%%%%%%%%%%%%%%%%%%%%%%%%%%%%%%%%%%%%%%%%%%%%%%%%%%%%%%%%%%%%%%

%%%%%%%%%%%%%%%%%%%%%%%%%%%%%%%%%%%%%%%%%%%%%%%%%%%%%%%%%%%%%%%%%%%%%%%%%
\begin{rmk} \label{rmk:6}
Hypothesis \ref{hp:3} implies the following properties:
\begin{eqnarray}
& & {\rm dist}(\Delta_{j}(\eps),\Delta_{j'}(\eps)) \ge \frac{C_{2}}{2}
\left( p_{j}(\eps) + p_{j'}(\eps) \right)^{\beta}
\qquad \forall j,j'\in\NNN \text{ such that } j\neq j' , \nonumber \\
& & {\rm diam}(\Delta_{j}(\eps))\le C_{1} C_{2} p_{j}^{\al+\beta}(\eps)
, \qquad \max_{\nn\in\Delta_{j}(\eps)} |\nn| \le 2 p_{j}(\eps)
\qquad \forall j\in\NNN , \nonumber
\end{eqnarray}
and, furthermore, we can always assume that $2^{c_{0}-1}
C_{1}C_{2}p_{j}^{\al+\beta} \le \zeta p_{j}$,
with $\zeta c_{3} < c_{1}/4$, where the
constants $c_{1}$ and $c_{3}$ are defined in Hypothesis \ref{hp:1}.
\end{rmk}
%%%%%%%%%%%%%%%%%%%%%%%%%%%%%%%%%%%%%%%%%%%%%%%%%%%%%%%%%%%%%%%%%%%%%%%%%

%%%%%%%%%%%%%%%%%%%%%%%%%%%%%%%%%%%%%%%%%%%%%%%%%%%%%%%%%%%%%%%%%%%%%%%%%
\begin{rmk} \label{rmk:7}
Given $N>0$ and for all $\eps$ outside a finite set (depending on $N$)
the sets $\Delta_{j}(\eps) \cap\{\nn:|\nn|\le N\}$ are locally constant,
namely for all $\bar\eps$ outside a finite set there exists an
interval $\gotI$ such that $\bar\eps\in\gotI$ with the following
property: There exists an $\eps$-independent numbering of the sets
$\Delta_{j}(\eps)$ contained in $\{\nn : |\nn|\le N\}$ so that
$\Delta_{j}(\eps)=\Delta_{j}(\bar\eps)$ for all $\eps\in\gotI$.
\end{rmk}
%%%%%%%%%%%%%%%%%%%%%%%%%%%%%%%%%%%%%%%%%%%%%%%%%%%%%%%%%%%%%%%%%%%%%%%%%%

We can now state our main result.

%%%%%%%%%%%%%%%%%%%%%%%%%%%%%%%%%%%%%%%%%%%%%%%%%%%%%%%%%%%%%%%%%%%%%%%%%
\begin{theorem} \label{thm:1}
Consider an equation in the class described by (\ref{eq:1.7})
and (\ref{eq:1.8}), such that the Hypotheses \ref{hp:1},
\ref{hp:2} and \ref{hp:3} hold.
There exist a positive constant $\eps_{0}$
and a Cantor set $\gotE \subset [0,\eps_{0}]$, such that
for all $\eps\in \gotE$ the equation admits a solution $u(x,t)$,
which is $2\pi$-periodic in time and Gevrey-smooth both
in time and in space, and such that
\begin{equation}
\left| u(x,t) -  q^{(0)}(x,t) \right|\le C \eps^{1/N} ,
\nonumber
\end{equation}
uniformly in $(x,t)$. The set $\gotE$ has positive Lebesgue measure and
\begin{equation}
\lim_{\eps\to 0^{+}} \frac{{\rm meas}(\gotE \cap [0,\eps])}{\eps}=1 ,
\label{eq:1.12} \end{equation}
where ${\rm meas}$ denotes the Lebesgue measure.
\end{theorem}
%%%%%%%%%%%%%%%%%%%%%%%%%%%%%%%%%%%%%%%%%%%%%%%%%%%%%%%%%%%%%%%%%%%%%%%%%

%%%%%%%%%%%%%%%%%%%%%%%%%%%%%%%%%%%%%%%%%%%%%%%%%%%%%%%%%%%%%%%%%%%%%%%%%
%%%%%%%%%%%%%%%%%%%%%%%%%%%%%%%%%%%%%%%%%%%%%%%%%%%%%%%%%%%%%%%%%%%%%%%%%
\zerarcounters
\section{Applications}
\label{sec:2}
%%%%%%%%%%%%%%%%%%%%%%%%%%%%%%%%%%%%%%%%%%%%%%%%%%%%%%%%%%%%%%%%%%%%%%%%%
%%%%%%%%%%%%%%%%%%%%%%%%%%%%%%%%%%%%%%%%%%%%%%%%%%%%%%%%%%%%%%%%%%%%%%%%%

%%%%%%%%%%%%%%%%%%%%%%%%%%%%%%%%%%%%%%%%%%%%%%%%%%%%%%%%%%%%%%%%%%%%%%%%%%
\subsection{Non-resonant equations}\label{sub:2.1}
%%%%%%%%%%%%%%%%%%%%%%%%%%%%%%%%%%%%%%%%%%%%%%%%%%%%%%%%%%%%%%%%%%%%%%%%%%

Let us prove that the equations (\ref{eq:1.1}), (\ref{eq:1.4}),
and (\ref{eq:1.5}) -- in particular the NLS, the NLB and
the NLW -- comply with all the Hypotheses and therefore
admit a periodic solution by Theorem \ref{thm:1}.

%%%%%%%%%%%%%%%%%%%%%%%%%%%%%%%%%%%%%%%%%%%%%%%%%%%%%%%%%%%%%%%%%%%%%%%%%%
\subsubsection{The NLS equation}
%%%%%%%%%%%%%%%%%%%%%%%%%%%%%%%%%%%%%%%%%%%%%%%%%%%%%%%%%%%%%%%%%%%%%%%%%%

%%%%%%%%%%%%%%%%%%%%%%%%%%%%%%%%%%%%%%%%%%%%%%%%%%%%%%%%%%%%%%%%%%%%%%%%%
\begin{theorem} \label{thm:2}
Consider the nonlinear Schr\"odinger equation in dimension $D$
$$ \ii \pr_{t} v - \Delta v + \mu \, v = f(x,v,\bar v) , $$
with Dirichlet boundary conditions on the square $[0,\pi]^{D}$,
where $\mu\in(0,\mu_{0})\subset \RRR$ and $f$ is given
according to (\ref{eq:1.2}) and (\ref{eq:1.3}),
with $N=2$, $a_{2,1}=1$ and $a_{r,s} =0$
for $r,s$ such that $r+s=3$ and $(r,s)\neq (2,1)$,
that is $f(x,v,\bar v)=|v|^{2}v+O(|v|^{4})$.
There exist a full measure set $\gotM \subset (0,\mu_{0})$
and a positive constant $\eps_{0}$ such that the following holds.
For all $\mu\in\gotM$ there exists a Cantor set $\gotE(\mu)
\subset [0,\eps_{0}]$, such that for all $\eps\in \gotE(\mu)$
the equation admits a solution $v(x,t)$, which is $2\pi/\om$-periodic
in time and Gevrey-smooth both in time and in space,
and such that
$$ \left| v(x,t) -  \sqrt{\eps} q_{0} {\rm e}^{\ii \om t}
\sin x_{1} \ldots \sin x_{D} \right| \le C \eps,
\qquad \om = D + \mu - \eps , \qquad \left| q_{0} \right| =
\Big( \frac{4}{3} \Big)^{D/2} , $$
uniformly in $(x,t)$. The set $\gotE=\gotE(\mu)$ has positive
Lebesgue measure and satisfies (\ref{eq:1.12}).
\end{theorem}
%%%%%%%%%%%%%%%%%%%%%%%%%%%%%%%%%%%%%%%%%%%%%%%%%%%%%%%%%%%%%%%%%%%%%%%%%

With the notations of Section \ref{sec:1} one has
$\de_{\nn}(\eps) = -\om n + |m|^{2} + \mu$, with $\om=\om_{0} - \eps$
and $\om_{0}=D+\mu$. Then it is easy to check that all items
of Hypothesis \ref{hp:1} are satisfied provided $\mu$ is chosen
in such a way that $|-\om_{0} n + |m|^{2}| \ge \g_{0}|n|^{-\tau_{0}}$.
This is possible for $\mu$ in a full measure set; cf. equation
(2.1) in \cite{GP2}. Then Hypothesis \ref{hp:1} holds
with $c_{0}=c_{2}=c_{3}=1$ and $c_{1}=1/\sqrt{1+4\om_{0}}$.

The subset $\gotQ$ is defined as $\gotQ := \{(n,m)\in \ZZZ^{1+D} :
n=1, \; |m_{i}|=1 \; \forall i=1,\dots D\}$, and one can assume
take $q_{0}$ to be real, so that,
by the Dirichlet boundary conditions, $\gotQ$ is
in fact one-dimensional, and $u_{n,m} = \pm q_{0}$
for all $(n,m)\in \gotQ$. The leading order of the $Q$ equation is
explicitly studied in \cite{GP2}, where it is proved that
Hypothesis \ref{hp:2} is satisfied.

Finally, Hypothesis \ref{hp:3} has been proven by Bourgain \cite{Bo2}
(see also Appendix A6 in \cite{GP2}).

Of course, Theorem \ref{thm:2} refers to solutions with
$m=(1,1,\ldots,1)$, but it easily extends to solutions
which continue other harmonics of the linear equation;
see comments in \cite{GP2}.

Also, the condition on the
nonlinearity can be weakened. In general $N$ can be any
integer $N>1$, and no other conditions must be assumed
on the functions $a_{r,s}(x)$ beyond those mentioned
after (\ref{eq:1.2}). In that case (for simplicity we consider
the same solution of the linear equation as in Theorem \ref{thm:2}),
the leading order of the $Q$ equation becomes $q_{0} =
{\rm sign}(\eps) A_{0} q_{0}^N $ (again by taking for simplicity's
sake $q_{0}$ to be real), where $A_{0}$ is a constant
depending on the nonlinearity. If $A_{0}$ is non-zero,
this surely has a non-trivial non-degenerate solution $q_{0}$
either for positive or negative values of $\eps$.
In general the non-degeneracy condition in item 2 of Hypothesis
\ref{hp:2} has to be verified case by case by computing $A_{0}$.

%%%%%%%%%%%%%%%%%%%%%%%%%%%%%%%%%%%%%%%%%%%%%%%%%%%%%%%%%%%%%%%%%%%%%%%%%%
\subsubsection{The NLW equation}
%%%%%%%%%%%%%%%%%%%%%%%%%%%%%%%%%%%%%%%%%%%%%%%%%%%%%%%%%%%%%%%%%%%%%%%%%%

%%%%%%%%%%%%%%%%%%%%%%%%%%%%%%%%%%%%%%%%%%%%%%%%%%%%%%%%%%%%%%%%%%%%%%%%%
\begin{theorem} \label{thm:3}
Consider the nonlinear wave equation in dimension $D$
$$ \pr_{tt} v - \Delta v + \mu \, v = f(x,v) , $$
with Dirichlet boundary conditions on the square $[0,\pi]^{D}$,
where $\mu\in(0,\mu_{0})\subset\RRR$ and $f$ is given according
to (\ref{eq:1.2}), with $s=0$, $N=2$, $a_{3,0}=1$,
that is $f(x,v)=v^{3}+O(v^{4})$.
There exist a full measure set $\gotM \subset (0,\mu_{0})$
and a positive constant $\eps_{0}$ such that the following holds.
For all $\mu\in\gotM$ there exists a Cantor set $\gotE(\mu)
\subset [0,\eps_{0}]$, such that for all $\eps\in \gotE(\mu)$
the equation admits a solution $v(x,t)$, which is $2\pi/\om$-periodic
in time and Gevrey-smooth both in time and in space,
and such that
$$ \left| v(x,t) - q_{0}
\sqrt{\eps}\cos \om t \sin x_{1} \ldots \sin x_{D} \right|
\le C \eps, \qquad \om = \sqrt{D + \mu - \eps}, \qquad
q_{0} = \left( \frac{4}{3} \right)^{(D+1)/2} ,$$
uniformly in $(x,t)$. The set $\gotE=\gotE(\mu)$ has positive
Lebesgue measure and satisfies (\ref{eq:1.12}).
\end{theorem}
%%%%%%%%%%%%%%%%%%%%%%%%%%%%%%%%%%%%%%%%%%%%%%%%%%%%%%%%%%%%%%%%%%%%%%%%%

In that case one has $\de_{\nn}(\eps) = -\om^{2} n^{2} +
|m|^{2} + \mu $, with $\om^{2}=\om_{0}^{2} - \eps$
and $\om^{2}_{0} = D^{2} +\mu$. Once more, it is easy to check
that Hypothesis \ref{hp:1} is satisfied provided $\mu$ is chosen
in a full measure set, with $c_{0}=c_{2}=c_{3}=1$ and
$c_{1}=1/(1+4\om^{2}_{0})$.

The subset $\gotQ$ is given by $\gotQ := \{(n,m)\in \ZZZ^{1+D}:
n=\pm 1, \; |m_{i}|=1 \; \forall {i}=1,\dots D\}$, and,
if one chooses to look for solutions that are even in time,
then $\gotQ$ is one-dimensional.
The $Q$ equation at $\eps=0$ can be discussed as in the case
of the nonlinear Schr\"odinger equation. For instance
for $f$ as in the statement of Theorem \ref{thm:3}
the non-degeneracy in item 2 of Hypothesis \ref{hp:2}
can be explicitly verified. Again, the
analysis easily extends to more general situations,
under the assumption that the $Q$ equation at $\eps=0$
admits a non-degenerate solution. For a fixed nonlinearity,
this can be easily checked with a simple computation.

Hypothesis \ref{hp:3} has been verified by Bourgain \cite{Bo1},
under some strong conditions on $\om$. Recently the same separation
estimates have been proved by Berti and Bolle \cite{BB3},
by only requiring that $\omega^{2}$ be Diophantine.

%%%%%%%%%%%%%%%%%%%%%%%%%%%%%%%%%%%%%%%%%%%%%%%%%%%%%%%%%%%%%%%%%%%%%%%%%%
\subsubsection{Other equations}
%%%%%%%%%%%%%%%%%%%%%%%%%%%%%%%%%%%%%%%%%%%%%%%%%%%%%%%%%%%%%%%%%%%%%%%%%%

Of course, the separation properties for the NLS equation
imply similar separation also for the nonlinear beam (NLB) equation
$$ \pr_{tt} v + \left( \Delta + \mu \right)^2 \, v = f(x,v) , $$
and in that case we can also consider nonlinearities with
one or two space derivatives.

As in the previous cases one restricts $\mu$ to some full measure set,
and Hypothesis \ref{hp:1} holds with $c_{0}=c_{3}=2$, $c_{2}=1$
and $c_{1}=1/\sqrt{1+2\om_{0}}$. This implies that the
subset $\gotQ$ is one-dimensional, provided we look for
real solutions which are even in time.

The same kind of arguments holds for all equations of the form
(\ref{eq:1.1}) and (\ref{eq:1.4}).
The separation of the points $(m,|m|^{2})$ in $\ZZZ^{D+1}$
implies, by convexity, also the separation of $(m,P(|m|^{2}))$,
with $P(x)$ defined after (\ref{eq:1.1}).

%%%%%%%%%%%%%%%%%%%%%%%%%%%%%%%%%%%%%%%%%%%%%%%%%%%%%%%%%%%%%%%%%%%%%%%%%%
\subsection{Completely resonant equations}\label{sub:2.2}
%%%%%%%%%%%%%%%%%%%%%%%%%%%%%%%%%%%%%%%%%%%%%%%%%%%%%%%%%%%%%%%%%%%%%%%%%%

Here we describe an application to completely resonant NLS and NLB
equations, namely equations (\ref{eq:1.1}) and (\ref{eq:1.4})
with $P(x)=x$ and $\mu=0$, and with Dirichlet boundary conditions
(the case of periodic boundary conditions is easier for fully
resonant equations). 
Since the equation is completely resonant we need some assumption
on the nonlinearity in order to comply with Hypothesis \ref{hp:2}.
We set $f(x,v,\bar v)=|v|^{2}v$ for the NLS and $f(x,v)=v^{3}$
for the NLB (the NLB falls in case (II) and we look for real solutions),
but our proofs extend easily to deal with higher order corrections
which are odd and do not depend explicitly on the space variables.
In the case of the NLS we say that the leading term of the nonlinearity
is cubic and gauge-invariant.\footnote{i.e. the equation up to the
third order is invariant under the transformation $v\to v^{i\al}$
for any $\al\in\RRR$.}

The validity of Hypothesis \ref{hp:1} can be discussed as in the
non-resonant equations of Subsections \ref{sub:2.1}.
The separation properties (Hypothesis \ref{hp:3}) do not change
in the presence of a mass term, and they have been already
discussed in the non-resonant examples of Subsection \ref{sub:2.1}.
Thus, we only need to prove the non-degeneracy of the solution of the
$Q$ equation. Since the nonlinearity does not depend explicitely on
$x$ we look for solutions such that $u_{\nn}\in \RRR$.
  We follow closely \cite{GP2}, but we set $\om_{0}=1$.
This is done for purely notational reasons, and is due to the fact
that a trivial rescaling of time allows us to put $\om_{0}=1$.

%%%%%%%%%%%%%%%%%%%%%%%%%%%%%%%%%%%%%%%%%%%%%%%%%%%%%%%%%%%%%%%%%%%%%%%%%%
\subsubsection{The NLS equation}
%%%%%%%%%%%%%%%%%%%%%%%%%%%%%%%%%%%%%%%%%%%%%%%%%%%%%%%%%%%%%%%%%%%%%%%%%%

The subset $\gotQ$ is infinite-dimensional, i.e.
$\gotQ:=\{ (n,m)\in \NNN \times \ZZZ^{D} : \; n = |m|^{2}\}$.
We set $u_{(n,m)} = q_{m} = a_{m} + O(\eps^{1/2})$ for $(n,m)\in \gotQ$
and restrict our attention to the case $q_{m}\in\RRR$.
At leading order, the $Q$ equation is (cf. \cite{GP2})
\begin{equation}
|m|^{2} a_{m}= \!\!\!\!\!\!\!\! \sum_{\substack{
m_{1},m_{2},m_{3} \\ m_{1}+m_{2}-m_{3}=m \\ 
\langle m_{1}-m_{3},m_{2}-m_{3}\rangle =0}} \!\!\!\!\!\!\!\!
a_{m_{1}} a_{m_{2}} a_{m_{3}} .
\label{eq:2.1} \end{equation}
Note that in the case of \cite{GP2}, the left hand side of
(\ref{eq:2.1}) was $|m|^{2+2s}D^{-1}a_{m}$, with $s$ a free parameter;
then (\ref{eq:2.1}) is recovered by setting $s=0$ and rescaling
by $1/\sqrt{D}$ the coefficients $q_{m}$.

By Lemma 17 of \cite{GP2} -- which holds for all values of $s$ --,
for each $N_{0}\ge 1$ there exist infinitely many finite
sets $\MM_{+} \subset \ZZZ^{D}_{+}$ with $N_{0}$ elements such that
equation (\ref{eq:2.1}) admits the solution (due to the Dirichlet
boundary conditions we describe the solution in $\ZZZ^{D}_{+}$)
\begin{equation}
a_{m} = \begin{cases}
0 , & m\in \ZZZ^{D}_{+} \setminus\MM_{+} \\
{\displaystyle \sqrt{\frac{1}{2^{D+1}-3^{D}} \Big( |m|^{2} - c_{1}
\sum_{m'\in\MM_+}|m'|^{2} \Big) } } , & m\in\MM_{+} ,
\end{cases} \nonumber
\end{equation}
with $c_{1}=2^{D+1}/(2^{D+1}(N_{0}-1)+3^{D})$. 
The set $\MM_{+}$ defines a matrix $J$ on $\ZZZ^{D}$ such that
\begin{equation}
\left( JQ \right)_{m} =  |m|^{2} - 2\!\!\!\!\!\!\!\!\!
\sum_{\substack{m_{1} , m_{2}, m_{3} \\
m_{1}+m_{2}-m_{3}=m \\ \laa m_{1}-m_{3},m_{2}-m_{3}\raa =0}}
\!\!\!\!\!\!\!\!\!\!\!
Q_{m_{1}} a_{m_{2}} a_{m_{3}} - 2\!\!\!\!\!\!\!
\sum_{\substack{m_{1}> m_{2},m_{3} \\ m_{1}+m_{2}-m_{3}=m \\ 
\laa m_{1}-m_{3},m_{2}-m_{3}\raa =0}}\!\!\!\!\!\!\!\!\!\!\!\!\!\!\!
a_{m_{1}} a_{m_{2}} Q_{m_{3}} ,
\label{eq:2.2} \end{equation}
where $m_{1}> m_{2}$ refers, say, to lexicographic ordering of $\ZZZ^{D}$;
see in particular equations (8.5) and (8.7) of \cite{GP2}.

Moreover we know (Lemma 18 of \cite{GP2}) that the matrix $J$
is block-diagonal with blocks of size depending only on $N_{0},D$:
we denote by $K(N_{0},D)$ the bound on such a size.
Whatever the block structure, the matrix $J$ has the form
${\rm diag}(|m|^{2})+ 2T$ where all the entries of $T$ are linear
combinations of terms $q_{m_{i}}q_{m_{j}}$ with integer coefficients.
If we multiply $J$ by $z:=(2^{D+1}-3^{D})(2^{D+1}(N_{0}-1)+3^{D})$
-- which is odd --, we obtain a matrix $J':= {\rm diag}(z|m|^{2})+ 2T'$,
where all the entries of $T'$ are integral linear combinations of the
square roots of a finite number of integers. Let us call the prime
factors of such integers $p_{0}=1,p_{1},p_{2},\ldots$. 

%%%%%%%%%%%%%%%%%%%%%%%%%%%%%%%%%%%%%%%%%%%%%%%%%%%%%%%%%%%%%%%%%%%%%%%%%%
\begin{defi} \label{def:3}
\textbf{ (The lattice $\boldsymbol \ZZZ_{\boldsymbol 1}^{\boldsymbol D}$).}
Let $\ZZZ^{D}_{1} := (1,0,\dots,0)+2\ZZZ^{D}$ be the affine
lattice of integer vectors such that the first component is odd and
the others even. Let $\ZZZ^{D}_{1,+}$ be its intersection
with $\ZZZ_{+}^{D}$. Of course, for all $m\in \ZZZ^{D}_{1}$
one has $|m|^{2}$ odd.
\end{defi}
%%%%%%%%%%%%%%%%%%%%%%%%%%%%%%%%%%%%%%%%%%%%%%%%%%%%%%%%%%%%%%%%%%%%%%%%%%

Since we are working with odd nonlinearities which do not depend
explicitly on the space variables we look for
solutions such that $u_{n,m}=0$ if $m\notin \ZZZ^{D}_{1}$.
\vskip5pt
Let $1,p_{1},\dots,p_{k}$ be prime numbers (as above), and let
$a_{1},\dots,a_{K}$ be the set of all products of square roots of
different numbers $p_{i}$, i.e. $a_{1}=1$, $a_{2}=\sqrt{p_{1}}$,
$a_{3}=\sqrt{p_{1}p_{2}}$, etc.
It is clear that the set of integral linear combinations of $a_{i}$ is a
ring (of algebraic integers). We denote it by $\gota$.
The following Lemma is a simple consequence of Galois theory \cite{Ar}.
For completeness, the proof is given in Appendix \ref{app:A}.

%%%%%%%%%%%%%%%%%%%%%%%%%%%%%%%%%%%%%%%%%%%%%%%%%%%%%%%%%%%%%%%%%%%%%%%%%%
\begin{lemma} \label{lem:1}
 The numbers $a_{i}$ are linearly independent over the rationals.
\end{lemma}
%%%%%%%%%%%%%%%%%%%%%%%%%%%%%%%%%%%%%%%%%%%%%%%%%%%%%%%%%%%%%%%%%%%%%%%%%%

Immediately we have the following corollary ($I$ denotes the identity).

%%%%%%%%%%%%%%%%%%%%%%%%%%%%%%%%%%%%%%%%%%%%%%%%%%%%%%%%%%%%%%%%%%%%%%%%%%
\begin{coro} \label{cor:1}
In $\gota$ consider  $2\gota$, i.e. the set of linear
combinations with even coefficients.
\begin{itemize}
\item $2\gota$ is a proper ideal,
and the quotient ring $\gota/2\gota$ is thus a non-zero ring.
\item if a matrix $M$ with entries in $\gota$ is such that $M-I$
has all entries in $2\gota$, then $M$ is invertible.
\end{itemize}
\end{coro}
%%%%%%%%%%%%%%%%%%%%%%%%%%%%%%%%%%%%%%%%%%%%%%%%%%%%%%%%%%%%%%%%%%%%%%%%%%

The point of Corollary \ref{cor:1} is that the determinant of
$M=I+2M'$, with the entries of $M'$ in $\gota$, is $1+2\alpha$,
with $\alpha\in\gota$.
Hence, by Lemma \ref{lem:1}, $2\alpha \neq \pm 1$.

%%%%%%%%%%%%%%%%%%%%%%%%%%%%%%%%%%%%%%%%%%%%%%%%%%%%%%%%%%%%%%%%%%%%%%%%%%
\begin{lemma} \label{lem:2}
For all $N_{0}$ and for all $\MM_{+}\subset \ZZZ^{D}_{1,+}$ the matrix
$J$  defined by $\MM_{+}$ is invertible. Its inverse is a block matrix
with blocks of dimension depending only on $N_{0},D$ so that for some
appropriate $C$ one has $(J^{-1})_{m,m'} \le C$ if $|m-m'|\le K(N_{0},D)$,
while $(J^{-1})_{m,m'}=0$ otherwise.
\end{lemma}
%%%%%%%%%%%%%%%%%%%%%%%%%%%%%%%%%%%%%%%%%%%%%%%%%%%%%%%%%%%%%%%%%%%%%%%%%%

%%%%%%%%%%%%%%%%%%%%%%%%%%%%%%%%%%%%%%%%%%%%%%%%%%%%%%%%%%%%%%%%%%%%%%%%%%
\prova We use Corollary \ref{lem:1}, the fact that the matrix $J'$
has entries in $\gota$ and the fact that $z|m|^{2}$ is odd for
all $m\in \ZZZ^{D}_{1,+}$. \EP
%%%%%%%%%%%%%%%%%%%%%%%%%%%%%%%%%%%%%%%%%%%%%%%%%%%%%%%%%%%%%%%%%%%%%%%%%%

Now, we can state our result on the completely resonant NLS.

%%%%%%%%%%%%%%%%%%%%%%%%%%%%%%%%%%%%%%%%%%%%%%%%%%%%%%%%%%%%%%%%%%%%%%%%%
\begin{theorem} \label{thm:4}
Consider the nonlinear Schr\"odinger equation in dimension $D$
$$ \ii \pr_{t} v - \Delta v = f(v,\bar v) , $$
with Dirichlet boundary conditions on the square $[0,\pi]^{D}$,
where $f$ is given according to (\ref{eq:1.2})
and (\ref{eq:1.3}), with $N=2$, $a_{2,1}=1$, $a_{r,s}=0$
for $r,s$ such that $r+s=3$ and $(r,s)\neq (2,1)$,
and $a_{r,s}(x)$ independent of $x$ for $r+s>3$
(so that in particular $a_{r,s}=0$ for even $r+s$).
For any $N_{0} \ge 1$ there exist sets $\MM_{+}$
of $N_{0}$ vectors in $\ZZZ_{+}^{D}$ and real amplitudes
$\{a_{m}\}_{m\in\MM_{+}}$ such that the following holds.
There exist a positive constant $\eps_{0}$ and a Cantor set $\gotE
\subset [0,\eps_{0}]$, such that for all $\eps\in \gotE$
the equation admits a solution $v(x,t)$, which is $2\pi/\om$-periodic
in time and Gevrey-smooth both in time and in space, and such that,
setting
\begin{equation}
q_{0}(x,t) = (2\ii)^D\sum_{m\in\MM_{+}}
a_{m} {\rm e}^{\ii |m|^{2} t}
\sin m_{1} x_{1} \ldots \sin m_{D} x_{D} , \qquad
\qquad \om = 1 - \eps ,
\label{eq:2.3} \end{equation}
one has
\begin{equation}
\left| v(x,t) - \sqrt{\eps} q_{0}(x,\om t)  \right| \le C \eps, \nonumber
\end{equation}
uniformly in $(x,t)$. The set $\gotE$ has positive
Lebesgue measure and satisfies (\ref{eq:1.12}).
\end{theorem}
%%%%%%%%%%%%%%%%%%%%%%%%%%%%%%%%%%%%%%%%%%%%%%%%%%%%%%%%%%%%%%%%%%%%%%%%%

%%%%%%%%%%%%%%%%%%%%%%%%%%%%%%%%%%%%%%%%%%%%%%%%%%%%%%%%%%%%%%%%%%%%%%%%%%
\subsubsection{The beam equation}
%%%%%%%%%%%%%%%%%%%%%%%%%%%%%%%%%%%%%%%%%%%%%%%%%%%%%%%%%%%%%%%%%%%%%%%%%%

We set $\om^{2}=\om_{0}^{2}-\eps=1-\eps$ (recall that we are
assuming $\om_{0}=1$ by a suitable time rescaling).
The subset $\gotQ$ is given by
$ \gotQ:=\{ (n,m)\in \NNN\times \ZZZ^{D}: \; |n| = |m|^{2}\}$.
We set $u_{n,m}=q^{+}_{m}$ for $n=|m|^{2}$ and $u_{n,m}=q^{-}_{m}$
for $n=-|m|^{2}$. We can require that $q^{+}_{m}=q^{-}_{m}\equiv q_{m}$
for all $m$ (we obtain a solution which is even in time).
Since we look  for real solutions, this implies that
$q_{m}\in\RRR$ if $D$ is even and $q_{m} \in \ii\RRR$
if $D$ is odd. Since the nonlinearity
does not depend explicitly on $x$, we can look for solutions
$u_{n,m}$ such that $m\in \ZZZ^{D}_{1}$ (see Definition \ref{def:3}).

Finally the separation properties of the small divisors do not
depend on the presence of the mass term, so that we only need
to prove the existence and non-degeneracy of the solutions
of the bifurcation equation.

The $Q$ equation at leading order is
$$ |m|^4 a_{m} = (-1)^{D} \!\!\!\!\!\!\!\!\!\!\!\!\!
\sum_{\substack{m_{1}+m_{2}+m_{3}=m \\
\pm |m_{1}|^{2} \pm |m_{2}|^{2} \pm |m_{3}|^{2}= \pm |m|^{2}}}
a_{m_{1}} a_{m_{2}} a_{m_{3}} , $$
where we have set $|q_{m} |= a_{m} + O(\eps^{1/2})$.

%%%%%%%%%%%%%%%%%%%%%%%%%%%%%%%%%%%%%%%%%%%%%%%%%%%%%%%%%%%%%%%%%%%%%%%%%%
\begin{lemma} \label{lem:3}
The condition $\pm |m_{1}|^{2} \pm |m_{2}|^{2} \pm |m_{3}|^{2} =
\pm |m|^{2}$, for $m_{i},m\in \ZZZ_{1}^{D}$, is equivalent to
$\laa m_{1} + m_{3}, m_{2}+m_{3}\raa=0$.
\end{lemma}
%%%%%%%%%%%%%%%%%%%%%%%%%%%%%%%%%%%%%%%%%%%%%%%%%%%%%%%%%%%%%%%%%%%%%%%%%%

%%%%%%%%%%%%%%%%%%%%%%%%%%%%%%%%%%%%%%%%%%%%%%%%%%%%%%%%%%%%%%%%%%%%%%%%%%
\prova The condition $|m_{1}|^{2} + |m_{2}|^{2} + |m_{3}|^{2}=
(m_{1}+m_{2}+m_{3})^{2}$ is equivalent to
$\laa m_{1},m_{2}+m_{3}\raa +\laa m_{2},m_{3}\raa=0$,
which is impossible since the left hand side is an odd integer.
The same happens with the condition $|m_{1}|^{2} - |m_{2}|^{2} -
|m_{3}|^{2}=(m_{1}+m_{2}+m_{3})^{2}$. Thus, we are left with $|m_{1}|^{2}
+ |m_{2}|^{2} - |m_{3}|^{2}=(m_{1}+m_{2}+m_{3})^{2}$,
which implies $\laa m_{1}+m_{3},m_{2}+m_{3}\raa=0$. \EP
%%%%%%%%%%%%%%%%%%%%%%%%%%%%%%%%%%%%%%%%%%%%%%%%%%%%%%%%%%%%%%%%%%%%%%%%%%

Lemma \ref{lem:3} implies that the bifurcation equation,
restricted to $\ZZZ^{D}_1$, is identical to that of a smoothing
NLS with $s=2$; cf. \cite{GP2}. Indeed by recalling that
$q_{m}= (-1)^{D} q_{-m}$ one has
\begin{equation}\label{b1}
|m|^{4} a_m = \sum_{\substack{m_{1}+m_{2}-m_{3}=m \\
\laa m_{1}-m_{3},m_{2}-m_{3}\raa = 0}}
a_{m_{1}}a_{m_{2}}a_{m_{3}} . 
\end{equation}
Then we can repeat the arguments of the previous subsection.
By Lemma 17 of \cite{GP2} -- which holds for all values of $s$ --
for each $N_{0}\ge 1$ there exist infinitely many finite
sets $\MM_{+} \subset \ZZZ^{D}_{1,+}$ with $N_{0}$ elements such that
the equation (\ref{b1}) has the solution
\begin{equation}
a_{m} = \begin{cases}
0 , & m\in \ZZZ^{D}_{+} \setminus\MM_{+} \\
{\displaystyle \sqrt{ \frac{1}{2^{D+1}-3^{D}} \Big( |m|^{4} - c_{1}
\sum_{m'\in\MM_+}|m'|^{4} \Big) } } , & m\in\MM_{+} ,
\end{cases} \nonumber
\end{equation}
with $c_{1}=2^{D+1}/(2^{D+1}(N_{0}-1)+3^{D})$.

The matrix $J$ is defined as in (\ref{eq:2.2}), only with $|m|^{4}$
on the diagonal. We know (Lemma 18 of \cite{GP2} does not depend
on the values of $s$) that the matrix $J$ is block-diagonal
with blocks of size bounded by $K(N_{0},D)$ (defined as in
subsection 2.2.1). Whatever the block structure, the matrix $J$
has the form ${\rm diag}(|m|^{4})+ 2T$, where all the entries of $T$
are linear combinations of terms $a_{m_{i}}a_{m_{j}}$
with integer coefficients. If we multiply $J$ by
$z:=(2^{D+1}-3^{D})(2^{D+1}(N_{0}-1)+3^{D})$ -- which is odd --, we
obtain a matrix $J':= {\rm diag}(z|m|^{4})+ 2T'$, where all the
entries of $T'$ are linear combinations of the square roots of a
finite number of integers; finally $z|m|^4$ is clearly odd and we can
apply Lemma \ref{lem:1} to obtain the analogous of Lemma \ref{lem:2}.
Thus, a theorem analogous to Theorem \ref{thm:4} is obtained,
with $q_{0}(x,t)$ in (\ref{eq:2.3}) replaced with
\begin{equation}
q_{0}(x,t) =2^{D+1} \sum_{m\in\MM_{+}}
a_{m} \cos |m|^{2} t \,
\sin m_{1} x_{1} \ldots \sin m_{D} x_{D} , \qquad
\qquad \om^{2} = 1 - \eps . \nonumber
\end{equation}
We leave the formulation to the reader.

%%%%%%%%%%%%%%%%%%%%%%%%%%%%%%%%%%%%%%%%%%%%%%%%%%%%%%%%%%%%%%%%%%%%%%%%%
%%%%%%%%%%%%%%%%%%%%%%%%%%%%%%%%%%%%%%%%%%%%%%%%%%%%%%%%%%%%%%%%%%%%%%%%%
\zerarcounters
\section{Technical set-up and propositions}
\label{sec:3}
%%%%%%%%%%%%%%%%%%%%%%%%%%%%%%%%%%%%%%%%%%%%%%%%%%%%%%%%%%%%%%%%%%%%%%%%%
%%%%%%%%%%%%%%%%%%%%%%%%%%%%%%%%%%%%%%%%%%%%%%%%%%%%%%%%%%%%%%%%%%%%%%%%%

%%%%%%%%%%%%%%%%%%%%%%%%%%%%%%%%%%%%%%%%%%%%%%%%%%%%%%%%%%%%%%%%%%%%%%%%%
\subsection{Renormalised $\boldsymbol P$-$\boldsymbol Q$ equations}
%%%%%%%%%%%%%%%%%%%%%%%%%%%%%%%%%%%%%%%%%%%%%%%%%%%%%%%%%%%%%%%%%%%%%%%%%

Group the equations (\ref{eq:1.10}) for $\nn\in\gotO$
as a matrix equation. Setting
\begin{equation}
\;\;\;\;\;U = \{u^{\s}_{\nn}\}_{\nn\in \gotO}^{\s=\pm} ,
\quad V = \{u^{\s}_{\nn}\}_{\nn\in \gotR}^{\s=\pm},
\quad Q = \{u^{\s}_{\nn}\}_{\nn\in \gotQ}^{\s=\pm},
\quad F = \{f^{\s}_{\nn}\}_{\nn\in \gotO}^{\s=\pm} ,
\quad \DD(\eps) = {\rm diag} \left\{\de_{\nn}(\eps)
\right\}_{\nn\in \gotO}^{\s=\pm} ,
\label{eq:3.1} \end{equation}
the $P$ equations spell
\begin{equation}
\begin{cases} \DD(\eps) \,  U = \eps F(U,V,Q,\eps^{1/N}) , & \\
u^{\s}_{\nn} = \eps \de_{\nn}^{-1}(\eps) \,
f^{\s}_{\nn}(U,V,Q,\eps^{1/N}) , & \nn\in\gotR , \end{cases}
\label{eq:3.2} \end{equation}
with a reordering of the arguments of the coefficients $f^{\s}_{\nn}$.

We want to introduce an appropriate ``correction'' to the left hand side
of (\ref{eq:3.2}). We shall consider self-adjoint matrices $\widehat
M(\eps):=\{\widehat M_{\nn,\nn'}^{\s,\s'}(\eps)\}_{\nn,\nn'\in
\gotO}^{\s,\s'=\pm}$, which for each fixed $\eps$ are
block-diagonal on the sets $\Delta_{j}(\eps)$ (cf. Definition \ref{def:2}),
namely $\widehat M^{\s,\s'}_{\nn,\nn'}(\eps)\neq 0$ can hold
only if $\nn,\nn'\in \Delta_{j}(\eps)$ for some $j$. Moreover we
require for $\widehat M_{\nn,\nn'}^{\s,\s'}(\eps)$ to
depend smoothly on $\eps$, at least in a large measure set.

We shall first introduce the self-adjoint matrices $\widehat M$ as
independent parameters, and eventually we shall manage to fix them
as functions of the parameter $\eps$. Note that in order to have
$u^{+}_{\nn}=\overline{u^{-}_{\nn}}$ we must require that
$\widehat M_{\nn,\nn'}^{\s,\s'}= \widehat M_{\nn',\nn}^{-\s',-\s} $.

%%%%%%%%%%%%%%%%%%%%%%%%%%%%%%%%%%%%%%%%%%%%%%%%%%%%%%%%%%%%%%%%%%%%%%%%%
\begin{defi} \label{def:4}
\textbf{ (The set $\boldsymbol \gotG$ and the matrix
$\widehat{\boldsymbol\chi}_{\boldsymbol 1}$).} 
Call $\gotG=\{ 1/4> \bar\g>0 : \left| |\de_{\nn}(0)|-\bar\g \right|
\ge \bar \g_{0}/|\nn|^{\bar\tau_{0}} \text{ for all }
\nn\in\ZZZ^{D+1}_{*}\}$, for suitable
constants $\bar\g_{0},\bar\tau_{0}>0$. For $\bar\g\in
\gotG$, we introduce the step function $\bar \chi_{1}(x)$ such that
$\bar\chi_{1}(x)=0$ if $|x| \geq\bar\g$ and $\bar\chi_{1}(x)=1$
if $|x|<\bar\g$, and set $\bar\chi_{0}(x)=1-\bar\chi_{1}(x)$.
We then introduce the ($\eps$-dependent) diagonal matrices
$\widehat\chi_{1}={\rm diag}\{\bar\chi_{1}(\de_{\nn}(\eps))
\}_{\nn\in\gotO}^{\s=\pm}$ and
$\widehat\chi_{0}={\rm diag}\{\bar\chi_{0}(\de_{\nn}(\eps))
\}_{\nn\in\gotO}^{\s=\pm}$.
\end{defi}
%%%%%%%%%%%%%%%%%%%%%%%%%%%%%%%%%%%%%%%%%%%%%%%%%%%%%%%%%%%%%%%%%%%%%%%%%

%%%%%%%%%%%%%%%%%%%%%%%%%%%%%%%%%%%%%%%%%%%%%%%%%%%%%%%%%%%%%%%%%%%%%%%%%%
\begin{rmk} \label{rmk:8}
One has $\gotG\neq\emptyset$. Moreover,
for any interval $\gotU\subset(0,1/4)$, the relative measure
of the set $\gotU\cap\gotG$ tends to 1 as $\bar\gamma_{0}$ tends to 0,
provided $\bar\tau_{0}$ is large enough
\end{rmk}
%%%%%%%%%%%%%%%%%%%%%%%%%%%%%%%%%%%%%%%%%%%%%%%%%%%%%%%%%%%%%%%%%%%%%%%%%%

%%%%%%%%%%%%%%%%%%%%%%%%%%%%%%%%%%%%%%%%%%%%%%%%%%%%%%%%%%%%%%%%%%%%%%%%%%
\begin{rmk} \label{rmk:9}
Note that
$\widehat\chi_{1}^{2}=\widehat\chi_{1}$ and
$\widehat\chi_{1}\widehat\chi_{0}=0$,
with $0$ the null matrix.
\end{rmk}
%%%%%%%%%%%%%%%%%%%%%%%%%%%%%%%%%%%%%%%%%%%%%%%%%%%%%%%%%%%%%%%%%%%%%%%%%%

%%%%%%%%%%%%%%%%%%%%%%%%%%%%%%%%%%%%%%%%%%%%%%%%%%%%%%%%%%%%%%%%%%%%%%%%%%
\begin{defi} \label{def:5}
\textbf{ (Resonant sets).}
A set $\NN=\{\nn_{1},\dots,\nn_{m}\}\subset \gotO$ is {\rm resonant}
if there exists $\eps\in\gotE_{0}$ and $j\in\NNN$ such that
$\nn_{1},\dots,\nn_{m}\in\Delta_{j}(\eps)$.
A resonant set $\{\nn_{1},\nn_{2}\}$ with $m=2$
will be called a {\rm resonant pair}. Given a resonant set
$\NN=\{\nn_{1},\dots,\nn_{m}\}$ we call $\CC_{\NN}$ the set of all
$\nn\in\gotO$ such that $\NN\cup\{\nn\}$ is still a resonant set.
Finally set $\overline\CC_{\NN}(\eps)
:=\{\nn'\in\CC_{\NN} : |\delta_{\nn'}(\eps)| < \bar\g\}$.
\end{defi}
%%%%%%%%%%%%%%%%%%%%%%%%%%%%%%%%%%%%%%%%%%%%%%%%%%%%%%%%%%%%%%%%%%%%%%%%%%

Define the {\it renormalised $P$ equation} as
\begin{equation}
\begin{cases}
\left( \DD(\eps) + \widehat M \right) U =
\h^{N} F(U,V,Q,\eta) + L \, U , & \\
u^{\s}_{\nn} = \eta^{N} \, \de_{\nn}^{-1}(\eps)\,
f^{\s}_{\nn}(U,V,Q,\eta) , &
\nu \in \gotR , \end{cases}
\label{eq:3.3}  \end{equation}
with $\widehat M = \widehat\chi_{1} M\widehat\chi_{1}$,
where $\h$ is a real parameter, while
$M=\{M^{\s,\s'}_{\nn,\nn'}\}_{\nn,\nn'\in \gotO}^{\s,\s'=\pm}$ and
$L=\{L^{\s,\s'}_{\nn,\nn'}\}_{\nn,\nn'\in \gotO}^{\s,\s'=\pm}$
are  self-adjoint matrices of free parameters with the properties:
\begin{enumerate}
\item $M^{\s,\s'}_{\nn,\nn'}=L^{\s,\s'}_{\nn,\nn'}=0$ if $\{\nn,\nn'\}$
is not a  resonant pair.
\item $M_{\nn,\nn'}^{\s,\s'}=M_{\nn',\nn}^{-\s',-\s}$ and
$L_{\nn,\nn'}^{\s,\s'} = L_{\nn',\nn}^{-\s',-\s} $.
\end{enumerate}

The {\it renormalised $Q$ equation} is defined as
\begin{equation}
\,u^{\s}_{\nn} = \sum_{\nn\in\gotQ} \sum_{\s'=\pm}
(J^{-1})^{\s,\s'}_{\nn,\nn'} f^{\s'}_{\nn'}(U,V,Q,\eta) ,
\qquad \nn \in \gotQ .
\label{eq:3.4}
\end{equation}
The parameter $\h$ and the {\it counterterms} $L$ will have to
satisfy eventually the identities ({\it compatibility equation})
\begin{equation} \h = \eps^{1/N} , \qquad
\widehat M = L .
\label{eq:3.5} \end{equation}

We proceed in the following way: first we solve the renormalised
$P$ and $Q$ equations (\ref{eq:3.3}) and (\ref{eq:3.4}),
then we impose the compatibility equation (\ref{eq:3.5}).

%%%%%%%%%%%%%%%%%%%%%%%%%%%%%%%%%%%%%%%%%%%%%%%%%%%%%%%%%%%%%%%%%%%%%%%%%
\subsection{Matrix spaces}
%%%%%%%%%%%%%%%%%%%%%%%%%%%%%%%%%%%%%%%%%%%%%%%%%%%%%%%%%%%%%%%%%%%%%%%%%%

Here we introduce some notations and properties that we shall need
in the following.

%%%%%%%%%%%%%%%%%%%%%%%%%%%%%%%%%%%%%%%%%%%%%%%%%%%%%%%%%%%%%%%%%%%%%%%%%%
\begin{defi} \label{def:6}
\textbf{ (The Banach space $\boldsymbol\BB_{\boldsymbol\kappa}$).}
We consider the space of infinite-dimensional self-adjoint matrices
$\{M^{\s,\s'}_{\nn,\nn'}\}^{\s,\s'=\pm}_{\nn,\nn'\in\gotO}$ such that
$M^{\s,\s'}_{\nn,\nn'}=0$ if $\{\nn,\nn'\}$ is not resonant.
For $\rho,\kappa>0$ we equip such a space with the norm
\begin{equation}
\left| M \right|_{\kappa} := \sup_{\nn,\nn'\in\gotO} \sup_{\s,\s'=\pm}
\left| M_{\nn,\nn'}^{\s,\s'} \right|
{\rm e}^{\kappa|\nn-\nn'|^{\rho}} , \nonumber
\end{equation}
so obtaining a Banach space that we call $\BB_{\kappa}$.
For $L$ a linear operator on $\BB_{\kappa}$ define the operator norm
$$ |L|_{\rm op} = \sup_{M\in\BB_{\kappa}}
\frac{|LM|_{\kappa}}{|M|_{\kappa}}. $$
\end{defi}
%%%%%%%%%%%%%%%%%%%%%%%%%%%%%%%%%%%%%%%%%%%%%%%%%%%%%%%%%%%%%%%%%%%%%%%%%%

%%%%%%%%%%%%%%%%%%%%%%%%%%%%%%%%%%%%%%%%%%%%%%%%%%%%%%%%%%%%%%%%%%%%%%%%%%
\begin{defi} \label{def:7}
\textbf{ (Matrix norms).}
Let $A$ be a $d\times d$ self-adjoint  matrix, and denote with $A(a,b)$
and $\la^{(a)}(A)$ its entries and its eigenvalues, respectively.
We define the norms
\begin{equation}
\left| A \right|_{\io} := \max_{1 \leq a,b\leq d} |A(a,b)| ,
\qquad \Vert A \Vert := \frac{1}{\sqrt{d}}
\sqrt{ {\rm tr}(A^2) } , \qquad
\left\| A \right\|_{2} := \max_{|x|_{2} \le 1} \left| Ax \right|_{2} ,
\nonumber
\end{equation}
where, given a vector $x\in\RRR^{d}$, we denote by $|x|_{2}$
its Euclidean norm.
\end{defi}
%%%%%%%%%%%%%%%%%%%%%%%%%%%%%%%%%%%%%%%%%%%%%%%%%%%%%%%%%%%%%%%%%%%%%%%%%%

%%%%%%%%%%%%%%%%%%%%%%%%%%%%%%%%%%%%%%%%%%%%%%%%%%%%%%%%%%%%%%%%%%%%%%%%%%
\begin{lemma} \label{lem:4}
Given $d\times d$ self-adjoint matrix $A$,
the following properties hold.
\begin{enumerate}
\item The norm $\Vert A\Vert$  depends smoothly on
the coefficients $A(a,b)$.
\item One has $\Vert A\Vert /\sqrt{d}\le |A|_\io\le
\sqrt{d}\Vert A\Vert $.
\item One has $\max_{1\le a \le d} |\la^{(a)}(A)| /\sqrt{d} \le
\Vert A\Vert \le \max_{1\le a \le d} |\la^{(a)}(A)|$.
\item For invertible $A$ one has
$\pr_{A(a,b)} A^{-1}(a',b')= - A^{-1}(a',a) \, A^{-1}(b,b')$ and
$\pr_{A(a,b)}\Vert A\Vert = A(a,b)/d \Vert A\Vert$.
\end{enumerate}
\end{lemma}
%%%%%%%%%%%%%%%%%%%%%%%%%%%%%%%%%%%%%%%%%%%%%%%%%%%%%%%%%%%%%%%%%%%%%%%%%%

Here and henceforth we shall write $A=\DD(\eps)+\widehat M$
in (\ref{eq:3.3}).

%%%%%%%%%%%%%%%%%%%%%%%%%%%%%%%%%%%%%%%%%%%%%%%%%%%%%%%%%%%%%%%%%%%%%%%%%%
\begin{defi} \label{def:8}
\textbf{ (Small divisors).}
For $\nn\in\gotO$ define $A^{\nn}(\eps)$ as the matrix with entries
$\bar\chi_{1}(\delta_{\nn}(\eps))\,A^{\s_{1},\s_{2}}_{\nn_{1},\nn_{2}}$
such that $\nn_{1},\nn_{2}\in \overline\CC_{\nn}(\eps)$
and $\s_{1},\s_{2}=\pm$. If $|\delta_{\nn}(\eps)| < \bar\g$,
define also $d^{\nn}(\eps) := 2| \overline\CC_{\nn}(\eps)|$ and
$p_{\nn}(\eps) = \min\{|\nn'| : \nn'\in \overline\CC_{\nn}(\eps)\}$.
 For real positive $\xi$, define the {\rm small divisor}
\begin{equation}
x_{\nn}(\eps) := \frac{1}{p_{\nn}^{\xi}(\eps)} \left\|
(A^{\nn}(\eps))^{-1} \right\|^{-1} , \nonumber
\end{equation}
if $A$ is invertible, and set $x_{\nn}(\eps)=0$ if $A$ is not invertible.
\end{defi}
%%%%%%%%%%%%%%%%%%%%%%%%%%%%%%%%%%%%%%%%%%%%%%%%%%%%%%%%%%%%%%%%%%%%%%%%%%

%%%%%%%%%%%%%%%%%%%%%%%%%%%%%%%%%%%%%%%%%%%%%%%%%%%%%%%%%%%%%%%%%%%%%%%%%%
\begin{rmk} \label{rmk:10}
Note that for $\nn\in\Delta_{j}(\eps)$ one has
$p_{\nn}(\eps)=p_{j}(\eps)$,
$d_{\nn}(\eps) \le 2 |\Delta_{j}(\eps)|$, and
$A^{\nn}(\eps)=A^{\nn'}(\eps)$ for all $\nn'\in\overline\CC_{\nn}(\eps)$.
This shows that $d_{\nn}(\eps)$, $x_{\nn}(\eps)$ and $p_{\nn}(\eps)$
are the same for all $\nn'\in\overline\CC_{\nn}(\eps)$.
Note also that, if $\nn\in\Delta_{j}(\eps)$ for some $j\in\NNN$, then
one has $\overline\CC_{\nn}(\eps) = \{\nn'\in\Delta_{j}(\eps) :
|\delta_{\nn'}(\eps)| < \bar\g\}$. Hypothesis \ref{hp:3}
implies $d_{\nn}(\eps) \le 2 C_{1} p_{\nn}^{\al}(\eps)$.
\end{rmk}
%%%%%%%%%%%%%%%%%%%%%%%%%%%%%%%%%%%%%%%%%%%%%%%%%%%%%%%%%%%%%%%%%%%%%%%%%%

%%%%%%%%%%%%%%%%%%%%%%%%%%%%%%%%%%%%%%%%%%%%%%%%%%%%%%%%%%%%%%%%%%%%%%%%%%
\begin{defi} \label{def:9}
\textbf{ (The sets $\boldsymbol\DDD_{\boldsymbol 0}$,
$\boldsymbol\DDD_{\boldsymbol 1} \boldsymbol(
\boldsymbol\g\boldsymbol)$, $\boldsymbol\DDD_{\boldsymbol 2}
\boldsymbol(\boldsymbol\g\boldsymbol)$, and
$\boldsymbol\DDD\boldsymbol(\boldsymbol\g\boldsymbol)$).}
We define $\DDD_{0} =\{(\eps,M) : \eps\in\gotE_{0} , \;
|M|_\kappa\le C_{0} \eps_{0} \}$, for a suitable positive constant $C_{0}$,
and, for fixed $\tau,\tau_{1}>0$ and $\g<\bar\g$, we set $\DDD_{1}(\g)
= \{ (\eps,M) \in \DDD_{0} : x_{\nn} \ge \g /p_{\nn}^{\tau}(\eps)
\text{ for all } j\in \NNN\}$, $\DDD_{2}(\g) =
\{ (\eps,M) \in \DDD_{0} : ||\delta_{\nn}(\eps)|-\bar\g| \ge
\g / |\nn|^{\tau_{1}} \text{ for all } \nn \in \gotO \}$,
and $\DDD(\g)=\DDD_{1}(\g) \cap \DDD_{2}(\g)$.
\end{defi}
%%%%%%%%%%%%%%%%%%%%%%%%%%%%%%%%%%%%%%%%%%%%%%%%%%%%%%%%%%%%%%%%%%%%%%%%%%

%%%%%%%%%%%%%%%%%%%%%%%%%%%%%%%%%%%%%%%%%%%%%%%%%%%%%%%%%%%%%%%%%%%%%%%%%%
\begin{defi} \label{def:10}
\textbf{ (The sets $\boldsymbol\II_{\boldsymbol \NN} \boldsymbol(
\boldsymbol\g\boldsymbol)$ and
$\overline{\boldsymbol\II}_{\boldsymbol \NN} \boldsymbol(
\boldsymbol\g\boldsymbol)$).}
Given a resonant set $\NN$ we define $\overline\II_{\NN}(\g) :=
\{ \eps \in \gotE_{0} : \exists \nn\in \CC_{\NN} \text{ such that }
||\delta_{\nn} (\eps)|-\bar\g| < \g |\nn|^{-\tau_{1}} \;\}$,
and set $\II_{\NN}(\g) :=\{ (\eps,M) \in \DDD_{0} : \eps \in
\overline\II_{\NN}(\g)\}$.
\end{defi}
%%%%%%%%%%%%%%%%%%%%%%%%%%%%%%%%%%%%%%%%%%%%%%%%%%%%%%%%%%%%%%%%%%%%%%%%%%

%%%%%%%%%%%%%%%%%%%%%%%%%%%%%%%%%%%%%%%%%%%%%%%%%%%%%%%%%%%%%%%%%%%%%%%%%%
\subsection{Main propositions}
%%%%%%%%%%%%%%%%%%%%%%%%%%%%%%%%%%%%%%%%%%%%%%%%%%%%%%%%%%%%%%%%%%%%%%%%%%

We state the propositions which represent our main technical
results. Theorem \ref{thm:1} is an immediate consequence of
Propositions \ref{prop:1} and \ref{prop:2} below.

%%%%%%%%%%%%%%%%%%%%%%%%%%%%%%%%%%%%%%%%%%%%%%%%%%%%%%%%%%%%%%%%%%%%%%%%%%
\begin{prop} \label{prop:1}
There exist positive constants
$K_{0},K_{1},\kappa,\rho,\h_{0}$ such that the following
holds true. For $(\eps,M)\in \DDD(\g)$, there exists a matrix
$L(\h,\eps,M)\in \BB_\kappa$, such that the following holds.
\begin{enumerate}
\item For each $\eps$ the matrix $L(\h,\eps,M)$ is block-diagonal
so as to satisfy $L(\h,\eps,M) = \widehat\chi_{1}
L(\h,\eps,M) \widehat\chi_{1}$.
\item There exists a unique solution $u^{\s}_{\nn}(\h,M,\eps)$,
with $\nn \in \ZZZ^{D+1}$, of equations (\ref{eq:3.3}) and
(\ref{eq:3.4}), which is analytic in $\h$ for $|\h|\leq \h_{0}$,
and such that for all $\nn\in\ZZZ^{D+1}$ and $\s = \pm$
\begin{equation}
\left| u^{\s}_{\nn}(\h,M,\eps) \right|
\le |\h| \, K_{0} {\rm e}^{- \kappa |\nn|^{1/2}} . \nonumber
\end{equation}
\item The matrix elements $L^{\s,\s'}_{\nn,\nn'}(\h,\eps,M)$
are analytic in $\h$ for $|\h|\leq \h_{0}$,
and uniformly bounded for $(\eps,M)\in \DDD(\g)$ as
\begin{equation}
\left| L(\h,\eps,M)\right|_{\kappa} \leq |\h|^N \, K_{0} . \nonumber
\end{equation}
\item The functions $u^{\s}_{\nn}(\h,\eps,M)$ can be extended on the
set $\DDD_{0}$ to $C^1$ functions $u^{E \, \s}_{\nn}(\h,\eps,M)$,
and the matrix  elements $L^{\s,\s'}_{\nn,\nn'}(\h,\eps,M)$ can be
extended on the set $\DDD_{0} \setminus \II_{\{\nn,\nn'\}}(\g)$
to $C^{1}$ functions $L^{E \, \s,\s'}_{\nn,\nn'}(\h,\eps,M)$, such that
$L^{E \, \s,\s'}_{\nn,\nn'}(\h,\eps,M) = L^{\s,\s'}_{\nn,\nn'}(\h,\eps,M)$
and $u_{\nn}^{E \, \s}(\h,\eps,M)=  u^{\s}_{\nn}(\h,\eps,M)$
for all $(\eps,M)\in \DDD(2\g)$.
\item The matrix elements $L^{E\,\s,\s'}_{\nn,\nn'}(\h,\eps,M)$ satisfy
for all $(\eps,M) \in \DDD_{0} \setminus \II_{\{\nn,\nn'\}}(\g)$
the bounds
\begin{eqnarray}
& & \left| L^{E \, \s,\s'}_{\nn,\nn'}(\h,\eps,M) \right| \le
{\rm e}^{-\kappa|\nn-\nn'|^\rho} |\h|^{N} K_{1} , \qquad
|\pr_\eps L^{E \, \s,\s'}_{\nn,\nn'}(\h,\eps,M)|
\le {\rm e}^{-\kappa|\nn-\nn'|^\rho} |\h|^{N} K_{1} |p_{\nn}|^{c_{0}} ,
\nonumber \\
& & 
|\pr_\eta L^{E\,\s,\s'}_{\nn,\nn'}(\h,\eps,M)|
\le{\rm e}^{-\kappa|\nn-\nn'|^\rho} N\,|\h|^{N-1} K_{1} , \nonumber 
\end{eqnarray}
for all $(\eps,M) \in \DDD_{0} \setminus \cup
\II_{\{\nn,\nn'\}}(\g)$, where the union is taken over all
the resonant pairs $\{\nn,\nn'\}$, one has
\begin{equation}
\left| \pr_{M} L^{E}(\h,\eps,M) \right|_{\rm{op}}
\le  \sum_{\nn\in \gotO} \sum_{\nn'\in \CC_{\nn}} \sum_{\s,\s'=\pm}
\left| \pr_{M^{\s,\s'}_{\nn,\nn'}} L^{ E}(\h,\eps,M) \right|_{\kappa}
\le |\h|^{N} K_{1} , \nonumber
\end{equation}
and, finally, one has
$$ \left| u^{E \, \s}_{\nn}(\h,\eps,M) \right|
\le |\h|^{N} K_{1} {\rm e}^{-\kappa |\nn|^{1/2}} , $$
uniformly for $(\eps,M)\in \DDD_{0}$.
\end{enumerate}
\end{prop}
%%%%%%%%%%%%%%%%%%%%%%%%%%%%%%%%%%%%%%%%%%%%%%%%%%%%%%%%%%%%%%%%%%%%%%%%%%

%%%%%%%%%%%%%%%%%%%%%%%%%%%%%%%%%%%%%%%%%%%%%%%%%%%%%%%%%%%%%%%%%%%%%%%%%%
\begin{rmk} \label{rmk:11}
In our analysis we choose $M\in B_{\kappa}$ because eventually
we obtain $L \in B_{\kappa}$, but -- as the bound on the
$M$-derivative in item 5 of Proposition \ref{prop:1} suggests -- we could
also take $M$ in a larger space, say $B_{\io}$ with norm $|M|_{\io}=
\sup_{\nn,\nn'\in\gotO}\sup_{\s,\s'=\pm}|M^{\s,\s'}_{\nn,\nn'}|$.
\end{rmk}
%%%%%%%%%%%%%%%%%%%%%%%%%%%%%%%%%%%%%%%%%%%%%%%%%%%%%%%%%%%%%%%%%%%%%%%%%%

Once we have proved Proposition \ref{prop:1}, we solve
the compatibility equation (\ref{eq:3.5}) for the extended
counterterms $L^{E}(\eps^{1/N},\eps,M)$, which are
well defined provided we choose $\eps<\eps_{0}$, with
$\eps_{0}=\h_{0}^{N}$.

%%%%%%%%%%%%%%%%%%%%%%%%%%%%%%%%%%%%%%%%%%%%%%%%%%%%%%%%%%%%%%%%%%%%%%%%%%
\begin{prop} \label{prop:2}
There exist  $C^1$ functions $\eps \to (\eps,M^{\s,\s'}_{\nn,\nn'}
(\eps))$ from $\gotE_{0} \setminus\overline\II_{\{\nn,\nn'\}}(\g)
\to \DDD_{0}$, with an appropriate choice of $C_0$ in
Definition \ref{def:9}, such that the following holds.
\begin{enumerate}
\item $M(\eps)$ verifies the equation
\begin{equation}
M^{\s,\s'}_{\nn,\nn'}(\eps)=  L^{E \, \s,\s'}_{\nn,\nn'}(\eps^{1/N},\eps,M(\eps)) , 
\label{eq:3.6} \end{equation}
and the bounds
$$ \left| M^{\s,\s'}_{\nn,\nn'}(\eps) \right|\leq K_{2}
\eps {\rm e}^{-\kappa|\nn-\nn'|^\rho} , \qquad
\left| \pr_\eps M^{\s,\s'}_{\nn,\nn'}(\eps) \right|
\le K_{2} \left( 1 + \eps p_{\nn}^{c_{0}}(\eps)
\right) {\rm e}^{-\kappa|\nn-\nn'|^{\rho}} , $$ 
for a suitable constant $K_{2}$.
\item The set $ \gotE(2\g) := \left\{ \eps \in \gotE_{0}:
(\eps,M(\eps))\in \DDD(2\g) \right\} $
has large relative Lebesgue measure, namely
$\lim_{\eps\to 0^{+}}\eps^{-1}{\rm meas}(\gotE(2\g) \cap (0,\eps))=1$.
\end{enumerate}
\end{prop}
%%%%%%%%%%%%%%%%%%%%%%%%%%%%%%%%%%%%%%%%%%%%%%%%%%%%%%%%%%%%%%%%%%%%%%%%%%

%%%%%%%%%%%%%%%%%%%%%%%%%%%%%%%%%%%%%%%%%%%%%%%%%%%%%%%%%%%%%%%%%%%%%%%%%%
\subsection{Proof of Theorem \ref{thm:1}}
%%%%%%%%%%%%%%%%%%%%%%%%%%%%%%%%%%%%%%%%%%%%%%%%%%%%%%%%%%%%%%%%%%%%%%%%%%

By item 1 in Proposition \ref{prop:1} for all $(\eps,M)\in \DDD(\g)$
we can find a matrix $L(\h,\eps,M)$
so that there exists a unique solution $u^{\s}_{\nn}(\h,\eps,M)$
of (\ref{eq:3.3}) and (\ref{eq:3.4}) for all $|\h|\le \h_{0}$,
for a suitable $\eta_{0}$, and for $\eps_{0}$ small enough.
By item 3 in Proposition \ref{prop:1} the matrix elements
$L^{\s,\s'}_{\nn,\nn'}(\h,\eps,M)$
and the solution $u^{\s}_{\nn}(\h,\eps,M)$ can be extended
to $C^1$ functions -- denoted by $L^{E \, \s,\s'}_{\nn,\nn'}(\h,\eps,M)$ and
$u^{E \, \s}_{\nn}(\h,\eps,M)$ -- for all
$(\eps,M)\in\DDD_{0}\setminus\II_{\{\nn,\nn'\}}(\g)$
and for all $(\eps,M)\in \DDD_{0}$, respectively.
Moreover $L^{E \, \s,\s'}_{\nn,\nn'}(\h,\eps,M)= L^{\s,\s'}_{\nn,\nn'}(\h,\eps,M)$ and 
$u^{E \, \s}_{\nn}(\h,\eps,M)= u^{\s}_{\nn}(\h,\eps,M)$
for all $(\eps,M)\in \DDD(2\g)$.

Equation (\ref{eq:3.3}) coincides with our original (\ref{eq:3.2})
provided the compatibility equation (\ref{eq:3.5}) is satisfied.
Now we fix $\eps_{0}<\h_{0}^{N}$ so that $L^{E}(\eps^{1/N},\eps,M)$
and  $u^{E \, \s}_{\nn}(\eps^{1/N},\eps,M)$ are well defined
for $|\eps|<\eps_{0}$. By item 1 in Proposition \ref{prop:2},
there exists a  matrix $M(\eps)$
which satisfies the extended compatibility equation (\ref{eq:3.6}).
Finally by item 2 in Proposition \ref{prop:2} the Cantor set 
$\gotE(2\g)$ is well defined and of large relative measure.
 
For all $\eps\in \gotE(2\g)$ the pair $(\eps,M(\eps))$ is by definition
in $\DDD(2\g)$, so that by item 3 in Proposition \ref{prop:1} one has
$L^{\s,\s'}_{\nn,\nn'}(\eps^{1/N},\eps,M(\eps))=
L^{E \, \s,\s'}_{\nn,\nn'}(\eps^{1/N},\eps,M(\eps))$
and $u^{\s}(\eps^{1/N},\eps,M(\eps);x,t)=u^{E \, \s}
(\eps^{1/N},\eps,M(\eps);x,t)$,
and hence $u^{\s}_{\nn}(\eps^{1/N},\eps,M(\eps))$ solves (\ref{eq:3.3}) for 
$\h=\eps^{1/N}$. So, by item 1 in Proposition \ref{prop:2},
$M(\eps)$ solves the true compatibility equation (\ref{eq:3.5})
for all $\eps\in \gotE(2\g)$.
Then $u^{\s}(\eps^{1/N},\eps,M(\eps);x,t)$ is a true nontrivial solution
of (\ref{eq:3.3}) and (\ref{eq:3.4}) in $\gotE(2\g)$.
Then by setting $\gotE=\gotE(2\g)$ the result follows.

%%%%%%%%%%%%%%%%%%%%%%%%%%%%%%%%%%%%%%%%%%%%%%%%%%%%%%%%%%%%%%%%%%%%%%%%%
%%%%%%%%%%%%%%%%%%%%%%%%%%%%%%%%%%%%%%%%%%%%%%%%%%%%%%%%%%%%%%%%%%%%%%%%%
\zerarcounters
\section{Tree expansion}
\label{sec:4}
%%%%%%%%%%%%%%%%%%%%%%%%%%%%%%%%%%%%%%%%%%%%%%%%%%%%%%%%%%%%%%%%%%%%%%%%%
%%%%%%%%%%%%%%%%%%%%%%%%%%%%%%%%%%%%%%%%%%%%%%%%%%%%%%%%%%%%%%%%%%%%%%%%%

%%%%%%%%%%%%%%%%%%%%%%%%%%%%%%%%%%%%%%%%%%%%%%%%%%%%%%%%%%%%%%%%%%%%%%%%%%
\subsection{Recursive equations}
%%%%%%%%%%%%%%%%%%%%%%%%%%%%%%%%%%%%%%%%%%%%%%%%%%%%%%%%%%%%%%%%%%%%%%%%%%

In this section we find a formal solution $u^{\s}_{\nn},L$
of (\ref{eq:3.3}) and (\ref{eq:3.4}) as a power series on $\h$;
the solution $u^{\s}_{\nn},L$ depends  on  the matrix 
$M$ and it will be written in the form of a tree expansion.

We assume for $u^{\s}_{\nn}(\h,\eps,M)$ for all $\nn\in\gotP$
and for the matrix $L(\h,\eps,M)$
a formal series expansion in $\h$:
\begin{equation}
u^{\s}_{\nn}(\h,\eps,M)= \sum_{k=N}^{\io}\h^{k} u_{\nn}^{(k)\s} ,
\qquad L(\h,\eps,M) = \sum_{k=N}^{\io}\h^{k} L^{(k)} ,
\label{eq:4.1} \end{equation}
with the Ansatz that $L^{(k)\s,\s'}_{\nn,\nn'}=0$ if either
$\bar\chi_{1}(\delta_{\nn}(\eps))\bar\chi_{1}(\delta_{\nn'}(\eps))=0$
or the pair $\{\nn,\nn'\}$ is not resonant, so that
$L=\widehat \chi_{1} L \widehat \chi_{1}$. 
We set also $u^{(k)\s}_{\nn}=0$  for all $k\le N$ and $\nn,\nn'\in\gotP$
same for $L^{(k)\s,\s'}_{\nn,\nn'}$ for $\nn,\nn'\in\gotO$.

For $\nn\in\gotQ$ we set
\begin{equation}
u^{\s}_{\nn}(\h,\eps,M) = u_{\nn}^{(0)\s} +
\sum_{k=N}^{\io} \h^{k} u_{\nn}^{(k)\s} .
\label{eq:4.2} \end{equation}
with $u^{(0)+}_{\nn}=u^{(0)}_{\nn}$ and $u^{(0)-}_{\nn}= \overline{
u^{(0)}_{\nn}}$ (cf. item 2 in Hypothesis \ref{hp:2} for notations).
Again we set $u^{(k)\s}_{\nn}=0$ for $0<k<N$ and $\nn\in\gotQ$.

Inserting the series expansions (\ref{eq:4.1}) and (\ref{eq:4.2})
into (\ref{eq:3.3}) we obtain
\begin{equation}
\begin{cases}
u^{(k)\s}_{\nn} = {\displaystyle
\frac{f_{\nn}^{(k-N)\s}}{\delta_{\nn}(\eps)} } , & \nn\in \gotR , \\
u^{(k)\s}_{\nn} = \!\!\!\!\!\!\!\!
{\displaystyle \sum_{\nn'\in\gotQ,\,\s'=\pm}
(J^{-1})^{\s,\s'}_{\nn,\nn'} f_{\nn'}^{(k)\s'} } , & \nn\in\gotQ , \\
\left( \DD(\eps)+\widehat M \right) U^{(k)} = F^{(k-N)} +
{\displaystyle \sum_{r=N}^{k-N} L^{(r)} U^{(k-r)} } . &
\end{cases}
\label{eq:4.3} \end{equation}
%
%where, for $\nn\in\gotQ\cup\gotR$, we have set $F^{(k)}_{\nn}=
%f^{(k)}_{\nn}$, with $f_{\nn}^{(k)}$ the $\h$-Taylor expansion
%at order $k$ of $f_{\nn}$.

%%%%%%%%%%%%%%%%%%%%%%%%%%%%%%%%%%%%%%%%%%%%%%%%%%%%%%%%%%%%%%%%%%%%%%%%%%
\subsection{Multiscale analysis}
%%%%%%%%%%%%%%%%%%%%%%%%%%%%%%%%%%%%%%%%%%%%%%%%%%%%%%%%%%%%%%%%%%%%%%%%%%

%%%%%%%%%%%%%%%%%%%%%%%%%%%%%%%%%%%%%%%%%%%%%%%%%%%%%%%%%%%%%%%%%%%%%%%%%%
\begin{defi} \label{def:11}
\textbf{ (The scale functions).}
Let $\chi$ be a non-increasing function $C^{\io}(\RRR_{+},[0,1])$,
such that $\chi(x)=0$ if $x\ge 2\g$ and $\chi(x)=1$ if $x\le \g$;
moreover one has $|\pr_{x}\chi(x)|\le \Gamma \g^{-1}$
for some positive constant $\Gamma$.
Let $\chi_{h}(x)=\chi(2^{h}x)-\chi(2^{h+1}x)$ for
$h \ge 0$, and $\chi_{-1}(x)=1-\chi(x)$.
\end{defi}
%%%%%%%%%%%%%%%%%%%%%%%%%%%%%%%%%%%%%%%%%%%%%%%%%%%%%%%%%%%%%%%%%%%%%%%%%%

Recall that for each $\eps$ the matrix $A=\DD(\eps)+ \widehat M$
is block diagonal with a diagonal part whose eigenvalues are
larger than $\bar\g > \g$ and a list of $C_{1} p_{\nn}^{\al}(\eps)
\times C_{1} p_{\nn}^{\al}(\eps)$ blocks  $A^\nn$ containing
small entries. In the following if $A^\nn$ is invertible -- i.e.
if $x_{\nn}\neq 0$ -- we will denote the entries of $(A^\nn)^{-1}$
by $(A^{-1})^{\s,\s'}_{\nn,\nn'}$ even though it may be
possible that the whole matrix $A$ is not invertible.

%%%%%%%%%%%%%%%%%%%%%%%%%%%%%%%%%%%%%%%%%%%%%%%%%%%%%%%%%%%%%%%%%%%%%%%%%%
\begin{defi} \label{def:12}
\textbf{ (Propagators).}
For $\nn,\nn'\in \gotO$, we define the {\rm propagators}
\begin{equation}
(G_{i,h})_{\nn,\nn'}^{\s,\s'} = \begin{cases}
\chi_{h}(x_{\nn}(\eps)) \,\bar\chi_{1}(\delta_{\nn}(\eps))
\bar\chi_{1}(\delta_{\nn'}(\eps)) (A^{-1})^{\s,\s'}_{\nn,\nn'} , &
\text{ if } i=1 \text{ and } \chi_{h}(x_{\nn}(\eps)) \neq 0 , \\
\bar\chi_{0}(\delta_{\nn}(\eps))\,\delta^{-1}_{\nn}(\eps) , &
\text{ if } i=0, \; \nn=\nn', \; \s=\s' \text{ and } h=-1 , \\
0 , & \text{ otherwise.} \end{cases} \nonumber
\end{equation}
\end{defi}
%%%%%%%%%%%%%%%%%%%%%%%%%%%%%%%%%%%%%%%%%%%%%%%%%%%%%%%%%%%%%%%%%%%%%%%%%%

In terms of the propagators we obtain
\begin{equation}
A^{-1}=\sum_{i=0,1} \sum_{h=-1}^{\io} G_{i,h},
\label{eq:4.4} \end{equation}
which provides the multiscale decomposition. Notice that if 
$(A^{-1})^{\s,\s'}_{\nn,\nn'}\neq 0$ then $x_{\nn}(\eps)=x_{\nn'}(\eps)$
(see Remark \ref{rmk:10}),
so that the matrices $G_{i,h}$ are indeed self-adjoint.

%%%%%%%%%%%%%%%%%%%%%%%%%%%%%%%%%%%%%%%%%%%%%%%%%%%%%%%%%%%%%%%%%%%%%%%%%%
\begin{rmk} \label{rmk:12}
Only the propagator $G_{1,h}$ can produce small divisors
while the propagator  $G_{0,-1}$ is diagonal and of order one.
Hence, there exists a positive constant $C$
such that we can bound the propagators as
\begin{equation}
\left| G_{0,-1} \right|_{\io} \le C \g^{-1} ,
\qquad \left| (G_{1,h})_{\nn,\nn'}^{\s,\s'} \right| \le
2^{h} C \g^{-1} p_{\nn}^{-\xi}(\eps) \sqrt{p_{\nn}^{\al}(\eps)} ,
\label{eq:4.5} \end{equation}
where the condition $d_{\nn}(\eps)\le 2 C_{1} p_{\nn}^{\al}(\eps)$
-- cf. Remark \ref{rmk:10} -- and
item 2 of Lemma \ref{lem:4} have been used.
\end{rmk}
%%%%%%%%%%%%%%%%%%%%%%%%%%%%%%%%%%%%%%%%%%%%%%%%%%%%%%%%%%%%%%%%%%%%%%%%%%

We write $L^{(k)}$ in (\ref{eq:4.1}) as
\begin{equation}
L^{(k)\s_1,\s_2}_{\nn_{1},\nn_{2}} = \sum_{h=-1}^{\io}
\chi_{h}(x_{\nn_{1}}(\eps))  L^{(k)\s_1,\s_2}_{h,\nn_{1},\nn_2} ,
\label{4.6} \end{equation}
for all resonant pairs $\{\nn_{1},\nn_{2}\}$; we denote by $L^{(k)}_{h}$
the matrix with entries $L^{(k)\s_1,\s_2}_{h,\nn_{1},\nn_2}$.
Finally we set
\begin{equation}
U^{(k)} = \sum_{i=0,1} \sum_{h=-1}^{\infty} U^{(k)}_{i,h} ,
\label{eq:4.7} \end{equation}
so that (\ref{eq:4.3}) gives
\begin{equation}
\begin{cases}
u^{(k)\s}_{\nn} = {\displaystyle \sum_{\nn'\in\gotQ}
(J^{-1})^{\s,\s'}_{\nn,\nn'} \, f^{(k)\s'}_{\nn'} } , & \nn \in \gotQ , \\
u^{(k)\s}_{\nn} = {\displaystyle
\frac{f^{(k-N)\s}_{\nn}}{\de_{\nn}(\eps)} } , & \nn \in \gotR , \\ 
U^{(k)}_{i,h} = {\displaystyle
G_{i,h} F^{(k-N)} + \de(i,1) \, G_{1,h}
\sum_{h_{1}=-1}^{\io} \sum_{r=N}^{k-N}L^{(r)}_{h}
U^{(k-r)}_{1,h_{1}} } , & i=0,1,\;h\ge-1, \end{cases}
\label{eq:4.8} \end{equation}
which are the recursive equations we want to study.

%%%%%%%%%%%%%%%%%%%%%%%%%%%%%%%%%%%%%%%%%%%%%%%%%%%%%%%%%%%%%%%%%%%%%%%%%%
\subsection{Diagrammatic rules}
%%%%%%%%%%%%%%%%%%%%%%%%%%%%%%%%%%%%%%%%%%%%%%%%%%%%%%%%%%%%%%%%%%%%%%%%%%

A connected graph $\GG$ is a collection of points (vertices)
and lines connecting all of them. We denote with
$V(\GG)$ and $L(\GG)$ the set of nodes and the set of lines,
respectively. A path between two nodes is the minimal subset of
$L(\GG)$ connecting the two nodes. A graph is planar if it
can be drawn in a plane without graph lines crossing.

%%%%%%%%%%%%%%%%%%%%%%%%%%%%%%%%%%%%%%%%%%%%%%%%%%%%%%%%%%%%%%%%%%%%%%%%%%
\begin{defi} \label{def:13}
\textbf{ (Trees).}
A tree is a planar graph $\GG$ containing no closed loops.
One can consider a tree $\GG$ with a single special node $v_{0}$:
this introduces a natural partial ordering on the set
of lines and nodes, and one can imagine that each line
carries an arrow pointing toward the node $v_{0}$.
We can add an extra (oriented) line $\ell_{0}$ exiting the special
node $v_{0}$; the added line $\ell_{0}$ will be called the root line
and the point it enters (which is not a node) will
be called the root of the tree. In this way we obtain
a rooted tree $\theta$ defined by $V(\theta)=V(\GG)$
and $L(\theta)=L(\GG)\cup\ell_{0}$. A labelled tree is a
rooted tree $\theta$ together with a label function defined on
the sets $L(\theta)$ and $V(\theta)$.
\end{defi}
%%%%%%%%%%%%%%%%%%%%%%%%%%%%%%%%%%%%%%%%%%%%%%%%%%%%%%%%%%%%%%%%%%%%%%%%%%

We shall call {\it equivalent} two rooted trees which can be transformed
into each other by continuously deforming the lines in the plane
in such a way that the latter do not cross each other
(i.e. without destroying the graph structure).
We can extend the notion of equivalence also to labelled trees,
simply by considering equivalent two labelled trees if they
can be transformed into each other in such a way that also
the labels match. 

Given two nodes $v,w\in V(\theta)$, we say that $v \prec w$
if $w$ is on the path connecting $v$ to the root line.
We can identify a line with the nodes it connects;
given a line $\ell=(w,v)$ we say that $\ell$
enters $w$ and exits (or comes out of) $v$,
and we write $\ell=\ell_{v}$. 
Given two comparable lines $\ell$ and $\ell_{1}$,
with $\ell_{1} \prec \ell$, we denote with $\PPP(\ell_{1},\ell)$
the path of lines connecting $\ell_{1}$ to $\ell$; by definition
the two lines $\ell$ and $\ell_{1}$ do not belong to $\PPP(\ell_{1},
\ell)$. We say that a node $v$ is along the path $\PPP(\ell_{1},
\ell)$ if at least one line entering or exiting $v$ belongs to the path.
If $\PPP(\ell_{1},\ell)=\emptyset$ there is only one node $v$
along the path (such that $\ell_{1}$ enters $v$ and $\ell$ exits $v$).

%%%%%%%%%%%%%%%%%%%%%%%%%%%%%%%%%%%%%%%%%%%%%%%%%%%%%%%%%%%%%%%%%%%%%%%%%%
\begin{defi} \label{def:14}
\textbf{ (Lines and nodes).}
We call {\rm internal nodes} the nodes such that there is at least
one line entering them; we call {\rm internal lines} the lines
exiting the internal nodes. We call {\rm end-nodes} the nodes
which have no entering line. We denote with $L(\theta)$, $V_{0}(\theta)$
and $E(\theta)$ the set of lines, internal nodes and end-nodes,
respectively. Of course $V(\theta)=V_{0}(\theta)\cup E(\theta)$.
\end{defi}
%%%%%%%%%%%%%%%%%%%%%%%%%%%%%%%%%%%%%%%%%%%%%%%%%%%%%%%%%%%%%%%%%%%%%%%%%%

We associate with the nodes (internal nodes and end-nodes) and lines
of any tree $\theta$ some labels, according to the following rules.

%%%%%%%%%%%%%%%%%%%%%%%%%%%%%%%%%%%%%%%%%%%%%%%%%%%%%%%%%%%%%%%%%%%%%%%%%%
\begin{defi} \label{def:15}
\textbf{ (Diagrammatic rules).}
\null\hfil\null
\begin{enumerate}
%\label 1
\item For each node $v$ there are $p_{v} \ge 0$ entering lines.
If $p_{v}=0$ then $v\in E(\theta)$, if $p_{v} >0$ then either $p_{v}=1$
or $p_{v} \ge N+1$ and $v\in V_0(\theta)$. If $L(v)$ is the set 
of lines entering $v$ one has $p_{v}=|L(v)|$.
%\label 2
\item With each internal line $\ell\in L(\theta)$
one associates a label $q$, $p$ or $r$. We say that
$\ell$ is a $p$-line, a $q$-line or an $r$-line, respectively,
and we call $L_{q}(\theta)$, $L_{p}(\theta)$ and $L_{r}(\theta)$
the set of internal lines $\ell\in L(\theta)$ which are $q$-lines,
$p$-lines and $r$-lines, respectively.
If $p_{v}=1$ then the line $\ell$ exiting $v$ and the line
$\ell_{1}$ entering $v$ are both $p$-lines.
%label 3
\item With each line $\ell\in L(\theta)$
one associates the {\rm type} label $i_{\ell}=0,1$.
%label 4
\item With each line $\ell\in L(\theta)$ except the root line $\ell_{0}$
one associates a {\rm sign} label $\s_\ell=\pm $.
%label 5
\item With each internal line $\ell\in L(\theta)$ one
associates the {\rm momenta} $(\nn_{\ell},\nn_{\ell}') \in \ZZZ^{D+1}
\times \ZZZ^{D+1}$.
%label 6
\item With each line $\ell\in L(\theta)$ exiting an end-node
one associates the {\rm momentum} $\nn_{\ell}$.
%label 7
\item With each line $\ell\in L(\theta)$
one associates the {\rm scale} label $h_{\ell}\in \NNN \cup\{-1,0\}$.
%label 8
\item With each end-node $v\in E(\theta)$ one associates
the {\rm mode} label $\nn_{v}\in \gotQ$, the {\rm order} label
$k_{v}=0$, and the {\rm sign} label $\s_{v}=\pm$.
%label 9
\item With each internal node $v\in V_{0}(\theta)$ one associates
the {\rm mode} label $m_{v}\in\ZZZ^{D}$,
the {\rm order} label $k_{v}\in\NNN$, and
the {\rm sign} label $\s_{v}=\pm $.
%label 10
\item For each internal node $v\in V_{0}(\theta)$
one defines $r_{v}$ as the number of lines  $\ell\in
L(v)$ with $\s_\ell=\s_{v}$, and one sets $s_{v}= p_{v}-r_{v}$.
%label 11
\item If a line $\ell\in L(\theta)$ is not a $p$-line one sets $i_{\ell}=0$
%label 12
\item If a line $\ell\in L(\theta)$ has $i_{\ell}=0$, then $h_{\ell}=-1$.
%label 13
\item Let $\ell\in L(\theta)$ be an internal line.
If $\ell$ is a $p$-line with $i_{\ell}=0$,
then $\nn_{\ell}=\nn'_{\ell}$.
If $\ell$ is a $p$-line with $i_{\ell}=1$,
then $\{\nn_{\ell},\nn'_{\ell}\}$ is a resonant pair.
If $\ell$ is a $q$-line,
then $\nn_{\ell},\nn'_{\ell}\in \gotQ$.
If $\ell$ is an $r$-line,
then $\nn_{\ell}=\nn_{\ell}' \in\gotR$.
%label 14
\item If $\ell$ exits an end-node $v$, then $\nn_{\ell}=\nn_{v}$.
%label 15
\item \label{itemscaledifference}
If two $p$-lines $\ell$ and $\ell'$
have $i_{\ell}=i_{\ell'}=1$ and are such that
$\{\nn_{\ell},\nn_{\ell}',\nn_{\ell'},\nn_{\ell'}'\}$
is a resonant set, then $|h_{\ell}-h_{\ell'}| \le 1$.
%label 16
\item If $\ell\in L(\theta)$ exits an end-node $v \in E(\theta)$,
then one sets $\s_{\ell}=\s_{v}$.
%label 17
\item \label{itemconservation}
If $\ell$ is the line exiting $v$ and $\ell_{1},\ldots,
\ell_{p_{v}}$ are the lines entering $v$ one has
$$ \nn_{\ell}' = (0,m_{v}) + \s_{v}(\s_{\ell_{1}} \nn_{\ell_{1}} +
\ldots + \s_{\ell_{p_{v}}} \nn_{\ell_{p_{v}}}) = (0,m_{v}) +
\s_{v}\sum_{\ell'\in L(v)} \s_{\ell'} \nn_{\ell'} , $$
which represents a conservation rule for the momenta.
%
%label 18
\item Given an internal node $v\in V_{0}(\theta)$,
if $p_{v}=1$ one has $k_{v} \ge N$, while if $p_{v} \ge N$
one has $k_{v}=p_{v}-1$.
%label 19
\item Given an internal node $v\in V_{0}(\theta)$,
if $p_{v}=1$, let $\ell_{1}$ be the line
entering $v$ and $\ell$ be the line exiting $v$.
One has $i_{\ell_{1}}=i_{\ell}=1$ and
$\{\nn'_{\ell},\nn_{\ell_{1}}\} $ is a resonant pair.
%label 20
\item With each end-node $v\in E(\theta)$ one associates
the {\rm node factor} $\h_{v}=u^{(0)\s_{v}}_{\nn_{v}}$;
cf. item 2 in Hypothesis \ref{hp:2} and (\ref{eq:4.2}) for notations.
%label 21
\item With each internal node $v\in V_{0}(\theta)$ with $p_{v}>1$ one
associates the {\rm node factor} $\h_{v}=a^{\s_{v}}_{r_{v},s_{v},m_{v}}$,
where $a^{\s}_{r,s,m}$ satisfies equation (\ref{eq:1.11}).
%label 22
\item With each internal node $v\in V_{0}(\theta)$ with $p_{v}=1$ one
associates the {\rm node factor} $\h_{v}=L^{(k_{v})\s_{v},
\s_{\ell_1}}_{{h_{\ell}},\nn_{\ell}',\nn_{\ell_{1}}}$,
still to be defined (see Definition \ref{def:25} below),
where $\ell$ and $\ell_{1}$ are the lines exiting and entering $v$,
respectively.
%label 23
\item One associates with each line $\ell\in L(\theta)$
a {\rm line propagator} $g_{\ell} \in \CCC$ with the following rules.
If $\ell$ is a $p$-line  exiting the internal node $v$ 
one sets $g_{\ell} := (G_{i_{\ell},h_{\ell}})^{\s_\ell,\s_{v}}
_{\nn_{\ell},\nn'_{\ell}}$, if $\ell$ is an $r$-line
one sets $g_{\ell}:=1/\de_{\nn_{\ell}}(\eps)$, if $\ell$ is
a $q$-line exiting the internal node $v$ one sets $g_{\ell} :=
(J^{-1})^{\s_\ell,\s_{v}}_{\nn_{\ell},\nn_{\ell}'}$,
if $\ell$ exits an end-node one sets $g_{\ell}=1$.
%label 24
\item One defines the {\rm order} of the tree $\theta$ as
$$ k(\theta) := \sum_{v\in V(\theta)} k_{v} , $$
the {\rm momentum} of $\theta$ as the momentum
$\nn_{\ell}$ of the root line $\ell$, and the {\rm sign} of $\theta$ as
the sign $\s_{v_{0}}$ of the node $v_{0}$ which the root line exits.
\end{enumerate}
\end{defi}
%%%%%%%%%%%%%%%%%%%%%%%%%%%%%%%%%%%%%%%%%%%%%%%%%%%%%%%%%%%%%%%%%%%%%%%%%%

%%%%%%%%%%%%%%%%%%%%%%%%%%%%%%%%%%%%%%%%%%%%%%%%%%%%%%%%%%%%%%%%%%%%%%%%%%
\begin{defi} \label{def:16}
\textbf{ (The sets of trees
$\boldsymbol\Theta_{\boldsymbol\nn}^{\boldsymbol( \boldsymbol k
\boldsymbol)\boldsymbol \s}$ and $\boldsymbol\Theta$).}
We call $\Theta^{(k)\s}_{\nn}$ the set of all the nonequivalent
trees of order $k$, momentum $\nn$ and sign $\s$, defined according
to the diagrammatic rules of Definition \ref{def:15}.
We call $\Theta$ the sets of trees belonging to $\Theta^{(k)\s}_{\nn}$
for some $k\ge 1$, $\s=\pm$ and $\nn\in\ZZZ^{D+1}$.
\end{defi}
%%%%%%%%%%%%%%%%%%%%%%%%%%%%%%%%%%%%%%%%%%%%%%%%%%%%%%%%%%%%%%%%%%%%%%%%%%

%%%%%%%%%%%%%%%%%%%%%%%%%%%%%%%%%%%%%%%%%%%%%%%%%%%%%%%%%%%%%%%%%%%%%%%%%%
\subsection{Clusters and resonances}
%%%%%%%%%%%%%%%%%%%%%%%%%%%%%%%%%%%%%%%%%%%%%%%%%%%%%%%%%%%%%%%%%%%%%%%%%%

%%%%%%%%%%%%%%%%%%%%%%%%%%%%%%%%%%%%%%%%%%%%%%%%%%%%%%%%%%%%%%%%%%%%%%%%%%
\begin{defi} \label{def:17}
\textbf{ (Clusters).}
Given a tree $\theta\in\Theta^{(k)\s}_{\nn}$ a {\rm cluster} $T$
on scale $h$ is a connected maximal set of nodes and lines such that
all the lines $\ell$ have a scale label $\le h$ and at least
one of them has scale $h$; we shall call $h_{T}=h$ the scale of
the cluster. We shall denote by $V(T)$, $V_{0}(T)$ and $E(T)$
the set of nodes, internal nodes and the set of end-nodes,
respectively, which are contained inside the cluster $T$,
and with $L(T)$ the set of lines connecting them.
Finally $k(T)=\sum_{v\in V(T)}k_{v}$ will be called the order of $T$.
\end{defi}
%%%%%%%%%%%%%%%%%%%%%%%%%%%%%%%%%%%%%%%%%%%%%%%%%%%%%%%%%%%%%%%%%%%%%%%%%%

An inclusion relation is established between clusters,
in such a way that the innermost clusters are the clusters
with lowest scale, and so on. A cluster $T$ can have an arbitrary
number of lines entering it ({\it entering lines}), but only one
or zero line coming out from it ({\it exiting line} or {\it root line}
of the cluster); we shall denote the latter (when it exists)
with $\ell_{T}^{}$. Notice that, by definition, $|V(T)|>1$
and all the entering and exiting lines have $i_{\ell}=1$. 

%%%%%%%%%%%%%%%%%%%%%%%%%%%%%%%%%%%%%%%%%%%%%%%%%%%%%%%%%%%%%%%%%%%%%%%%%%
\begin{defi} \label{def:18}
\textbf{ (Resonances).}
We call {\rm resonance} on scale $h$ a cluster $T$ on scale $h_{T}=h$
such that
\begin{enumerate}
\item the cluster has only one entering line $\ell_{T}^{1}$ and
one exiting line $\ell_{T}^{}$ of scale $h_{\ell_{T}^{}} \ge h+2$,
\item one has that $\{\nn_{\ell_{T}^{}}',\nn_{\ell_{T}^{1}}\}$ is
a resonant pair and $\min\{|\nn_{\ell_{T}^{1}}|,|\nn_{\ell_{T}^{}}'|\}
\ge 2^{(h-2)/\tau}$,
\item for all $\ell\in\PPP(\ell_{T}^{1},\ell_{T}^{})$ with $i_{\ell}=1$
the pair $\{\nn_{\ell}',\nn_{\ell_{T}^{1}}\}$ is not resonant,
\item for all $\ell\in L(T) \setminus \PPP(\ell_{T}^{1},\ell_{T}^{})$
the pair $\{\nn_{\ell}',\nn_{\ell_{T}^{1}}\}$ is not resonant.
\end{enumerate}
The line $\ell_{T}^{}$ of a resonance will be
called the root line of the resonance.
\end{defi}
%%%%%%%%%%%%%%%%%%%%%%%%%%%%%%%%%%%%%%%%%%%%%%%%%%%%%%%%%%%%%%%%%%%%%%%%%%

%%%%%%%%%%%%%%%%%%%%%%%%%%%%%%%%%%%%%%%%%%%%%%%%%%%%%%%%%%%%%%%%%%%%%%%%%%
\begin{defi} \label{def:19}
\textbf{ (The sets of trees $\boldsymbol\RR_{\boldsymbol h,\boldsymbol
\nn,\nn'}^{\boldsymbol( \boldsymbol k \boldsymbol)
\boldsymbol\s,\boldsymbol\s'}$ and
$\boldsymbol\RR$).}
For $k\ge N$, $h \ge 1$ and a resonant pair $\{\nn,\nn'\}$ such that
$\min\{|\nn|,|\nn'|\} \ge 2^{(h-2)/\tau}$, we define
$\RR_{h,\nn,\nn'}^{(k)\s\s'}$ as the set of trees with
the following differences with respect to $\Theta_{\nn}^{(k)\s}$.
\begin{enumerate} 
\item There is a single end-node, called $e$,
with node factor $\h_{e}=1$
(but no label no labels $\nn_{e}$ nor $\s_{e}$).
\item The line $\ell_{e}$ exiting $e$ is a $p$-line. We associate
with $\ell_{e}$ the labels $\nn_{\ell_{e}}=\nn'$, $\s_{\ell}=\s'$,
and $i_{\ell_{e}}=1$ (but no labels $\nn_{\ell}'$ nor $h_{\ell}$),
and the corresponding line propagator is $g_{\ell_{e}}=
\bar\chi_1(\de_{\nn'}(\eps))$.
\item The root line $\ell_{0}$ is a $p$-line. We associate with
$\ell_{0}$ the labels $i_{\ell_{0}}=1$ and $\nn_{\ell_{0}}'=\nn$
(but no labels $\nn_{\ell_{0}}$ nor $h_{\ell_{0}}$),
and the corresponding line propagator is $g_{\ell_{0}}=
\bar\chi_{1}(\de_{\nn}(\eps))$. Let $v_{0}$ be the node which
the line $\ell_{0}$ exits: we set $\s_{v_{0}}=\s$.
\item One has $\max_{\ell\in L(\theta)\setminus\{\ell_{0},
\ell_{e}\}} h_{\ell}=h$.
\item If $\ell\in\PPP(\ell_{e},\ell_{0})$ is such that
$\{\nn_{\ell}',\nn'\}$ is resonant, then $i_{\ell}=0$.
\item For $\ell\notin\PPP(\ell_{e},\ell_{0})$
one has that $\{\nn_{\ell}',\nn'\}$ is not a resonant pair.
\end{enumerate}
We call $\RR$ the sets of trees belonging to
$\RR_{h,\nn,\nn'}^{(k)\s\s'}$ for some $k\ge 1$, $h\ge 1$,
$\s,\s'=\pm$, and $\nn,\nn\in\gotO$ such that
$\{\nn,\nn'\}$ is resonant and
$\min\{|\nn|,|\nn'|\} \ge 2^{(h-2)/\tau}$.
\end{defi}
%%%%%%%%%%%%%%%%%%%%%%%%%%%%%%%%%%%%%%%%%%%%%%%%%%%%%%%%%%%%%%%%%%%%%%%%

%%%%%%%%%%%%%%%%%%%%%%%%%%%%%%%%%%%%%%%%%%%%%%%%%%%%%%%%%%%%%%%%%%%%%%%%%%
\begin{defi} \label{def:20}
\textbf{ (Clusters for trees in $\boldsymbol\RR$).}
Given a tree $\theta \in \RR$, a cluster $T$
on scale $h_{T} \le h$ is a connected maximal set of nodes
$v\in V(\theta)$ and lines $\ell\in L(\theta) \setminus \{\ell_{0},
\ell_{e}\}$ such that all the lines $\ell$ have a scale label
$\le h_{T}$ and at least one of them has scale $h_{T}$.
\end{defi}
%%%%%%%%%%%%%%%%%%%%%%%%%%%%%%%%%%%%%%%%%%%%%%%%%%%%%%%%%%%%%%%%%%%%%%%%%%

Note that if $\theta\in\RR_{h,\nn,\nn'}^{(k)\s,\s'}$, then
for any cluster $T$ in $\theta$ one necessarily has $h_{T} \le h$.

%%%%%%%%%%%%%%%%%%%%%%%%%%%%%%%%%%%%%%%%%%%%%%%%%%%%%%%%%%%%%%%%%%%%%%%%%%
\begin{defi} \label{def:21}
\textbf{ (Resonances for trees in $\boldsymbol\RR$).}
Given a tree $\theta \in \RR$, a cluster $T$ is
a resonance if the four items of Definition \ref{def:18} are satisfied.
\end{defi}
%%%%%%%%%%%%%%%%%%%%%%%%%%%%%%%%%%%%%%%%%%%%%%%%%%%%%%%%%%%%%%%%%%%%%%%%%%

%%%%%%%%%%%%%%%%%%%%%%%%%%%%%%%%%%%%%%%%%%%%%%%%%%%%%%%%%%%%%%%%%%%%%%%%%%
\begin{rmk} \label{rmk:13}
There is a one-to-one correspondence between resonances $T$
of order $k$ and scale $h$ with $\nn_{\ell_{T}^{1}}=\nn'$,
$\nn'_{\ell_{T}^{}}=\nn$, $\s_{v_{0}}=\s$, $\s_{\ell_{T}^{1}}=\s'$
(here $v_{0}$ is the node which $\ell_{T}^{}$ exits)
and trees $\theta\in\RR^{(k)\s,\s'}_{h,\nn,\nn'}$;
cf. \cite{GP2}, Section 3.4 and Figure 7.
\end{rmk}
%%%%%%%%%%%%%%%%%%%%%%%%%%%%%%%%%%%%%%%%%%%%%%%%%%%%%%%%%%%%%%%%%%%%%%%%%%

%%%%%%%%%%%%%%%%%%%%%%%%%%%%%%%%%%%%%%%%%%%%%%%%%%%%%%%%%%%%%%%%%%%%%%%%%%
\begin{defi} \label{def:22}
\textbf{ (The sets of renormalised trees $\boldsymbol\Theta_{ \boldsymbol
R,\boldsymbol\nn}^{\boldsymbol( \boldsymbol k \boldsymbol)
\boldsymbol\s}$, 
$\boldsymbol\RR_{\boldsymbol R,\boldsymbol h,\boldsymbol \nn,
\boldsymbol \nn'}^{\boldsymbol( \boldsymbol k \boldsymbol)
\boldsymbol\s,\boldsymbol\s'}$,
$\boldsymbol\Theta_{ \boldsymbol R}$ and
$\boldsymbol\RR_{ \boldsymbol R}$).}
We define the set of {\rm re\-norm\-alis\-ed trees}
$\Theta_{R,\nn}^{(k)\s}$ and
$\RR_{R,h,\nn,\nn'}^{(k)\s,\s'}$ as the set of trees defined as
$\Theta_{\nn}^{(k)\s}$ and $\RR_{h,\nn,\nn'}^{(k)\s,\s'}$, respectively,
but with no resonances and no nodes $v$ with $p_{v}=1$.
Analogously we define the sets $\Theta_{R}$ and $\RR_{R}$.
\end{defi}
%%%%%%%%%%%%%%%%%%%%%%%%%%%%%%%%%%%%%%%%%%%%%%%%%%%%%%%%%%%%%%%%%%%%%%%%

In the following it will turn out to be convenient to
introduce also the following set of trees.

%%%%%%%%%%%%%%%%%%%%%%%%%%%%%%%%%%%%%%%%%%%%%%%%%%%%%%%%%%%%%%%%%%%%%%%%%%
\begin{defi} \label{def:23}
\textbf{ (The set of renormalised trees $\boldsymbol\SSSS_{\boldsymbol R,
\boldsymbol h,\boldsymbol \nn,\boldsymbol \nn'}^{\boldsymbol(
\boldsymbol k \boldsymbol)
\boldsymbol\s,\boldsymbol\s'}$ and
$\boldsymbol\SSSS_{ \boldsymbol R}$).}
For $k\ge N$, $h \ge 1$ and $\nn,\nn'\in \gotO$ such that
$|\nn'|\ge 2^{(h-2)/\tau}$ we define the
set of renormalised trees $\SSSS_{R,h,\nn,\nn'}^{(k)\s,\s'}$
as the set of trees with the following differences with respect
to $\RR^{(k)\s,\s'}_{R,h,\nn,\nn'}$ (see Definition \ref{def:19}).

\np Items $1$ and $2$ are unchanged.
%\begin{itemize}
\item{$\quad$3$\,'$} One assigns to the line
$\ell_{0}$ the further label $h_{\ell_{0}} \le h$, and
requires $|\nn| \ge 2^{(h_{\ell_{0}}-2)/\tau}$.
\item{$\quad$4$\,'$} One has $ \max_{\ell\in L(\theta)\setminus
\{\ell_e\}}h_\ell=h$
%\end{itemize}

\np Items $5$ and $6$ are unchanged.

\np The set $\SSSS_{R}$ is defined analogously as $\RR_{R}$.
\end{defi}
%%%%%%%%%%%%%%%%%%%%%%%%%%%%%%%%%%%%%%%%%%%%%%%%%%%%%%%%%%%%%%%%%%%%%%%%%%

%%%%%%%%%%%%%%%%%%%%%%%%%%%%%%%%%%%%%%%%%%%%%%%%%%%%%%%%%%%%%%%%%%%%%%%%%%
\begin{rmk} \label{rmk:14}
Note that if $\theta \in \RR_{R,h,\nn,\nn'}^{(k)\s,\s'}$ then
$\Val(\theta)=\Val(\theta')$ with
$\theta' \in \SSSS_{R,h,\nn,\nn'}^{(k)\s,\s'}$ such that
$h_{\ell_{0}}=h-1$. Thus, it is enough to study the set $\SSSS_{R}$
in order to obtain bounds for trees in $\RR_{R}$.
\end{rmk}
%%%%%%%%%%%%%%%%%%%%%%%%%%%%%%%%%%%%%%%%%%%%%%%%%%%%%%%%%%%%%%%%%%%%%%%%%%

%%%%%%%%%%%%%%%%%%%%%%%%%%%%%%%%%%%%%%%%%%%%%%%%%%%%%%%%%%%%%%%%%%%%%%%%%%
\begin{defi} \label{def:24}
\textbf{ (Tree values).}
For any tree or renormalised tree $\theta$ call
\begin{equation}
\Val(\theta) =
\Big( \prod_{\ell\in L(\theta)} g_{\ell} \Big)
\Big( \prod_{v \in V(\theta)} \h_{v} \Big) \nonumber
\end{equation}
the {\rm value} of the tree $\theta$. To make explicit the dependence
of the tree value on $\eps$ and $M$, sometimes we shall write
$\Val(\theta)=\Val(\theta;\eps,M)$.
\end{defi}
%%%%%%%%%%%%%%%%%%%%%%%%%%%%%%%%%%%%%%%%%%%%%%%%%%%%%%%%%%%%%%%%%%%%%%%%%%

%%%%%%%%%%%%%%%%%%%%%%%%%%%%%%%%%%%%%%%%%%%%%%%%%%%%%%%%%%%%%%%%%%%%%%%%%%
\begin{defi} \label{def:25}
\textbf{ (Counterterms).}
We define the node factors $L^{(k)\s,\s'}_{h,\nn,\nn'}$
(cf. item 21 in Definition \ref{def:15}) by setting
\begin{equation}
L^{(k)\s,\s'}_{h,\nn,\nn'} = \sum_{h' < h-1}
\sum_{\theta\in \RR^{(k)\s,\s'}_{R,h',\nn,\nn'}} \Val(\theta) ,
\qquad \s,\s'=\pm ,
\label{eq:4.9} \end{equation}
for all $k\ge N$, all $h\ge 1$, and all resonant pairs $\{\nn,\nn'\}$.
The counterterms $L$ are then expressed in terms
of (\ref{eq:4.9}) through (\ref{eq:4.1}) and (\ref{4.6}).
\end{defi}
%%%%%%%%%%%%%%%%%%%%%%%%%%%%%%%%%%%%%%%%%%%%%%%%%%%%%%%%%%%%%%%%%%%%%%%%%%

%%%%%%%%%%%%%%%%%%%%%%%%%%%%%%%%%%%%%%%%%%%%%%%%%%%%%%%%%%%%%%%%%%%%%%%%%%
\begin{lemma} \label{lem:5}
For any tree $\theta\in \RR_{R,h,\nn,\nn'}^{(k)\s,\s'}$
there exists a tree $\theta'\in \RR_{R,h,\nn',\nn}^{(k)-\s',-\s}$
such that $\Val(\theta)=\Val(\theta')$.
\end{lemma}
%%%%%%%%%%%%%%%%%%%%%%%%%%%%%%%%%%%%%%%%%%%%%%%%%%%%%%%%%%%%%%%%%%%%%%%%%%

%%%%%%%%%%%%%%%%%%%%%%%%%%%%%%%%%%%%%%%%%%%%%%%%%%%%%%%%%%%%%%%%%%%%%%%%%%
\prova Given a tree $\theta \in \RR_{R,h,\nn,\nn'}^{(k)\s,\s'}$,
consider the path $\PPP=\PPP(\ell_{e},\ell_{0})$, and set
$\PPP=\{\ell_{1},\ldots,\ell_{N}\}$, with $\ell_{0} \succ
\ell_{1} \succ \ldots \succ \ell_{N} \succ \ell_{N+1}=\ell_{e}$
(if $\PPP=\emptyset$, set $N=0$ in the forthcoming discussion).
For $k=0,\ldots,N$, denote by $v_{k}$ the node which the line
$\ell_{k}$ exits and by $L_{0}(v_{k})$ the set $L(v_{k})
\setminus\{\ell_{k+1}\}$
(cf. item 1 in Definition \ref{def:15}).

We construct a tree $\theta' \in \RR_{R,h,\nn',\nn}^{(k)-\s',-\s}$
in the following way.
\begin{enumerate}
\item We shift the sign labels down the path $\PPP$ and change their
sign, so that $\s_{\ell_{k}} \to - \s_{v_{k}}$
and $\s_{v_{k}} \to - \s_{\ell_{k+1}}$ for $k=0,\ldots,N$.
In particular $\ell_{0}$ acquires the label $-\s_{v_{0}}$, while
$\ell_{e}$ loses its label $\s_{\ell_{e}}$ (which with the opposite sign
becomes associated with the node $v_{N}$).
\item The end-node $e$ becomes the root, and the root line
becomes the end-node $e$. In particular
the line $\ell_{e}$ becomes the root line,
and the line $\ell_{0}$ becomes the entering line,
so that the arrows of all the lines $\ell\in\PPP$ are reverted,
while the ordering of all the lines and nodes outside
$\PPP$ is not changed.
\item For all the lines $\ell\in\PPP$ we exchange the labels
$\nn_{\ell},\nn'_{\ell}$, so that $\nn_{\ell_{k}} \to \nn_{\ell_{k}}'$
and $\nn_{\ell_{k}}' \to \nn_{\ell_{k}}$ for $k=1,\ldots,N$,
and we set $\nn_{\ell_{e}}'=\nn'$ and $\nn_{\ell_{0}} = \nn$.
\item For all $k=0,\ldots,N$ we replace
$m_{v_{k}} \to -\s_{v_{k}}\s_{\ell_{k+1}} m_{v_{k}}$.
\end{enumerate}
By construction, the tree $\theta'$ belongs to
$\RR_{R,h,\nn',\nn}^{(k)-\s',-\s}$, and all line propagators
and node factors of the lines and nodes,
respectively, which do not belong to $\PPP$ remain the same. 

Moreover, the line propagator of each $\ell_{k} \in \PPP$ in $\theta'$
is $(G_{i_{\ell_{k}},h_{\ell_{k}}})^{-\s_{v_{k}},-\s_{\ell_{k}}}_{
\nn_{\ell},\nn_{\ell_{k}'}}=
(G_{i_{\ell_{k}},h_{\ell_{k}}})^{\s_{\ell_{k}},\s_{v_{k}}}_{
\nn_{\ell_{k}'},\nn_{\ell}}$, hence it does not change with
respect with the line propagator of the corresponding line in $\theta$.
For each node $v_{k}$, the conservation law
\begin{equation}
\nn_{\ell_{k+1}} = (0,- \s_{v_{k}} \s_{\ell_{k+1}} m_{v_{k}})
- \s_{\ell_{k+1}} \Big( -\s_{v_{k}} \nn_{\ell_{k}}' +
\sum_{\ell' \in L_{0}(v_{k})} \s_{\ell'} \nn_{\ell'} \Big)
\label{eq:4.10} \end{equation}
is assured by the conservation law (cf. item \ref{itemconservation}
in Definition \ref{def:15})
\begin{equation}
\nn_{\ell_{k}}' = (0,m_{v_{k}}) +
\s_{v_{k}} \Big( \s_{\ell_{k+1}} \nn_{\ell_{k+1}} +
\sum_{\ell' \in L_{0}(v_{k})} \s_{\ell'} \nn_{\ell'} \Big)
\label{eq:4.11} \end{equation}
for the corresponding node $v_{k}$ in $\theta$: simply
multiply (\ref{eq:4.11}) times $\s_{v_{k}} \s_{\ell_{k+1}}$
in order to obtain (\ref{eq:4.10}).

Finally we want to show that the product of the combinatorial factors
times the node factors of the nodes $v_{0},\ldots,v_{N}$ do not change.
Take a node $v=v_{k}$, for $k=0,\ldots,N$, and call $r_{v}'$ and $s_{v}'$
the number of lines $\ell'\in L_{0}(v)$ with
$\s_{\ell'}=\s_{v}$ and $\s_{\ell'}=-\s_{v}$, respectively.
Set $\s_{v}=\s$ and $\s_{\ell_{k+1}}=\s'$.

Consider first the case $\s'=\s$. In that case in $\theta$ one has
$r_{v}=r_{v}'+1$ and $s_{v}=s_{v}'$, and the combinatorial factor
contains a factor $r_{v}$ because there are $r_{v}$ lines
$\ell$ entering $v$ with $\s_{\ell}=\s$. In $\theta'$
one has $\s_{v} \to -\s$, $r_{v} \to s_{v}'+1$, $s_{v} \to
r_{v}'$ and $m_{v} \to - m_{v}$ (because $\s\s'=1$).
Moreover the corresponding combinatorial factor contains a factor
$(s_{v}+1)$ because there are $s_{v}+1$ lines $\ell$
entering $v$ with $\s_{\ell}=-\s$. Therefore, taking into
account also the combinatorics, the node factor associated with
the node $v$ in $\theta$ is
$(s_{v}+1)a^{-\s}_{s_{v}+1,r_{v}-1,-m_{v}} =
r_{v}\,a^{\s}_{r_{v},s_{v},m_{v}}$, i.e. the same as
in $\theta$, by the condition (\ref{eq:1.11}).

Now, we pass to the case $\s=-\s'$. In that case in $\theta$
one has $r_{v}=r_{v}'$, $s_{v}=s_{v}'+1$. In $\theta'$
one has the same values for $r_{v}$, $s_{v}$ and $\s_{v}$,
so that, by using also that $-\s\s'm_{v}=m_{v}$ in such a case,
the node factors $a^{\s_{v}}_{r_{v},s_{v},m_{v}}$ do not change.
Of course the combinatorial factors do not change either.

In conclusion, one has $\Val(\theta)=\Val(\theta')$,
which yields the assertion.\EP
%%%%%%%%%%%%%%%%%%%%%%%%%%%%%%%%%%%%%%%%%%%%%%%%%%%%%%%%%%%%%%%%%%%%%%%%%%

%%%%%%%%%%%%%%%%%%%%%%%%%%%%%%%%%%%%%%%%%%%%%%%%%%%%%%%%%%%%%%%%%%%%%%%%%%
\begin{rmk} \label{rmk:15}
By Lemma \ref{lem:5} we have that the matrix $L^{(k)}_h$
is self-adjoint, and the Definition \ref{def:25} together
with (\ref{4.6}) implies that we can write
\begin{equation}
L^{(k)\s,\s'}_{\nn,\nn'} = \sum_{h=-1}^{\io} C_{h}(x_{\nn}(\eps)) \!\!\!\!
\sum_{\theta\in \RR^{(k)\s,\s'}_{R,h,\nn,\nn'}} \!\!\!\!
\Val(\theta) , \qquad
C_{h}(x) = \sum_{h'=h+2}^{\io} \chi_{h}(x) ,
\qquad \s=\pm , \nonumber
\end{equation}
for all $k\ge N$, all $h\ge 1$, and all resonant pairs $\{\nn,\nn'\}$.
By construction $x_{\nn}(\eps)=x_{\nn'}(\eps)$
whenever $L^{(k)\s,\s'}_{h,\nn,\nn'}\neq 0$, so that also
$L^{(k)}$ is self-adjoint.
Finally we have that $L^{(k)} =
\widehat \chi_{1} L^{(k)} \widehat \chi_{1}$ 
(cf. the definition of the line propagators $g_{\ell_{0}}$ and
$g_{\ell_{e}}$ for trees $\theta\in \RR^{(k)\s,\s'}_{R,h,\nn,\nn'}$ in
Definition \ref{def:19}).
\end{rmk}
%%%%%%%%%%%%%%%%%%%%%%%%%%%%%%%%%%%%%%%%%%%%%%%%%%%%%%%%%%%%%%%%%%%%%%%%%%

%%%%%%%%%%%%%%%%%%%%%%%%%%%%%%%%%%%%%%%%%%%%%%%%%%%%%%%%%%%%%%%%%%%%%%%%%%
\begin{lemma} \label{lem:6}
One has
\begin{equation}
u^{(k)\s}_{\nn} = \sum_{\theta\in\Theta_{R,\nn}^{(k)\s}} \Val(\theta) ,
\qquad \s = \pm ,
\label{eq:4.12} \end{equation}
for all $k\ge 1$ and all $\nn\in\ZZZ^{D+1}$.
\end{lemma}
%%%%%%%%%%%%%%%%%%%%%%%%%%%%%%%%%%%%%%%%%%%%%%%%%%%%%%%%%%%%%%%%%%%%%%%%%%

%%%%%%%%%%%%%%%%%%%%%%%%%%%%%%%%%%%%%%%%%%%%%%%%%%%%%%%%%%%%%%%%%%%%%%%%%%
\prova For any given counterterm $L$, the coefficients $u^{(k)\s}_{\nn}$
can be written as sums over tree values
\begin{equation}
u^{(k)\s}_{\nn} = \sum_{\theta\in\Theta_{\nn}^{(k)\s}} \Val(\theta) .
\nonumber \end{equation}
This can be easily proved by induction, using the diagrammatic rules
and definitions given in this section; we refer to Lemma 3.6
of \cite{GP2} for details. Then, defining the counterterms
according to Definition \ref{def:25}, all contributions
arising from trees belonging to the set $\Theta_{\nn}^{(k)\s}$
but not to the set $\Theta_{R,\nn}^{(k)\s}$ cancel out exactly
-- see Lemma 3.13 of \cite{GP2} for further details --
and hence the assertion follows. \EP
%%%%%%%%%%%%%%%%%%%%%%%%%%%%%%%%%%%%%%%%%%%%%%%%%%%%%%%%%%%%%%%%%%%%%%%%%%

%%%%%%%%%%%%%%%%%%%%%%%%%%%%%%%%%%%%%%%%%%%%%%%%%%%%%%%%%%%%%%%%%%%%%%%%%
%%%%%%%%%%%%%%%%%%%%%%%%%%%%%%%%%%%%%%%%%%%%%%%%%%%%%%%%%%%%%%%%%%%%%%%%%
\zerarcounters
\section{Bryuno lemma and bounds}
\label{sec:5}
%%%%%%%%%%%%%%%%%%%%%%%%%%%%%%%%%%%%%%%%%%%%%%%%%%%%%%%%%%%%%%%%%%%%%%%%%
%%%%%%%%%%%%%%%%%%%%%%%%%%%%%%%%%%%%%%%%%%%%%%%%%%%%%%%%%%%%%%%%%%%%%%%%%

Given a tree $\theta\in\Theta_{R}$, call $\gotS(\theta,\g)$
the set of  $(\eps,M)\in \DDD_{0}$ such that for all
$\ell\in L_{p}(\theta)$ with $i_{\ell}=1$ one has
\begin{equation}
\begin{cases}
2^{-h_{\ell}-1} \g \le |x_{\nn_{\ell}}(\eps)| \le 2^{-h_{\ell}+1}\g ,
& h_{\ell} \neq -1 , \\
|x_{\nn_{\ell}}(\eps)| \ge \g ,
& h_{\ell}= -1 , \end{cases}
\label{eq:5.1} \end{equation}
and for all $\ell\in L_{p}(\theta)$ one has
\begin{equation}
\begin{cases}
|\de_{\nn_{\ell}}(\eps)| \le \bar\g,
\quad|\de_{\nn'_{\ell}}(\eps)| \le \bar\g ,
& i_{\ell}=1 , \\
\bar\g \le |\de_{\nn_\ell}(\eps)| ,
& i_{\ell}= 0.
\end{cases}
\label{eq:5.2} \end{equation}

Define also $\DDD(\theta,\g)\subset \DDD_{0}$ as the set
of $(\eps,M) \in\DDD_{0}$ such that for all
$\ell\in L_{p}(\theta)$ with $i_{\ell}=0$ one has
$ |\de_{\nn_\ell}(\eps)\pm \bar\g| \ge \g/|\nn_{\ell}|^{\tau_{1}}$,
while for all $\ell\in L_{p}(\theta)$ with $i_{\ell}=1$ one has
\begin{equation}
x_{\nn_{\ell}}(\eps) \ge \frac{\g}{p_{\nn_{\ell}}^{\tau}(\eps) } ,
\qquad |\de_{\nn}(\eps)\pm \bar\g| \ge
\frac{\g}{|\nn|^{\tau_{1}} }\,\;\forall\nn\in
\CC_{\nn_\ell} \cup \CC_{\nn_{\ell}'} ,
\label{eq:5.3} \end{equation}
for some $\tau,\tau_{1}>0$. Note that the second condition
in (\ref{eq:5.3}) does not depend on $M$.

Analogously, given a tree $\theta\in\SSSS_{R}$, we
call $\widetilde\gotS(\theta,\g)$ the set of $(\eps,M)\in \DDD_{0}$
such that (\ref{eq:5.1}) holds for all $\ell\in L_{p}(\theta)
\setminus \{\ell_{e},\ell_{0}\}$ with $i_{\ell}=1$ and (\ref{eq:5.2})
holds for all $\ell\in L_{p}(\theta)$, and we call
$\widetilde\DDD(\theta,\g)$ as the set of $(\eps,M) \in \DDD_{0}$
such that (\ref{eq:5.3}) holds for all $\ell\in L_{p}(\theta)
\setminus \{\ell_{e},\ell_{0}\}$ with $i_{\ell}=1$, while for all
$\ell\in L_{p}(\theta)$ with $i_{\ell}=0$ one has
$|\de_{\nn_\ell}(\eps) \pm \bar\g| \ge \g/|\nn_{\ell}|^{\tau_{1}}$.

%%%%%%%%%%%%%%%%%%%%%%%%%%%%%%%%%%%%%%%%%%%%%%%%%%%%%%%%%%%%%%%%%%%%%%%%%%
\begin{rmk} \label{rmk:16}
If  $(\eps,M)\in\gotS(\theta,\g)$
then $\Val(\theta;\eps,M) \neq 0$, while $(\eps,M)\in\DDD(\theta,\g)$
means that we can use the bounds (\ref{eq:5.3}) to estimate
$\Val(\theta;\eps,M)$. Analogous considerations hold for
trees $\theta\in\SSSS_{R}$.
\end{rmk}
%%%%%%%%%%%%%%%%%%%%%%%%%%%%%%%%%%%%%%%%%%%%%%%%%%%%%%%%%%%%%%%%%%%%%%%%%%

%%%%%%%%%%%%%%%%%%%%%%%%%%%%%%%%%%%%%%%%%%%%%%%%%%%%%%%%%%%%%%%%%%%%%%%%%%
\begin{rmk} \label{rmk:17}
If for some $\eps$ one has $\Val(\theta;\eps,M) \neq 0$
and for two comparable lines $\ell,\ell'\in L(\theta)$
the pair $\{\nn_{\ell}',\nn_{\ell'}\}$ is resonant,
then all the set $\{\nn_{\ell},\nn_{\ell}',\nn_{\ell'},\nn_{\ell'}'\}$
is resonant. This motivates the condition in
item \ref{itemscaledifference} in Definition \ref{def:15}.
\end{rmk}
%%%%%%%%%%%%%%%%%%%%%%%%%%%%%%%%%%%%%%%%%%%%%%%%%%%%%%%%%%%%%%%%%%%%%%%%%%

%%%%%%%%%%%%%%%%%%%%%%%%%%%%%%%%%%%%%%%%%%%%%%%%%%%%%%%%%%%%%%%%%%%%%%%%%%
\begin{rmk} \label{rmk:18}
If $\theta\in\RR^{(k)\s,\s'}_{R,h,\nn,\nn'}$ is such that
$\Val(\theta;\eps,M) \neq 0$, then $\nn,\nn'\in \Delta_j(\eps)$
for some $j$, so that $p_{\nn}(\eps)=p_{\nn'}(\eps)$ and
$|\nn-\nn'| \le C_{1} C_{2} p_{\nn}^{\al+\beta}(\eps) \le
C_{1} C_{2} p_{\nn}^{2\al}(\eps)$.
Moreover $p_{\nn}(\eps) \le |\nn|,|\nn'| \le 2p_{\nn}(\eps)$.
Such properties follow from Hypothesis \ref{hp:3} -- cf.
also Remark \ref{rmk:6}.
\end{rmk}
%%%%%%%%%%%%%%%%%%%%%%%%%%%%%%%%%%%%%%%%%%%%%%%%%%%%%%%%%%%%%%%%%%%%%%%%%%

%%%%%%%%%%%%%%%%%%%%%%%%%%%%%%%%%%%%%%%%%%%%%%%%%%%%%%%%%%%%%%%%%%%%%%%%%%
\begin{defi} \label{def:26}
\textbf{ (The quantity $\boldsymbol N_{\boldsymbol h}\boldsymbol(
\boldsymbol\theta \boldsymbol)$).}
Define $N_{h}(\theta)$ as the set of lines
$\ell\in L(\theta)$
with $i_{\ell}=1$ and scale $h_{\ell}\ge h$.
\end{defi}
%%%%%%%%%%%%%%%%%%%%%%%%%%%%%%%%%%%%%%%%%%%%%%%%%%%%%%%%%%%%%%%%%%%%%%%%%%

%%%%%%%%%%%%%%%%%%%%%%%%%%%%%%%%%%%%%%%%%%%%%%%%%%%%%%%%%%%%%%%%%%%%%%%%%%
\begin{defi} \label{def:27}
\textbf{ (The quantity $\boldsymbol K \boldsymbol(
\boldsymbol\theta \boldsymbol)$).}
Define
$$ K(\theta) = k(\theta) + \sum_{v\in V_{0}(\theta)} |m_{v}| +
\sum_{\ell\in L_{q}(\theta)} |\nn_{\ell}-\nn_{\ell}'| +
\sum_{v\in  E(\theta)} |\nn_{v}| , $$
where $k(\theta)$ is the order of $\theta$.
\end{defi}
%%%%%%%%%%%%%%%%%%%%%%%%%%%%%%%%%%%%%%%%%%%%%%%%%%%%%%%%%%%%%%%%%%%%%%%%%%]

%%%%%%%%%%%%%%%%%%%%%%%%%%%%%%%%%%%%%%%%%%%%%%%%%%%%%%%%%%%%%%%%%%%%%%%%%%
\begin{lemma} \label{lem:7}
There exists a constant $B$ such the following holds.
\enumerate
\item \label{lem:71}
For all $\theta\in\Theta_{R}$ and all lines $\ell\in L(\theta)$
one has $|\nn_{\ell}| \le B (K(\theta))^{1+4\al}$.
\item \label{lem:72}
If $\theta\in \SSSS_{R}$,
for all lines $\ell\in L(\theta)\setminus (\PPP(\ell_{e},\ell_{0})
\cup \{\ell_{0},\ell_{e}\})$
one has $|\nn_{\ell}| \le B (K(\theta))^{1+4\al}$,
while for all lines $\ell\in\PPP(\ell_{e},\ell_{0})\cup\{\ell_{0}\}$
one has $|\nn_{\ell}'| \le B (|\nn_{\ell_{e}}|+K(\theta))^{1+4\al}$.
\item \label{lem:73}
Given a tree $\theta$ let $\ell,\ell'\in L(\theta)$ be two
comparable lines, with $\ell\prec\ell'$,
such that $i_{\ell}=i_{\ell'}=1$ and $i_{\ell''}=0$ for all
the lines $\ell''\in \PPP(\ell,\ell')$.
If $|\nn_{\ell}'-\nn_{\ell'}| \ge B K(\theta)^{1+4\alpha}$,
then one has $\Val(\theta)=0$ for all $\eps$. 
\item \label{lem:74}
If $\theta\in \SSSS_{R}$, $\ell\in\PPP(\ell_{e},\ell_{0})\cup
\{\ell_{0}\}$ and, moreover, $i_{\ell'}=0$ for all lines
$\ell'\in\PPP(\ell_{e},\ell)$, then $|\nn_{\ell}'|$ $\le$
$|\nn_{\ell_{e}}|+B (K(\theta))^{1+4\al}$.
\end{lemma}
%%%%%%%%%%%%%%%%%%%%%%%%%%%%%%%%%%%%%%%%%%%%%%%%%%%%%%%%%%%%%%%%%%%%%%%%%%

%%%%%%%%%%%%%%%%%%%%%%%%%%%%%%%%%%%%%%%%%%%%%%%%%%%%%%%%%%%%%%%%%%%%%%%%%%
\prova Let us consider first trees $\theta\in\Theta_{R}$.
The proof is by induction on the order of the tree $k=k(\theta)$.
For $k=1$ the bound is trivial. If the root line $\ell_{0}$ is either
a $q$-line or an $r$-line or a $p$-line with $i_{\ell_{0}}=0$, again
the bound follows trivially from the inductive bound. If $\ell_{0}$
is a $p$-line with $i_{\ell_{0}}=1$, call $v_{0}$ the node such
that $\ell_{0}=\ell_{v_{0}}$ and $\theta_{1},\ldots,\theta_{s}$
the subtrees with root in $v_{0}$. By the inductive hypothesis
and Hypothesis \ref{hp:3} one obtains, for a suitable constant $C$
and taking $B$ large enough,
$ |\nn_{\ell}| \le |m_{v_{0}}| + B \left(
K(\theta) - 1 - |m_{v_{0}}| \right)^{1+4\al} + C
\left( |m_{v_{0}}| + B( K(\theta) - 1 - |m_{v_{0}}|)
\right)^{2\al(1+4\al)} \le B (K(\theta))^{1+4\al}$,
which proves the assertion for $\Theta_{R}$ in item \ref{lem:71}.

As a byproduct also the bound for $\SSSS_{R}$
is obtained, as far as lines $\ell\notin\PPP(\ell_{e},\ell_{0})
\cup \{\ell_{0},\ell_{e}\}$ are concerned.
The bound $|\nn_{\ell}'| \le B (|\nn_{\ell_{e}}|+K(\theta))^{1+4\al}$
for the lines $\ell\in\PPP(\ell_{e},\ell_{0})
\cup\{\ell_{0}\}$ can be proved similarly by induction.
Thus, also item \ref{lem:72} is proved.

Given two comparable lines $\ell,\ell'$ such that
$i_{\ell''}=0$ for all lines $\ell''\in\PPP(\ell,\ell')$,
then by momentum conservation one has $\min\{|\nn_{\ell}'-\nn_{\ell'}|,
|\nn_{\ell}'+\nn_{\ell'}|\}\le B(K(\theta))^{1+4\alpha}$ in case (I) and
$|\nn_{\ell}'-\nn_{\ell'}| \le B(K(\theta))^{1+4\alpha}$ in case (II).
This proves the bounds in item \ref{lem:73} in case (II)
and in item \ref{lem:74} for both cases (I) and (II).

In case (I), if $i_{\ell}=i_{\ell'}=1$ and
$\max\{|\delta_{\nn_{\ell}'}(\eps)|,|\delta_{\nn_{\ell'}}(\eps)|
\}< 1/2$, then $|\nn_{\ell}'-\nn_{\ell'}| \le |\nn_{\ell}'+\nn_{\ell'}|$
by item 5 in Hypothesis \ref{hp:1}. On the other hand
if $i_{\ell}=i_{\ell'}=1$ and $\max\{|\delta_{\nn_{\ell}'}(\eps)|,
|\delta_{\nn_{\ell'}}(\eps)|\} \ge 1/2$, one has
$\Val(\theta;\eps,M)=0$. Hence item \ref{lem:73} follows
also in case (I). \EP
%%%%%%%%%%%%%%%%%%%%%%%%%%%%%%%%%%%%%%%%%%%%%%%%%%%%%%%%%%%%%%%%%%%%%%%%%%

%%%%%%%%%%%%%%%%%%%%%%%%%%%%%%%%%%%%%%%%%%%%%%%%%%%%%%%%%%%%%%%%%%%%%%%%%%
\begin{lemma} \label{lem:8}
Given a tree $\theta \in \Theta_{R}$ such that
$\DDD(\theta,\g)\cap \gotS(\theta,\g)\neq \emptyset$,
for all $h\ge1$ one has
$$ N_{h}(\theta) \le  \max\{ 0 , c \, K(\theta)
2^{(2-h)\beta/2\tau} - 1 \} , $$
where $c$ is a suitable constant.
\end{lemma}
%%%%%%%%%%%%%%%%%%%%%%%%%%%%%%%%%%%%%%%%%%%%%%%%%%%%%%%%%%%%%%%%%%%%%%%%%%

%%%%%%%%%%%%%%%%%%%%%%%%%%%%%%%%%%%%%%%%%%%%%%%%%%%%%%%%%%%%%%%%%%%%%%%%%%
\prova Define $E_{h}:=c^{-1}2^{(h-2)\beta/2\tau}$. So, we have to prove
that $N_{h}(\theta)\le \max\{0,K(\theta)E_{h}^{-1}-1\}$.

If a line $\ell$ is on scale $h\ge0$ then $\g/p_{\nn_{\ell}}^{\tau}
(\eps) < x_{\nn_{\ell}}(\eps)\le 2^{-h+1}\g$ by (\ref{eq:5.1})
and (\ref{eq:5.3}). Hence $B(K(\theta))^{2} \ge
B (K(\theta))^{1+4\al} \ge |\nn_{\ell}| \ge p_{\nn_{\ell}}(\eps) >
2^{(h-1)/\tau}$, by Lemma \ref{lem:7}, so that
$K(\theta) E_{h}^{-1} \ge c B^{-1/2} 2^{(h-1)/2\tau}
2^{(2-h)\beta/2\tau} \ge 2$ for $c$ suitably large.
Therefore if a tree $\theta$ contains a line $\ell$ on scale $h$
one has $\max\{0,K(\theta)E_{h}^{-1}-1\}=K(\theta)E_{h}^{-1}-1\ge 1$.

The bound $N_{h}(\theta)\le \max\{0,K(\theta)E_{h}^{-1}-1\}$ will
be proved by induction on the order of the tree. Let $\ell_{0}$ be
the root line of $\theta$ and call $\theta_{1},\ldots,\theta_{m}$
the subtrees of $\theta$ whose root lines $\ell_{1},\ldots,\ell_{m}$
are the lines on scale $h_{\ell_{i}}\ge h-1$ and $i_{\ell_{i}}=1$
which are the closest to $\ell_{0}$.

If $h_{\ell_{0}}<h$  we can write 
$N_{h}(\theta)=N_{h}(\theta_{1})+
\ldots+N_{h}(\theta_{m})$, and the bound follows by induction.
If $h_{\ell_{0}}\ge h$
then $\ell_{1},\ldots,\ell_{m}$ are the
entering lines of a cluster $T$ with exiting line $\ell_{0}$;
in that case we have $N_{h}(\theta)=1+N_{h}(\theta_{1})+
\ldots+N_{h}(\theta_{m})$. Again the bound follows by induction
for $m=0$ and $m\ge2$. The case $m=1$ can be dealt with as follows.

If $\{\nn_{\ell_{0}}',\nn_{\ell_{1}}\}$ is a resonant pair, then
either there exists a line $\ell\in \PPP(\ell_{1},\ell_{0})$
with $i_{\ell}=1$ such that $\{\nn_{\ell}',\nn_{\ell_{1}}\}$
is a resonant pair or there must be a line $\ell\in L(T)\setminus
\PPP(\ell_{1},\ell_{0})$ with $\{\nn_{\ell}',\nn_{\ell_{1}}\}$
a resonant pair. In fact, the first case is not possible: indeed,
also $\{\nn_{\ell_{0}}',\nn_{\ell}'\}$ would be resonant
(cf. Remark \ref{rmk:17}), so that $|h_{\ell}-h_{\ell_{0}}| \le 1$
(cf. item \ref{itemscaledifference} in Definition \ref{def:15}),
and hence the contradiction $h-2 \ge h_{\ell} \ge h_{\ell_{0}}-1
\ge h-1$ would follow. In the second case,
one has $|\nn_{\ell}'| \ge p_{\nn_{\ell_{1}}}(\eps) > 2^{(h-2)/\tau}$,
hence if $\theta'$ is the subtree with root line $\ell$, then one has
$K(\theta)-K(\theta_{1}) > K(\theta') > 2 E_{h}$, and the
bound follows once more by the inductive hypothesis.

If $\{\nn_{\ell_{0}}',\nn_{\ell_{1}}\}$ is not a resonant pair,
call $\bar\ell$ the line along the path $\PPP(\ell_{1},\ell_{0})
\cup\{\ell_{1}\}$ with $i_{\bar\ell} =1$ closest to $\ell_{0}$.
Since $i_{\bar\ell} =1 $ and by hypothesis $h_{\bar\ell} <h-1$ then
$\{\nn_{\bar\ell}, \nn_{\ell_{0}}\}$ is not a resonant pair
(see item \ref{itemscaledifference} in Definition \ref{def:15}).
Call $\tilde T$ the set of nodes and lines preceding $\ell_{0}$
and following $\bar\ell$, and define $K(T)=K(\theta)-K(\theta_{1})$
and $K(\tilde T)=K(\theta)-K(\bar\theta)$, where $\bar\theta$
is the tree with root line $\bar\ell$. Set also 
$\bar\nn=\nn_{\bar\ell}$ and $\nn_{0}=\nn_{\ell_{0}}'$. 
One has $2 |\bar\nn-\nn_{0}| \ge C_{2} (p_{\bar\nn}(\eps)+
p_{\nn_{0}}(\eps))^{\beta}\ge C_{2} p_{\nn_{0}}^{\beta}(\eps)$
(see Remark \ref{rmk:6}), so that by Lemma \ref{lem:7}
one finds $B(K(\theta)-K(\theta_{1})^{2}
\ge B (K(\tilde T))^{2} \ge |\bar\nn-\nn_{0}| \ge
C_{2}\,p_{\nn_{0}}^{\beta}(\eps)/2 \ge C_{2} 2^{(h-1)\beta/\tau}/2$.
Hence $(K(\theta)-K(\theta_{1})) E_{h}^{-1} \ge K(T)E_{h}^{-1}
\ge K(\tilde T)E_{h}^{-1} \ge 2$, provided $c$ is large enough.
This proves the bound. \EP
%%%%%%%%%%%%%%%%%%%%%%%%%%%%%%%%%%%%%%%%%%%%%%%%%%%%%%%%%%%%%%%%%%%%%%%%%%

%%%%%%%%%%%%%%%%%%%%%%%%%%%%%%%%%%%%%%%%%%%%%%%%%%%%%%%%%%%%%%%%%%%%%%%%%%
\begin{lemma} \label{lem:9}
There exists positive constants $\xi_{0}$ and $D_{0}$
such that, if $\xi>\xi_{0}$ in Definition \ref{def:8}, then
for all trees $\theta\in \Theta_{R}$ and for all 
$(\eps,M)\in \DDD(\theta,\g) \cap \gotS(\theta,\g)$ one has
\begin{subequations}
\begin{align}
& \left| \Val(\theta) \right| \le D_{0}^{k} {\rm e}^{-\kappa K(\theta)}
\prod_{\substack{\ell \in L(\theta)\\ i_{\ell} = 1}}
p_{\nn_{\ell}}^{-(\xi-\xi_{0})}(\eps) ,
\label{eq:5.4a} \\
& \left| \pr_{\eps} \Val(\theta) \right| \le
D_{0}^{k} {\rm e}^{-\kappa K(\theta)}
\prod_{\substack{\ell \in L(\theta)\\ i_{\ell} = 1}}
p_{\nn_{\ell}}^{-(\xi-\xi_{0})}(\eps) ,
\label{eq:5.4b} \\
& \sum_{\nn\in\gotO}\sum_{\nn'\in \CC_{\nn}} \sum_{\s,\s'=\pm}
\left|\pr_{M^{\s,\s'}_{\nn,\nn'}}\Val(\theta) \right| \le D_{0}^{k}
{\rm e}^{-\kappa K(\theta)}
\prod_{\substack{\ell \in L(\theta)\\ i_{\ell} = 1}}
p_{\nn_{\ell}}^{-(\xi-\xi_{0})}(\eps) .
\label{eq:5.4c}
\end{align}
\label{eq:5.4}
\end{subequations}
\vskip-.5truecm
\noindent 
\end{lemma}
%%%%%%%%%%%%%%%%%%%%%%%%%%%%%%%%%%%%%%%%%%%%%%%%%%%%%%%%%%%%%%%%%%%%%%%%%%

%%%%%%%%%%%%%%%%%%%%%%%%%%%%%%%%%%%%%%%%%%%%%%%%%%%%%%%%%%%%%%%%%%%%%%%%%%
\prova The propagators are bounded according to (\ref{eq:4.5}),
so that for all trees $\theta\in\Theta^{(k)}_{R,\nn}$ one has
\begin{eqnarray}
& & \left| \Val(\theta) \right| \le
C^{k} \Big( \prod_{v\in V_{0}(\theta)} {\rm e}^{-A_{2}|m_{v}|} \Big)
\Big( \prod_{\ell \in L_{q}(\theta)}
{\rm e}^{-\lambda_{0}|\nn_{\ell}- \nn_{\ell}'|} \Big) \times
\nonumber \\
& & \qquad \qquad \qquad \times
\Big( \prod_{v\in E(\theta)} {\rm e}^{-\lambda_{0}|\nn_{v}|} \Big)
2^{k h_{0}} \Big( \prod_{h=h_{0}+1}^{\io} 2^{h N_{h}(\theta)} \Big)
\prod_{\substack{\ell \in L(\theta)\\ i_{\ell} = 1}}
p_{\nn_{\ell}}^{-\xi}(\eps) \, p_{\nn_{\ell}}^{a_{0}}(\eps) , \nonumber
\end{eqnarray}
for arbitrary $h_{0}$ and for suitable constants $C$ and $a_{0}$.
For $(\eps,M)\in \DDD(\theta,\g) \cap \gotS(\theta,\g)$
one can bound $N_{h}(\theta)$ through Lemma \ref{lem:8}.
Therefore, by choosing $h_{0}$ large enough the bound (\ref{eq:5.4a})
follows, provided $\xi - a_{0}>0$ and $\kappa$ is suitably chosen.

When bounding $\pr_{\eps}\Val(\theta)$, one has to consider
derivatives of the line propagators, i.e. $\pr_{\eps}g_{\ell}$.
If $\ell$ is an $r$-line then $|\pr_{\eps}g_{\ell}|$ is bounded
proportionally to $|\nn_{\ell}|^{c_{0}}$,
whereas if $\ell$ is a $p$-line, then the derivative
produces factors which admit bounds of the form
\begin{equation}
C p_{\nn_{\ell}}^{a_{1}}(\eps) \, 2^{2h_{\ell}}
p_{\nn_{\ell}}^{c_{0}}(\eps)
p_{\nn_{\ell}}^{-\xi}(\eps) , 
\label{eq:5.5} \end{equation}
for suitable constants $C$ and $a_{1}$; see the proof
of Lemma 4.2 in \cite{GP2} for details
(and use item 3 in Hypothesis \ref{hp:1}).

The extra factor $2^{h_{\ell}}$ can be taken into account
by bounding the product of line propagators with
$$ 2^{2h_{0}k} \prod_{h=h_{0}+1}^{\io} 2^{2hN_{h}(\theta)} . $$
One can bound $|\nn_{\ell}| \le B (K(\theta))^{2}$, and use part of
the exponential decaying factors ${\rm e}^{-A_{2}|m_{v}|}$,
${\rm e}^{-\lambda_{0}|\nn_{\ell}-\nn_{\ell}'|}$, and
${\rm e}^{-\lambda_{0}|\nn_{v}|}$, to control the contribution
$\sum_{v\in V_{0}(\theta)}|m_{v}| + 
\sum_{\ell\in L_{q}(\theta)} |\nn_{\ell}-\nn_{\ell}'| +
\sum_{v\in  E(\theta)} |\nn_{v}|$
to $K(\theta)$ (cf. Definition \ref{def:27}).
Then, if $\xi$ is large enough, so that
$\xi-a_{1}>0$ for all possible values of $a_{1}$ in (\ref{eq:5.5}),
the bound (\ref{eq:5.4b}) follows.

Also the bound (\ref{eq:5.4c}) can be discussed in the same way.
We refer again to \cite{GP2} for the details. \EP
%%%%%%%%%%%%%%%%%%%%%%%%%%%%%%%%%%%%%%%%%%%%%%%%%%%%%%%%%%%%%%%%%%%%%%%%%%

%%%%%%%%%%%%%%%%%%%%%%%%%%%%%%%%%%%%%%%%%%%%%%%%%%%%%%%%%%%%%%%%%%%%%%%%%%
\begin{rmk} \label{rmk:19}
Note that for $(\eps,M)\in \DDD(\theta,\g)$ the
singularities of the functions $\bar\chi_{1}$ are avoided,
so that $\pr_{\eps}\bar\chi_{1}(\de_{\nn_{\ell}}(\eps))=0$
for all $\ell\in L(\theta)$. Note also that the bound
(\ref{eq:5.4c}) is not really needed in the following.
\end{rmk}
%%%%%%%%%%%%%%%%%%%%%%%%%%%%%%%%%%%%%%%%%%%%%%%%%%%%%%%%%%%%%%%%%%%%%%%%%%

%%%%%%%%%%%%%%%%%%%%%%%%%%%%%%%%%%%%%%%%%%%%%%%%%%%%%%%%%%%%%%%%%%%%%%%%%%
\begin{lemma} \label{lem:10}
There are two positive constants $B_{2}$ and $B_{3}$
such that the following holds.
\begin{enumerate}
\item \label{lem:10.1}
Given a tree $\theta\in\SSSS_{R}$ such that $\Val(\theta;\eps,M)\neq0$,
if $K(\theta) \le B_{2} p_{\nn_{\ell_{e}}}^{\beta/2}(\eps)$
then for all lines $\ell\in\PPP(\ell_{e},\ell_{0})$ one has $i_\ell=0$.
Moreover for all such lines $\ell$, if $\{\nn_{\ell}',\nn_{\ell_{e}}\}$
is not a resonant pair, then one has $|\delta_{\nn_{\ell}}(\eps)| \ge 1/2$.
\item \label{lem:10.2}
Given a tree $\theta\in \RR_{R}$ such that $\Val(\theta;\eps,M)\neq0$,
one has $\left| \nn_{\ell_{0}}'-\nn_{\ell_{e}} \right| \le
B_{3} (K(\theta))^{1/\rho}$,
with $\rho$ depending on $\al$ and $\beta$.
\end{enumerate}
\end{lemma}
%%%%%%%%%%%%%%%%%%%%%%%%%%%%%%%%%%%%%%%%%%%%%%%%%%%%%%%%%%%%%%%%%%%%%%%%%%

%%%%%%%%%%%%%%%%%%%%%%%%%%%%%%%%%%%%%%%%%%%%%%%%%%%%%%%%%%%%%%%%%%%%%%%%%%
\prova Suppose that $\theta \in \SSSS^{(k)\s,\s'}_{R,h,\nn,\nn'}$ and
$\PPP(\ell_{e},\ell_{0})$ contains lines $\ell$ with $i_{\ell}=1$
and consequently with $\{\nn_{\ell}',\nn'\}$ not resonant
(cf. Definition \ref{def:23}).
Let $\bar\ell$ be the one closest to $\ell_{e}$; thus,
one has $|\nn_{\bar\ell}'-\nn'| \ge C_{3} (|\nn_{\bar\ell}'|+
|\nn'|)^{\beta} \ge C_{3} p_{\nn'}^{\beta}(\eps)=
C_{3} p_{\nn}^{\beta}(\eps)$,
so that we can apply item \ref{lem:73} in Lemma \ref{lem:7} to obtain
$B (K(\theta))^{2} \ge C p_{\nn}^{\beta}(\eps)$,
for some positive constant $C$.
This proves the first statement in item \ref{lem:10.1}.
The proof of the second statement is identical,
since $|\delta_{\nn_\ell}(\eps)| <1/2$
implies that $\nn_{\ell}\in \Delta_{j_{1}}(\eps)$ for some
$j_{1}$, so that if $\{\nn_{\ell}',\nn'\}$ is not a resonant pair
then $\nn'\notin \Delta_{j_{1}}(\eps)$, and therefore $|\nn_{\ell}'-\nn'|
\ge C_{3} p_{\nn'}^{\beta}(\eps)$.

To prove item \ref{lem:10.2}, notice that $|\nn-\nn'| \le C_{1}C_{2}
p_{\nn}^{\al+\beta}(\eps)$ (cf. Remark \ref{rmk:18}).
If $K(\theta) > B_{2} p_{\nn}^{\beta/2}(\eps)$ then
$K(\theta) \ge C|\nn-\nn'|^{\beta/2(\al+\beta)}$.
If $K(\theta) \le  B_{2} p_{\nn}^{\beta/2}(\eps)$ then
$\PPP(\ell_{e},\ell_{0})$ has only lines with $i_\ell=0$, so that by
item \ref{lem:73} in Lemma \ref{lem:7}
one finds $|\nn-\nn'| \le B K(\theta)^{2}$. \EP
%%%%%%%%%%%%%%%%%%%%%%%%%%%%%%%%%%%%%%%%%%%%%%%%%%%%%%%%%%%%%%%%%%%%%%%%%%

%%%%%%%%%%%%%%%%%%%%%%%%%%%%%%%%%%%%%%%%%%%%%%%%%%%%%%%%%%%%%%%%%%%%%%%%%%
\begin{lemma} \label{lem:11}
Given a tree $\theta \in \SSSS_{R}$
such that $\widetilde\DDD(\theta,\g)
\cap\widetilde\gotS(\theta,\g)
\neq \emptyset$, if $N_{h}(\theta) \ge 1$ for some $h \ge 1$,
then $c K(\theta) 2^{(2-h)\beta/2\tau} \ge 1$,
with $c$ the same constant as in Lemma \ref{lem:8}.
\end{lemma}
%%%%%%%%%%%%%%%%%%%%%%%%%%%%%%%%%%%%%%%%%%%%%%%%%%%%%%%%%%%%%%%%%%%%%%%%%%

%%%%%%%%%%%%%%%%%%%%%%%%%%%%%%%%%%%%%%%%%%%%%%%%%%%%%%%%%%%%%%%%%%%%%%%%%%
\prova Consider a tree $\theta\in\SSSS^{(k)\s,\s'}_{R,\bar h,\nn,\nn'}$
for some $k\ge1$, $\bar h\ge 1$, $\s,\s'=\pm$ and $\nn,\nn'\in\gotO$
such that $|\nn'|$ $\ge$ $2^{(\bar h-2)/\tau}$.
Assume $N_{h}(\theta) \ge 1$ for some $\bar h \ge h \ge 1$.

If there is a line $\ell \in L(\theta)$, which does not belong
to $\PPP:=\PPP(\ell_{e},\ell_{0})$, such that $h_{\ell}\ge h$,
then one can reason as at the beginning of the proof of Lemma \ref{lem:8}
to obtain $K(\theta) E_{h}^{-1} \ge 2$,
with $E_{h}=c^{-1} 2^{(h-2)\beta/2\tau} \ge 1$.

Otherwise, there are lines $\ell\in\PPP$ on scale $h_{\ell}\ge h$,
and hence such that $i_{\ell}=1$ and, consequently,
$\{\nn_{\ell}',\nn'\}$ is not a resonant pair.
Let $\bar\ell$ be the one closest to $\ell_{e}$ among such lines;
thus, one has $|\nn_{\bar\ell}'-\nn'| \ge C_{3} p_{\nn'}^{\beta}(\eps)$,
so that one obtains $B (K(\theta))^{2} \ge C p_{\nn'}^{\beta}(\eps)
\ge C 2^{(\bar h-2)\beta/\tau}$, for some positive constant $C$.
So, the desired bound follows once more. \EP
%%%%%%%%%%%%%%%%%%%%%%%%%%%%%%%%%%%%%%%%%%%%%%%%%%%%%%%%%%%%%%%%%%%%%%%%%%

%%%%%%%%%%%%%%%%%%%%%%%%%%%%%%%%%%%%%%%%%%%%%%%%%%%%%%%%%%%%%%%%%%%%%%%%%%
\begin{lemma} \label{lem:12}
Given a tree $\theta \in \SSSS_{R}$ such that $\widetilde
\DDD(\theta,\g) \cap\widetilde\gotS(\theta,\g) \neq \emptyset$,
for all $h \ge 1$ one has
$$ N_{h}(\theta) \le   c \, K(\theta)
2^{(2-h)\beta/2\tau} , $$
where $c$ is the same constant as in Lemma \ref{lem:8}.
\end{lemma}
%%%%%%%%%%%%%%%%%%%%%%%%%%%%%%%%%%%%%%%%%%%%%%%%%%%%%%%%%%%%%%%%%%%%%%%%%%

%%%%%%%%%%%%%%%%%%%%%%%%%%%%%%%%%%%%%%%%%%%%%%%%%%%%%%%%%%%%%%%%%%%%%%%%%%
\prova Consider a tree $\theta\in\SSSS^{(k)\s,\s'}_{R,\bar h,\nn,\nn'}$
for some $k\ge1$, $\bar h\ge 1$, $\s,\s'=\pm$ and $\nn,\nn'\in \gotO$
such that $|\nn'|$ $\ge$ $2^{(\bar h-2)/\tau}$.

For $k(\theta)=1$ one has $N_{h}(\theta) \le 1$, so that
the bound follows from Lemma \ref{lem:11}.

For $k(\theta)>1$ one can proceed as follows. 
Let $\ell_{0}$ be the root line
of $\theta$ and call $\theta_{1},\ldots,\theta_{m}$
the subtrees of $\theta$ whose root lines $\ell_{1},\ldots,\ell_{m}$
are the lines on scale $h_{\ell_{i}}\ge h-1$ and $i_{\ell_{i}}=1$
which are the closest to $\ell_{0}$.
All the trees $\theta_{i}$ such that
$\ell_{i} \notin\PPP(\ell_{e},\ell_{0})$ belong to some
$\Theta_{R,\nn_i}^{(k_{i})\pm}$ with $k_i<k$.
If  $K(\theta)\ge B_2p_{\nn'}^{\beta/2}(\eps)$
(cf. Lemma \ref{lem:10}) it may be possible that
a line, say $\ell_{1}$, belongs to $\PPP(\ell_e,\ell_0)$,
so that $\Val(\theta_{1})=g_{\ell_{1}}\,\Val(\theta_{1}')$, with
$\theta_{1}' \in \SSSS_{R,h_1,\nn_1,\nn'}^{(k_1),\s_1,\s'}$
with $h_{1}\le \bar h$, $\s_{1}=\pm$ and $k_1<k$.

If $h_{\ell_{0}}<h$ one has
$N_{h}(\theta)= N_{h}(\theta_{1}) + \ldots + N_{h}(\theta_{m})$,
so that the bound $N_{h}(\theta)\le K(\theta)E_{h}^{-1}$
follows by the inductive hypothesis.

If $h_{\ell_{0}} \ge h$ one has $N_{h}(\theta)= 1 +
N_{h}(\theta_{1}) + \ldots + N_{h}(\theta_{m})$.
For $m=0$ the bound can be obtained once more from Lemma \ref{lem:11},
while for $m \ge 2$ at least one tree, say $\theta_{m}$,
belongs to $\Theta_{R,\nn'}^{(k')\pm}$
for some $k'$ and $\nn'$ so that we can apply Lemma \ref{lem:8}
and the inductive hypothesis to obtain
\begin{eqnarray}
& & N_{h}(\theta) \le 1 + \left( K(\theta_{1}) + \ldots +
K(\theta_{m-1}) \right) E_{h}^{-1} +
\left( K(\theta_{m}) E_{h}^{-1} - 1 \right) \nonumber \\
& & \hskip1.0truecm
\le \left( K(\theta_{1}) + \ldots + K(\theta_{m-1}) \right) E_{h}^{-1}
+ K(\theta_{m}) E_{h}^{-1} \le
K(\theta) E_{h}^{-1} , \nonumber
\end{eqnarray}
which yields the bound.

Finally if $m=1$ one has $N_{h}(\theta)=1+N_{h}(\theta_{1})$.
Hence, if $\ell_{1}\notin\PPP(\ell_{e},\ell_{0})$,
again the bound follows from Lemma \ref{lem:8}. If on the contrary
$\ell_{1} \in \PPP(\ell_{e},\ell_{0})$, one can adapt the discussion
of the case $m=1$ in the proof of Lemma \ref{lem:8}. \EP
%%%%%%%%%%%%%%%%%%%%%%%%%%%%%%%%%%%%%%%%%%%%%%%%%%%%%%%%%%%%%%%%%%%%%%%%%%

%%%%%%%%%%%%%%%%%%%%%%%%%%%%%%%%%%%%%%%%%%%%%%%%%%%%%%%%%%%%%%%%%%%%%%%%%%
\begin{lemma} \label{lem:13}
There exists positive constants $\kappa$, $\xi_{1}$ and $D_{1}$
such that, if $\xi>\xi_{1}$ in Definition \ref{def:8},
then for all trees $\theta\in \RR_{R}$ and for all 
$(\eps,M)\in \widetilde\DDD(\theta,\g) \cap \widetilde\gotS(\theta,\g)$,
by setting $\nn=\nn_{\ell_{0}}'$ and $\nn'=\nn_{\ell_{e}}$, one has
\begin{subequations}
\begin{align}
& \left| \Val(\theta) \right| \le
D_{1}^{k} 2^{-h} {\rm e}^{-\kappa|\nn-\nn'|^{\rho}}
\prod_{\substack{\ell \in L(\theta)\\ i_{\ell} = 1}}
p_{\nn_{\ell}}^{-(\xi-\xi_{1})}(\eps) ,
\label{eq:5.6a} \\
& \left| \pr_{\eps} \Val(\theta) \right| \le D_{1}^{k}
2^{-h} p_{\nn}^{c_{0}}(\eps) \, {\rm e}^{-\kappa|\nn-\nn'|^{\rho}}
\prod_{\substack{\ell \in L(\theta)\\ i_{\ell} = 1}}
p_{\nn_{\ell}}^{-(\xi-\xi_{1})}(\eps) ,
\label{eq:5.6b} \\
& \sum_{\nn_{1}\in\gotO}\sum_{\nn_{2}\in\CC_{\nn_{1}}}
\sum_{\s_{1},\s_{2}=\pm}\left|\pr_{M^{\s_{1},\s_{2}}_{\nn_{1},\nn_{2}}}
\Val(\theta) \right|
\le D_{1}^{k} 2^{-h} {\rm e}^{-\kappa|\nn-\nn'|^{\rho}}
\prod_{\substack{\ell \in L(\theta)\\ i_{\ell} = 1}}
p_{\nn_{\ell}}^{-(\xi-\xi_{1})}(\eps) ,
\label{eq:5.6c} 
\end{align}
\label{eq:5.6}
\end{subequations}
\vskip-.3truecm
\noindent with $\rho$ as in Lemma \ref{lem:10}.
\end{lemma}
%%%%%%%%%%%%%%%%%%%%%%%%%%%%%%%%%%%%%%%%%%%%%%%%%%%%%%%%%%%%%%%%%%%%%%%%%%

%%%%%%%%%%%%%%%%%%%%%%%%%%%%%%%%%%%%%%%%%%%%%%%%%%%%%%%%%%%%%%%%%%%%%%%%%%
\prova Set for simplicity $\PPP=\PPP(\ell_{e},\ell_{0})$ and
\begin{eqnarray}
\Sigma(\theta) & \!\!\! = \!\!\! &
\sum_{v\in V_{0}(\theta)} |m_{v}| +
\sum_{\ell\in L_{q}(\theta)} |\nn_{\ell}-\nn_{\ell}'| +
\sum_{v\in  E(\theta)} |\nn_{v}| ,
\nonumber \\
\Pi(\theta) & \!\!\! = \!\!\! &
\Big( \prod_{v\in V_{0}(\theta)} {\rm e}^{A_{2}|m_{v}|/8} \Big)
\Big( \prod_{\ell\in L_{q}(\theta)}
{\rm e}^{\lambda_{0}|\nn_{\ell}-\nn_{\ell}'|}
\Big) \Big( \prod_{v\in  E(\theta)}
{\rm e}^{\lambda_{0}|\nn_{v}|} \Big) .
\nonumber
\end{eqnarray}
If $\theta\in\RR^{(k)\s,\s'}_{R,h,\nn,\nn'}$ for
some $k\ge1$, $h\ge 1$, $\s,\s'=\pm$ and $\{\nn,\nn'\}$ resonant,
then $N_{h}(\theta) \ge 1$,
so that $ K(\theta) = k + \Sigma(\theta) > C 2^{h\beta/2\tau}$,
for some constant $C$, which imply
$1 \le 2^{-h} C^{k} \Pi(\theta)$, for some constant $C$.
This produces the extra factor $2^{-h}$.

By item \ref{lem:10.2} in Lemma \ref{lem:10} one has
$(B_{3}^{-1}|\nn-\nn'|)^{\rho} \le K(\theta)$,
so that $1 \le {\rm e}^{-|\nn-\nn'|^{\rho}}C^{k}
\Pi(\theta)$, for some constant $C$. The factor $\Pi(\theta)$
can be bounded by using part of the factors ${\rm e}^{-A_{2}|m_{v}|}$,
${\rm e}^{-\lambda_{0}|\nn_{v}|}$, and
${\rm e}^{-\lambda_{0}|\nn_{\ell}-\nn_{\ell}'|}$,
associated with the nodes and with the $q$-lines.
This proves the bound (\ref{eq:5.6a}),

To prove the bound (\ref{eq:5.6b}) one has to take into account
the further $\eps$-derivative acting on the line propagator
$g_{\ell}$, for some $\ell\in L(\theta)$. 
If the line $\ell$ does not belong to $\PPP$ then one can
reason as in the proof of (\ref{eq:5.4b}) in Lemma \ref{lem:9}.
If $\ell\in\PPP$ one has to distinguish between two cases.
If there exists a line $\bar\ell\in\PPP$ such that
$i_{\bar\ell}=1$ then $K(\theta) > B_{2} p_{\nn}^{\beta/2}(\eps)$
by item \ref{lem:10.1} in Lemma \ref{lem:10}, so that,
by item \ref{lem:72} in Lemma \ref{lem:7}, one has
$p_{\nn_{\ell}}(\eps) \le |\nn_{\ell}'| \le
B (|\nn_{\ell_{e}}| + K(\theta))^{1+4\al} \le
B (2p_{\nn}(\eps) + K(\theta))^{1+4\al} \le C (K(\theta))^{4/\beta}$,
for some constant $C$. If $i_{\ell}=0$ for all
lines $\ell\in\PPP$ then, by item \ref{lem:73} in
Lemma \ref{lem:7}, one has $p_{\nn_{\ell}}(\eps) \le
|\nn_{\ell}'| \le |\nn_{\ell_{e}}| +
B (K(\theta))^{2}$. Then item 3 in Hypothesis \ref{hp:1}
implies the bound (\ref{eq:5.6b}).

To prove (\ref{eq:5.6c}) one has to study a sum
of terms each containing a derivative
$\partial_{M_{\nn_{1},\nn_{2}}^{\s_{1},\s_{2}}}
g_{\ell}$, for some $\ell\in L(\theta)$. 
If $\ell\in\PPP$ we distinguish between the two cases.
If $K(\theta) > B_{2}p_{\nn}^{\beta/2}(\eps)$, the sum over
$\nn_{1},\nn_{2}$ has the limitations $|\nn_{1}-\nn_{2}|\le
C p_{\nn_{1}}^{\al+\beta}(\eps)$, $|\nn_{1}-\nn_{\ell}|\le
C p_{\nn_{1}}^{\al+\beta}(\eps)$ and $|\nn_{\ell}|
\le (|\nn_{\ell_{e}}| + B K(\theta))^{1+4\al}
\le C(K(\theta))^{4/\beta}$, for some constant $C$:
hence the sum over $\nn_{1},\nn_{2}$ produces a factor
$C (K(\theta))^{C'}$ for suitable constants $C$ and $C'$,
and one has $(K(\theta))^{C'} \le C^{k} \Pi(\theta)$,
for some constant $C$. If $K(\theta) \le B_{2}
p_{\nn}^{\beta/2}(\eps)$, then $i_{\ell}=0$ 
for all lines $\ell\in\PPP$, so that the line propagators
$g_\ell$ do not depend on $M$.
Finally if $\ell\not\in\PPP$ then one has $|\nn_{\ell}| \le
B (K(\theta))^{1+4\al}$, so that the sum over $\nn_{1},\nn_{2}$
is bounded once more proportionally to $(K(\theta))^{C'}$,
for some constant $C'$, and again one can bound $(K(\theta))^{C'}
\le C^{k} \Pi(\theta)$, for some constant $C$. \EP
%%%%%%%%%%%%%%%%%%%%%%%%%%%%%%%%%%%%%%%%%%%%%%%%%%%%%%%%%%%%%%%%%%%%%%%%%%

%%%%%%%%%%%%%%%%%%%%%%%%%%%%%%%%%%%%%%%%%%%%%%%%%%%%%%%%%%%%%%%%%%%%%%%%%%
\begin{rmk} \label{rmk:20}
Both Lemma \ref{lem:12} and \ref{lem:16} deal with the first
derivatives of $\Val(\theta)$. One can easily extend the analysis
so to include derivatives of arbitrary order, at the price
of allowing larger constants $\xi_{1}$ and $D_{1}$ --
and a factor $p_{\nn}^{c_{0}}(\eps)$ for any further $\eps$-derivative.
Therefore, one can prove that the function $\Val(\theta)$ is $C^{r}$
for any integer $r$, in particular for $r=1$, which we shall need
in the following -- cf. in particular the forthcoming
Lemma \ref{lem:14}.
\end{rmk}
%%%%%%%%%%%%%%%%%%%%%%%%%%%%%%%%%%%%%%%%%%%%%%%%%%%%%%%%%%%%%%%%%%%%%%%%%%

%%%%%%%%%%%%%%%%%%%%%%%%%%%%%%%%%%%%%%%%%%%%%%%%%%%%%%%%%%%%%%%%%%%%%%%%%
%%%%%%%%%%%%%%%%%%%%%%%%%%%%%%%%%%%%%%%%%%%%%%%%%%%%%%%%%%%%%%%%%%%%%%%%%
\zerarcounters
\section{Proof of Proposition \ref{prop:1}}
\label{sec:6}
%%%%%%%%%%%%%%%%%%%%%%%%%%%%%%%%%%%%%%%%%%%%%%%%%%%%%%%%%%%%%%%%%%%%%%%%%
%%%%%%%%%%%%%%%%%%%%%%%%%%%%%%%%%%%%%%%%%%%%%%%%%%%%%%%%%%%%%%%%%%%%%%%%%

%%%%%%%%%%%%%%%%%%%%%%%%%%%%%%%%%%%%%%%%%%%%%%%%%%%%%%%%%%%%%%%%%%%%%%%%%%
\begin{defi} \label{def:28}
\textbf{ (The extended tree values).}
Let the function $\chi_{-1}$ be as in Definition \ref{def:11}.
Define
\begin{eqnarray}
& & \Val^{E}(\theta) =
\Big( \prod_{\substack{\ell\in L(\theta) \\ i_{\ell}=1}}
\chi_{-1}(|x_{\nn_{\ell}}(\eps)|\,p_{\nn_{\ell}}^{\tau}(\eps)) \Big)
\Big( \prod_{\substack{\ell\in L(\theta) \\ i_{\ell}=1}}
\prod_{\nn\in \CC_{\{\nn_{\ell},\nn_{\ell}'\}}}
\chi_{-1}(||\delta_{\nn}(\eps)|-\bar\g|\,|\nn|^{\tau_{1}}) \Big)
\times \nonumber \\
& & \qquad \qquad \qquad \times
\Big(\prod_{\substack{\ell\in L_{p}(\theta) \\ i_{\ell}=0}}
\chi_{-1}(||\delta_{\nn_{\ell}}(\eps)|-\bar\g|\,
|\nn_{\ell}|^{\tau_{1}})\Big)
\Val(\theta)
\label{eq:6.1} \end{eqnarray}
for $\theta\in\Theta^{(k)}_{R,\nn}$, and
\begin{eqnarray}
& & \Val^{E}(\theta) =
\Big( \prod_{\substack{\ell\in L(\theta)\setminus
\{\ell_{0},\ell_{e}\} \\ i_{\ell}=1}}
\chi_{-1}(|x_{\nn_{\ell}}(\eps)|\,p_{\nn_{\ell}}^{\tau}(\eps)) \Big)
\Big( \prod_{\substack{\ell\in L(\theta)\setminus
\{\ell_{0},\ell_{e}\} \\ i_{\ell}=1}}
\prod_{\nn\in \CC_{\{\nn_{\ell},\nn_{\ell}'\}}}
\chi_{-1} (||\delta_{\nn}(\eps)|-\bar\g|\,p_{\nn_{\ell}}^{\tau_{1}}(\eps))
\Big) \times \nonumber \\
& & \qquad \qquad \qquad \times
\Big(\prod_{\substack{\ell\in L_{p}(\theta) \\
i_{\ell}=0,\;\nn_{\ell} \notin \CC_{\{\nn,\nn'\}}}}
\chi_{-1}(||\delta_{\nn_\ell}(\eps)|-\bar\g|\,|\nn_\ell|^{\tau_{1}})\Big)
\Val(\theta)
\label{eq:6.2} \end{eqnarray}
for $\theta\in \RR^{(k)}_{R,h,\nn,\nn'}$. We call $\Val^{E}(\theta)$
the {\rm extended value} of the tree $\theta$.
\end{defi}
%%%%%%%%%%%%%%%%%%%%%%%%%%%%%%%%%%%%%%%%%%%%%%%%%%%%%%%%%%%%%%%%%%%%%%%%%%

The following result proves Proposition \ref{prop:1}.

%%%%%%%%%%%%%%%%%%%%%%%%%%%%%%%%%%%%%%%%%%%%%%%%%%%%%%%%%%%%%%%%%%%%%%%%%%
\begin{lemma} \label{lem:14}
Given $\theta\in\RR^{(k)\s,\s'}_{R,h,\nn,\nn'}$, the function $\Val(\theta)$
can be extended to the function (\ref{eq:6.1}) defined
and $C^{1}$ in $\DDD_{0}\setminus \II_{\{\nn,\nn'\}}(\g)$, such that,
defining the ``extended'' counterterm $L^{E\,\s,\s'}_{\nn,\nn'}$ according
to Definition \ref{def:25}, with $\Val(\theta)$
replaced with $\Val^{E}(\theta)$, the following holds.
\begin{enumerate}
\item Possibly with different constants $\xi_{1}$ and $K_{0}$,
$\Val^{E}(\theta)$ satisfies for all $(\eps,M)\in\DDD_{0}\setminus
\II_{\{\nn,\nn'\}}(\g)$ the same bounds in Lemma \ref{lem:13}
as $\Val(\theta)$ in $\DDD(\g)$.
\item There exist constants $\xi_{1}$, $K_{1}$, $\kappa$, $\rho$
and $\eta_{0}$, such that, if $\xi>\xi_{1}$ in Definition \ref{def:8},
$L^{E\,\s,\s'}_{\nn,\nn'}$ satisfies, for all $(\eps,M)\in\DDD_{0}
\setminus\II_{\{\nn,\nn'\}}(\g)$ and $|\eta| \le \eta_{0}$, the bounds
\begin{eqnarray}
& & \left| L^{E\,\s,\s'}_{\nn,\nn'} \right| \le |\eta|^{N}
K_{1} {\rm e}^{-\kappa|\nn-\nn'|^{\rho}} , \qquad
\left| \partial_{\eps} L^{E\,\s,\s'}_{\nn,\nn'} \right| \le |\eta|^{N}
K_{1} p_{\nn}^{c_{0}} {\rm e}^{-\kappa|\nn-\nn'|^{\rho}} , \nonumber \\
& & \left| \partial_{\eta} L^{E\,\s,\s'}_{\nn,\nn'} \right|
\le N\,|\eta|^{N-1} K_{1} {\rm e}^{-\kappa|\nn-\nn'|^{\rho}} , \nonumber \\
& & \sum_{\nn_{1}\in\gotO\,,\s_1=\pm}
\sum_{\nn_{2} \in \CC_{\nn_{1}}\,,\s_2=\pm} \left|
\partial_{M^{\s_1,\s_2}_{\nn_{1},\nn_{2}}} L^{E\, \s,\s'}_{\nn,\nn'} \right|
{\rm e}^{\kappa|\nn-\nn'|^{\rho}}
\le |\eta|^{N} K_{1} . \nonumber
\end{eqnarray}
\item $\Val^{E}(\theta) = \Val(\theta)$ for $(\eps,M)\in\DDD(2\g)$
and $\Val^{E}(\theta)=0$ for $(\eps,M)\in \DDD_{0}\setminus\DDD(\g)$.
\end{enumerate}
Analogously, given $\theta\in\Theta^{(k)\s}_{R,\nn}$,
the function $\Val(\theta)$ can be extended to the function
(\ref{eq:6.2}) defined and $C^{1}$ in $\DDD_{0}$, such that,
defining $u_{\nn}^{E\,(k)}$ as in Lemma \ref{lem:6}
with $\Val(\theta)$ replaced with $\Val^{E}(\theta)$,
the following holds.
\begin{enumerate}
\item Possibly with different constants $\xi_{1}$ and $K_{0}$,
$\Val^{E}(\theta)$ satisfies for all $(\eps,M)\in\DDD_{0}$
the same bounds in Lemma \ref{lem:9} as $\Val(\theta)$ in $\DDD(\g)$.
\item There exist constants $\xi_{1}$, $K_{1}$, $\kappa$ and $\eta_{0}$
such that, if $\xi>\xi_{1}$ in Definition \ref{def:8}, $u_{\nn}^{E\,\s}$
satisfies, for all $(\eps,M) \in \DDD_{0}$ and $|\eta|\le \eta_{0}$,
the bounds
\begin{equation}
\left| u_{\nn}^{E\,\s} \right| \le |\eta|^{N}
K_{1} {\rm e}^{-\kappa|\nn|^{1/2} } \nonumber
\end{equation}
for all $\nn\in\ZZZ^{D+1}$.
\item $\Val^{E}(\theta) = \Val(\theta)$ for $(\eps,M)\in\DDD(2\g)$
and $\Val^{E}(\theta)=0$ for $(\eps,M)\in \DDD_{0}\setminus\DDD(\g)$.
\end{enumerate}
\end{lemma}
%%%%%%%%%%%%%%%%%%%%%%%%%%%%%%%%%%%%%%%%%%%%%%%%%%%%%%%%%%%%%%%%%%%%%%%%%%

%%%%%%%%%%%%%%%%%%%%%%%%%%%%%%%%%%%%%%%%%%%%%%%%%%%%%%%%%%%%%%%%%%%%%%%%%%
\prova We shall consider explicitly the case
of trees $\theta \in \RR^{(k)\s,\s'}_{R,h,\nn,\nn'}$. The case of
trees $\theta \in \Theta^{(k)\s}_{R,\nn}$
can be discussed in the same way.

Item 3 follows from the very definition. The bounds of item 1
can be proved by reasoning as in Section \ref{sec:5}, by taking
into account the further derivatives which arise because of
the compact support functions $\chi_{-1}$ in (\ref{eq:6.2}).
On the other hand all such derivatives produce factors
proportional to $p_{\nn_{\ell}}^{a_{2}}(\eps)$
for some constant $a_{2}$ (again we refer to \cite{GP2} for details);
in particular we are using item 2 in Hypothesis \ref{hp:1} 
to bound the derivatives of $\delta_{\nn_{\ell}}(\eps)$ with respect
to $\eps$. Therefore by using Lemma \ref{lem:8} and 
possibly taking larger constants $\xi_{1}$ and $K_{0}$
the bounds of Lemma \ref{lem:13} follow also for the extended
function (\ref{eq:6.2}).

Finally the bounds on $L^{E}$ in item 2 come directly
from the definition. Indeed, the counterterms $L^{E\,\s,\s'}_{\nn,\nn'}$
are expressed in terms of the values $\Val(\theta)$ according to
Remark \ref{rmk:15}, and the factor $2^{-h}$ is used to perform
the summation over the scale labels. Hence we have to control
the sum over the trees.

Let us fix $\eps$. For each $v\in E(\theta)$ the sum over
$|\nn_{v}|$ is controlled
by using the exponential factors ${\rm e}^{-\lambda_{0}|\nn_{v}|}$.
For each line $\ell\in L(\theta)$ the labels $\nn_{\ell}'$
are fixed by the conservation rule of item 12 in
Definition \ref{def:15}, while the sum over $\nn_{\ell}$
gives a factor $C_{1}p_{\nn_{\ell}}^{\al}(\eps)$ for the $p$-lines
(see item 2 in Hypothesis \ref{hp:3}), and it is controlled 
by using the exponential factors ${\rm e}^{-\lambda_{0}|\nn_{\ell}-
\nn_{\ell}'|}$ for the $q$-lines.
The sums over $i_{\ell}$ and $h_{\ell}$ can be bounded by a factor $4$.
Finally the sum over all the unlabelled trees of order $k$
is bounded by $C^{k}$ for some constant $C$.
Thus, the bounds on $L^{E}_{\nn,\nn'}$ are proved.

Finally, the $C^{1}$ smoothness follows from Remark \ref{rmk:20}.\EP
%%%%%%%%%%%%%%%%%%%%%%%%%%%%%%%%%%%%%%%%%%%%%%%%%%%%%%%%%%%%%%%%%%%%%%%%%%

%%%%%%%%%%%%%%%%%%%%%%%%%%%%%%%%%%%%%%%%%%%%%%%%%%%%%%%%%%%%%%%%%%%%%%%%%
%%%%%%%%%%%%%%%%%%%%%%%%%%%%%%%%%%%%%%%%%%%%%%%%%%%%%%%%%%%%%%%%%%%%%%%%%
\zerarcounters
\section{Proof of Proposition \ref{prop:2}}
\label{sec:7}
%%%%%%%%%%%%%%%%%%%%%%%%%%%%%%%%%%%%%%%%%%%%%%%%%%%%%%%%%%%%%%%%%%%%%%%%%
%%%%%%%%%%%%%%%%%%%%%%%%%%%%%%%%%%%%%%%%%%%%%%%%%%%%%%%%%%%%%%%%%%%%%%%%%

The following result proves item 1 in Proposition \ref{prop:2}.
Here and henceforth we write $L=L(\eta,\eps,M)$ and
$L^{E}=L^{E}(\eta,\eps,M)$, and we fix $\eta=\eps^{1/N}$.

%%%%%%%%%%%%%%%%%%%%%%%%%%%%%%%%%%%%%%%%%%%%%%%%%%%%%%%%%%%%%%%%%%%%%%%%%%
\begin{lemma} \label{lem:15}
There exists constants $\eps_{0}> 0$ such that there exist
functions $M^{\s,\s'}_{\nn,\nn'}(\eps)=M^{\s,\s'}_{\nn',\nn}(\eps)$
well defined and $C^{1}$ for $\eps\in \gotE_{0} \setminus \overline
\II_{\{\nn,\nn'\}}(\g)$, 
such that the ``extended'' compatibility equation
$$ M^{\s,\s'}_{\nn,\nn'}(\eps) =
L^{E\,\s,\s'}_{\nn,\nn'}(\eps^{1/N},\eps, M(\eps)) $$
holds for all $\eps\in(0,\eps_0)\setminus
\overline\II_{\{\nn,\nn'\}}(\g)$.
\end{lemma}
%%%%%%%%%%%%%%%%%%%%%%%%%%%%%%%%%%%%%%%%%%%%%%%%%%%%%%%%%%%%%%%%%%%%%%%%%%

%%%%%%%%%%%%%%%%%%%%%%%%%%%%%%%%%%%%%%%%%%%%%%%%%%%%%%%%%%%%%%%%%%%%%%%%%%
\prova By definition we set $M^{\s,\s'}_{\nn,\nn'}(\eps)=0$ for all $\eps$
such that $\bar\chi_{1}(\delta_{\nn}(\eps))
\bar\chi_{1}(\delta_{\nn'}(\eps))=0$. Consider the Banach space
$\overline\BB$ of lists $\{M^{\s,\s'}_{\nn,\nn'}(\eps)\}$,
with $\{\nn,\nn'\}$ a resonant pair, such that each
$M^{\s,\s'}_{\nn,\nn'}(\eps)$ is well defined and $C^1$
in $\eps\in\gotE_{0}\setminus \overline \II_{\{\nn,\nn'\}}(\g)$
and $M^{\s,\s'}_{\nn,\nn'}(\eps)=0$ for $\eps \in
\overline \II_{\{\nn,\nn'\}}(\g)$. By definition
$\{L^{E\,\s_{1},\s_{2}}_{\nn_{1},\nn_{2}}
(\eps^{1/N},\eps, \{M^{\s,\s'}_{\nn,\nn'}(\eps)\})\}$ is well defined as a
continuously differentiable application from $\overline\BB$ in itself,
since, for each tree $\theta\in\RR^{(k)\s_{1},\s_{2}}_{
R,h,\nn_{1},\nn_{2}}$, the value $\Val^{E}(\theta)$
by definition smoothes out to zero
the value of each line propagator $g_{\ell}$ in the corresponding
intervals $\overline\II_{\{\nn_{\ell},\nn'_{\ell}}\}(2\g)\setminus
\overline\II_{\{\nn_\ell,\nn'_\ell\}}(\g)$. Again by definition
$L^{E}(0,0,0)=0$ and $|\partial_{M} L(0,0,0)|_{\rm op}=0$,
so that we can  apply the implicit function theorem. \EP
%%%%%%%%%%%%%%%%%%%%%%%%%%%%%%%%%%%%%%%%%%%%%%%%%%%%%%%%%%%%%%%%%%%%%%%%%%

Now we pass to the proof of item 2 in Proposition \ref{prop:2}.
We need some preliminary results.

%%%%%%%%%%%%%%%%%%%%%%%%%%%%%%%%%%%%%%%%%%%%%%%%%%%%%%%%%%%%%%%%%%%%%%%%%%
\begin{lemma} \label{lem:16}
Let $A=A(\eps)$ a self-adjoint matrix piecewise differentiable in the
parameter $\eps$. Then, if $\la^{(a)}(A)$ and $\phi^{(a)}(A)$ denote the
eigenvalues and the (normalised) eigenvectors of $A$, respectively,
the following holds.
\begin{enumerate}
\item One has  $|\la^{(a)}(A(\eps))| \le \|A(\eps)\|_{2}$.
\item The eigenvalues $\la^{(a)}(A(\eps))$ are piecewise
differentiable in $\eps$.
\item One has $|\partial_{\eps}\la^{(a)}(A(\eps))| \le
\|\partial_{\eps}A(\eps)\|_{2}$.
\end{enumerate}
\end{lemma}
%%%%%%%%%%%%%%%%%%%%%%%%%%%%%%%%%%%%%%%%%%%%%%%%%%%%%%%%%%%%%%%%%%%%%%%%%%

%%%%%%%%%%%%%%%%%%%%%%%%%%%%%%%%%%%%%%%%%%%%%%%%%%%%%%%%%%%%%%%%%%%%%%%%%%
\prova See \cite{Ka} for items 1 and 2. Moreover, for each
interval in which $A$ is differentiable, let $A_{n}$ be
an analytic approximation of $A$ in such an interval,
with $A_{n}\to A$ as $n\to\io$: then the eigenvalues
$\phi^{(a)}(A_{n})$ are piecewise differentiable \cite{Ka}, and one has
$$ \pr_{\eps} \la^{(a)}(A_{n}) = \pr_{\eps}
\left( \phi^{(a)},A_{n}\phi^{(a)} \right) =
\la^{(a)}(A_{n}) \pr_{\eps} \left( \phi^{(a)},\phi^{(a)} \right) +
\left( \phi^{(a)},\pr_{\eps} A_{n}\phi^{(a)} \right) =
\left( \phi^{(a)},\pr_{\eps} A_{n}\phi^{(a)} \right) , $$
which yields item 3 when the limit $n\to\io$ is taken. \EP
%%%%%%%%%%%%%%%%%%%%%%%%%%%%%%%%%%%%%%%%%%%%%%%%%%%%%%%%%%%%%%%%%%%%%%%%%%

For $M\in\BB_{\kappa}$ we can write $\widehat M = \bigoplus_{j} M_{j}$,
where $M_{j}$ are block matrices, so that we can define
$\|\widehat M\|_{2} = \sup_{j}\|M_{j}\|_{2}$, with
$\|M_{j}\|_{2}$ given as in Definition \ref{def:7}.

%%%%%%%%%%%%%%%%%%%%%%%%%%%%%%%%%%%%%%%%%%%%%%%%%%%%%%%%%%%%%%%%%%%%%%%%%%
\begin{lemma} \label{lem:17}
For $M\in \BB_{\kappa}$ one has $\|\widehat M\|_{2} \le C \eps_{0}$
for some constant $C$ depending on $\kappa$ and $\rho$.
\end{lemma}
%%%%%%%%%%%%%%%%%%%%%%%%%%%%%%%%%%%%%%%%%%%%%%%%%%%%%%%%%%%%%%%%%%%%%%%%%%

%%%%%%%%%%%%%%%%%%%%%%%%%%%%%%%%%%%%%%%%%%%%%%%%%%%%%%%%%%%%%%%%%%%%%%%%%%
\prova If $M\in\BB_{\kappa}$ then $\widehat M=\bigoplus_{j} M_{j}$,
with $M_{j}$ a block matrix with dimension $d_{j}$ depending on $j$,
and $M_{j}(a,b)=M^{\s,\s'}_{\nn,\nn'}$, for suitable $\nn,\nn',\s,\s'$
such that $|M^{\s,\s'}_{\nn,\nn'}| \le
D \eps_{0} {\rm e}^{-\kappa |\nn-\nn'|^{\rho}}$ for
some constant $D$. Therefore
\begin{eqnarray}
\left\| M_{j} \right\|_{2}^{2} & \!\! = \!\! &
\max_{|x|_{2} \le 1}
\left| M_{j} x \right|_{2}^{2} \le \max_{|x|_{2} \le 1} 
\sum_{a,b,c=1}^{d_{j}} \left| M_{j}(a,b)\right| \left| x(b) \right|
\left| M_{j}(a,c) \right| \left| x(c) \right|
\nonumber \\
& \!\! \le \!\! & \frac{1}{2} \max_{|x|_{2} \le 1} 
\sum_{a,b,c=1}^{d_{j}}
\left| M_{j}(a,b) \right| \left| M_{j}(a,c) \right|
\left( \left| x(b) \right|^{2} + \left| x(c) \right|^{2} \right)
\nonumber \\
& \!\! \le \!\! & \max_{|x|_{2} \le 1} 
\sum_{a=1}^{d_{j}} \left| M_{j}(a,b) \right|
\sum_{c=1}^{d_{j}} \left| M_{j}(a,c) \right|
\sum_{b=1}^{d_{j}} \left| x(b) \right|^{2} \le
\left( \sum_{a=1}^{d_{j}} \left| M_{j}(a,b) \right| \right)^{2} , \nonumber
\end{eqnarray}
which yields the assertion. \EP
%%%%%%%%%%%%%%%%%%%%%%%%%%%%%%%%%%%%%%%%%%%%%%%%%%%%%%%%%%%%%%%%%%%%%%%%%%

%%%%%%%%%%%%%%%%%%%%%%%%%%%%%%%%%%%%%%%%%%%%%%%%%%%%%%%%%%%%%%%%%%%%%%%%%%
\begin{lemma} \label{lem:18}
Let $A,B$ be two self-adjoint $d\times d$ matrices. Then
$$ \left| \la^{(a)}(A+B)-\la^{(a)}(A) \right| \le
\sum_{b=1}^{d} \left| \la^{(b)}(B) \right| $$
for all $a=1,\ldots,d$.
\end{lemma}
%%%%%%%%%%%%%%%%%%%%%%%%%%%%%%%%%%%%%%%%%%%%%%%%%%%%%%%%%%%%%%%%%%%%%%%%%%

%%%%%%%%%%%%%%%%%%%%%%%%%%%%%%%%%%%%%%%%%%%%%%%%%%%%%%%%%%%%%%%%%%%%%%%%%%
\prova The result follows from Lidskii's lemma; cf. \cite{Ka}. \EP
%%%%%%%%%%%%%%%%%%%%%%%%%%%%%%%%%%%%%%%%%%%%%%%%%%%%%%%%%%%%%%%%%%%%%%%%%%

Define $\gotE_{1}=\{\eps\in\gotE_{0} : x_{\nn}(\eps) \ge
2\g/p_{\nn}^{\tau}(\eps) \; \forall \nn\in\gotO\}$ and
$\gotE_{2}=\{\eps\in\gotE_{0} : ||\de_{\nn}(\eps)|-\bar\g| \ge
2\g/|\nn|^{\tau_{1}} \; \forall \nn\in\gotO\}$, and set
$\gotE = \gotE_{1} \cap \gotE_{2}$.

We can denote by $\la^{\s}_{\nn}(A)$, with $\nn\in\gotO$ and
$\s=\pm$, the eigenvalues of the block matrix $A=\DD +\widehat M$.
If $|\delta_{\nn}(\eps)| \ge \bar\g$ then $\la^{\s}_{\nn}(\eps)=
\delta_{\nn}(\eps)$. Moreover for each $\eps\in\gotE_{0}$
and each $\nn\in\gotO$ such that $|\de_{\nn}(\eps)|<\bar\g$,
there exists a block $A^{\nn}(\eps)$ of the matrix $A$, of size
$d_{\nn}(\eps) \le 2C_{1} p_{\nn}^{\al}(\eps)$ such that
$\la^{\pm}_{\nn}(A)$ depends only on the entries of such a block.
This follows from Remarks \ref{rmk:6} and \ref{rmk:12}.

Therefore we have to discard from $\gotE_{0}$ only
values of $\eps$ such that $|\de_{\nn}(\eps)| < \bar\g$
for some $\nn\in\gotO$: for all such $\nn$ the matrix
$A^{\nn}(\eps)$ is well defined, and one has
$\la^{\s}_{\nn}(A)=\la^{\s}_{\nn}(A_{\nn}(\eps))$.

One has, by item 3 in Lemma \ref{lem:4},
\begin{equation}
x_{\nn}(\eps) \ge \frac{1}{p_{\nn}^{\xi}(\eps)}
\min_{a=1,\ldots,d_{\nn}(\eps)}
\left| \la^{(a)} (A^{\nn}(\eps)) \right| \ge
\frac{1}{p_{\nn}^{\xi}(\eps)} \min_{\nn'\in\overline\CC_{\nn}(\eps)}
\min_{\s=\pm} \left| \la^{\s}_{\nn'}(A^{\nn}(\eps)) \right| ,
\label{eq:7.1} \end{equation}
so that, by using that
$\la^{\s}_{\nn'}(A^{\nn}(\eps)) =\la^{\s}_{\nn'}(A^{\nn'}(\eps)) =
\la^{\s}_{\nn'}(A)$ for all $\nn'\in\overline\CC_{\nn}(\eps)$,
we shall impose the conditions
\begin{equation}
\left| \la^{\s}_{\nn}(A^{\nn}(\eps)) \right| \ge
\frac{\g_{2}}{|\nn|^{\tau_{2}}} , \qquad \nn \in \gotO ,
\qquad \s = \pm , \label{eq:7.2} \end{equation}
for suitable $\g_{2}>2\g$.
Thus, the conditions (\ref{eq:7.2}), together with the bound
$|\nn| \le 2 p_{\nn}(\eps)$ (cf. Remark \ref{rmk:18}), will imply
through (\ref{eq:7.1}) the bounds (\ref{eq:5.3}) for $x_{\nn}(\eps)$.

Define
\begin{equation}
\gotK_{\nn}^{\s} = \left\{ \eps \in \gotE_{0} :
\left| \la^{\s}_{\nn} (A) \right|
\le \frac{\g_{2}}{|\nn|^{\tau_{2}}} \right\} ,
\qquad \nn\in\gotO , \qquad \s = \pm , 
\label{eq:7.3} \end{equation}
with $\tau_{2}=\tau-\xi$, so that we can estimate
\begin{equation}
{\rm meas}(\gotE_{0}\setminus\gotE_{1}) \le
\sum_{\nn\in\gotO} \sum_{\s=\pm} {\rm meas}(\gotK^{\s}_{\nn}) .
\label{eq:7.4} \end{equation}
Moreover, by defining
\begin{equation}
\gotH_{\nn,\s} = \left\{ \eps \in \gotE_{0} :
\left| \de_{\nn}(\eps) - \s \bar\g \right|
\le \frac{2\g}{|\nn|^{\tau_{1}}}
\right\} , \qquad \nn\in\CC_{j} , \quad j \in \NNN, \quad \s=\pm ,
\label{eq:7.5} \end{equation}
with $\tau_{1}$ to be determined, one has
\begin{equation}
{\rm meas}(\gotE_{0}\setminus\gotE_{2}) \le \sum_{\nn\in\ZZZ^{D+1}}
\sum_{\sigma=\pm} {\rm meas}(\gotH_{\nn,\s}) .
\label{eq:7.6}
\end{equation}
%

%%%%%%%%%%%%%%%%%%%%%%%%%%%%%%%%%%%%%%%%%%%%%%%%%%%%%%%%%%%%%%%%%%%%%%%%%%
\begin{lemma} \label{lem:19}
There exists constants $w_{0}$ and $w_{1}$ such that
$\gotK^{\pm}_{\nn}=\emptyset$ for all $\nn$ such that
$|\nn| \le w_{0} /\eps_{0}^{w_{1}}$.
There exists constants $y_{0}$ and $y_{1}$ such that
$\gotH_{\nn,\pm}=\emptyset$ for all $\nn$ such that
$|\nn| \le y_{0} /\eps_{0}^{y_{1}}$.
\end{lemma}
%%%%%%%%%%%%%%%%%%%%%%%%%%%%%%%%%%%%%%%%%%%%%%%%%%%%%%%%%%%%%%%%%%%%%%%%%%

%%%%%%%%%%%%%%%%%%%%%%%%%%%%%%%%%%%%%%%%%%%%%%%%%%%%%%%%%%%%%%%%%%%%%%%%%%
\prova We start by considering the sets $\gotK^{\s}_{\nn}$
for $\nn\in\gotO$ and $\s=\pm$.
If $|\de_{\nn}(\eps)| < \bar\g$ one can write $A^{\nn}(\eps)
={\rm diag}\{\de_{\nn'}(0),\de_{\nn'}(0)\}_{\nn'\in\overline
\CC_{\nn}(\eps)} + B^{\nn}(\eps)$,
which defines the matrix $B^{\nn}(\eps)$ as
\begin{equation}
B^{\nn}(\eps) = {\rm diag}\{\de_{\nn'}(\eps)- \de_{\nn'}(0)
\}_{\nn'\in\overline\CC_{\nn}(\eps)}^{\s=\pm} + M^{\nn}(\eps) , \nonumber
\end{equation}
where $M^{\nn}(\eps)$ is the block of $M(\eps)$ with
entries $M^{\s_{1},\s_{2}}_{\nn_{1},\nn_{2}}(\eps)$ such that
$\nn_{1},\nn_{2} \in \overline\CC_{\nn}(\eps)$.
By Lemma \ref{lem:18}, one has 
\begin{equation}
\left| \la^{\s}_{\nn}(A) - \de_{\nn}(0) \right| 
\le \sum_{a=1}^{d_{\nn}(\eps)} \left| \la^{(a)}(B^{\nn}(\eps)) \right|
\le 2 C_{1} p_{\nn}^{\al}(\eps) \|B^{\nn}(\eps)\|_{2} ,
\qquad \nn\in\gotO , \qquad \s = \pm , 
\label{eq:7.7} \end{equation}
where we have used Remark \ref{rmk:10} to bound $d_{\nn}(\eps)$.

One has $|\de_{\nn}(0)| \ge \g_{0}/|\nn|^{\tau_{0}}
\ge \g_{0}/(2p_{\nn}(\eps))^{\tau_{0}}$ by item 2 in
Hypothesis \ref{hp:1}, whereas $\|B^{\nn}(\eps)\|_{2}
\le c_{2} (2p_{\nn}(\eps))^{c_{0}} \eps_{0} + \|M^{\nn}(\eps)\|_{2}$,
by items 1 and 2 in Hypothesis \ref{hp:1},
and $\|M^{\nn}(\eps)\|_{2} \le \|M(\eps)\|_{2} \le
C_{0}\eps_{0}$ by Lemma \ref{lem:17}.
Therefore (\ref{eq:7.7}) implies
\begin{equation}
\left| \la^{\s}_{\nn}(A) \right| \ge
\frac{\g_{0}}{(2p_{\nn}(\eps))^{\tau_{0}}} - C
p_{\nn}^{c_{0}+1}(\eps)\, \eps_{0} , \nonumber
\end{equation}
for a suitable constant $C$, so that, by setting $w_{1}=
c_{0}+1+\tau_{0}$ and choosing suitably the constants
$\g_{2}$, $\tau$ and $w_{0}$,
one has $| \la^{\pm}_{\nn}(A) | \ge \g_{0}/2(2p_{\nn}(\eps))^{\tau_{0}}
\ge \g_{2}/ p_{\nn}^{\tau_{2}}(\eps)$ for all $\nn$ such that
$|\nn| \le w_{0}/\eps_{0}^{w_{1}}$.

For the sets $\gotH_{\nn,\s}$, one can reason in the same way,
by using that $\bar\g\in\gotG$ (cf. Definition \ref{def:4}). \EP
%%%%%%%%%%%%%%%%%%%%%%%%%%%%%%%%%%%%%%%%%%%%%%%%%%%%%%%%%%%%%%%%%%%%%%%%%%

%%%%%%%%%%%%%%%%%%%%%%%%%%%%%%%%%%%%%%%%%%%%%%%%%%%%%%%%%%%%%%%%%%%%%%%%%%
\begin{lemma} \label{lem:20}
Let $\xi>\xi_{1}$ and $\eps_{0}=\eta_{0}^{N}$
be fixed as in Lemma \ref{lem:14}.
There exist constants $\g$, $\tau$ and $\tau_{1}$ such that
${\rm meas}(\gotE_{0}\setminus\gotE) = o(\eps_{0})$.
\end{lemma}
%%%%%%%%%%%%%%%%%%%%%%%%%%%%%%%%%%%%%%%%%%%%%%%%%%%%%%%%%%%%%%%%%%%%%%%%%%

%%%%%%%%%%%%%%%%%%%%%%%%%%%%%%%%%%%%%%%%%%%%%%%%%%%%%%%%%%%%%%%%%%%%%%%%%%
\prova First of all we have to discard from $\gotE_{0}$ the sets
$\gotH_{\nn,\s}$. It is easy to see that one has
\begin{equation}
{\rm meas}(\gotH_{\nn,\s}) \le \frac{2\g}{|\nn|^{\tau_{1}}} \,
\frac{2}{c_{1}|\nn|^{c_{0}}} , \nonumber
\end{equation}
for some positive constant $C$, so that,
by using the second assertion in Lemma \ref{lem:19}, we find
\begin{equation}
\sum_{\nn\in\gotO} \sum_{\sigma=\pm1}
{\rm meas}(\gotH_{\nn,\s}) \le
\sum_{\substack{\nn\in\gotO \\ |\nn| \ge y_{0}/\eps_{0}^{y_{1}} }}
\sum_{\sigma=\pm 1} {\rm meas}(\gotH_{\nn,\s}) \le
C \eps_{0}^{y_{1}(\tau_{1}+c_{0}-D-1)} , \nonumber
\end{equation}
for some constant $C$, provided $\tau_{1}+c_{0}-D>1$,
so that we shall require for $\tau_{1}$
to be such that $\tau_{1}+c_{0}-D>1$ and
$y_{1}(\tau_{1}+c_{0}-D-1)>1$.

Next, we consider the sets $\gotK^{\pm}_{\nn}$. For all $\nn\in\gotO$
consider $A^{\nn}(\eps)$ and write $A^{\nn}(\eps)=\de_{\nn}(\eps)I +
B^{\nn}(\eps)$, which defines the matrix $B^{\nn}(\eps)$ as
\begin{equation}
B^{\nn}(\eps) = {\rm diag}\{\de_{\nn'}(\eps)-
\de_{\nn}(\eps)\}_{\nn\in\overline\CC_{\nn}(\eps)}^{\s=\pm} +
M^{\nn}(\eps) , \nonumber
\end{equation}
with $M^{\nn}(\eps)$ defined as in the proof of Lemma \ref{lem:19}.

Then the eigenvalues of $A^{\nn}(\eps)$ are of the
form $\la^{(a)}(A^{\nn}(\eps)) =\de_{\nn}(\eps) +
\la^{(a)}(B^{\nn}(\eps))$, so that for all $\eps\in\gotE_{0}
\setminus \overline \II_{\nn}(\g)$ one has
$$ \left| \partial_{\eps}\la^{(a)}(A^{\nn}) \right| \ge
\left| \partial_{\eps} \de_{\nn}(\eps) \right| -
\left\| \partial_{\eps} B^{\nn}(\eps) \right\|_{2} , $$
where item 3 in Lemma \ref{lem:16} has been used.
One has $|\partial_{\eps} \de_{\nn}(\eps)| \ge
c_{1} |\nn|^{c_{0}}$, by item 2 in Hypothesis \ref{hp:1}, and
$\| \partial_{\eps} B^{\nn}(\eps) \|_{2} \le \max_{\nn'\in\overline
\CC_{\nn}(\eps)} |\pr_{\eps} (\de_{\nn'}(\eps)-\de_{\nn}(\eps))| +
\|\pr_{\eps} M^{\nn}(\eps)\|_{2} \le \zeta c_{3} p_{\nn(\eps)}\,
p_{\nn}^{c_{0}-1}(\eps) + \eps_{0} C p_{\nn}^{c_{0}}(\eps)$,
for a suitable constant $C$, as follows from item 4 in
Hypothesis \ref{hp:1}, from Hypothesis \ref{hp:3}
(see Remark \ref{rmk:6} for the definition of $\zeta$),
from Lemma \ref{lem:14}, and from Lemma \ref{lem:15}.
Hence we can bound $| \partial_{\eps}\la^{(a)}(A^{\nn})| \ge
c_{1} |\nn_{0}|^{c_{0}}/2$ for $\eps_{0}$ small enough.

Therefore one has
\begin{equation}
{\rm meas}(\gotK^{\s}_{\nn}) \le
\frac{2\g_{2}}{|\nn|^{\tau_{2}}(\eps)}
\frac{2}{c_{1} |\nn|^{c_{0}}} \left( C |\nn|^{(\al+\be)(D+1)} \right) ,
\label{eq:7.8} \end{equation}
for some constant $\overline C$, where the last factor
$\overline C|\nn|^{(\al+\be)(D+1)}$ arises for the following reason.
The eigenvalues $\la^{\s}_{\nn}(A)$ are differentiable in $\eps$
except for those values $\eps$ such that for some $\nn'\in\CC_{\nn}$
one has $|\de_{\nn'}(\eps)|=\bar\g$ and $|\de_{\nn}(\eps)|<\bar\g$.
Because of item 3 in Hypothesis \ref{hp:1} all functions
$\de_{\nn'}(\eps)$ are monotone in $\eps$ as far as
$|\de_{\nn'}(\eps)|<1/2$, so that for each $\nn'\in\CC_{\nn}$
the condition $|\de_{\nn'}(\eps)|=\bar\g$ can occur at most twice.
The number of $\nn'\in\CC_{\nn}$ such that the conditions
$|\de_{\nn'}(\eps)|=\bar\g$ and $|\de_{\nn}(\eps)|<\bar\g$
can occur for some $\eps\in\gotE_{0}$ is bounded by
the volume of a sphere of centre $\nn$ and radius
proportional to $|\nn|^{\al+\be}$ (cf. Remark \ref{rmk:6}).
Hence $\overline C|\nn|^{(\al+\be)(D+1)}$ counts the number of
intervals in $\gotE_{0}\setminus \overline\II_{\nn}(\g)$.

Thus, (\ref{eq:7.8}) yields, by making use of the first assertion
of Lemma \ref{lem:19},
\begin{equation}
\sum_{\nn\in\gotO} \sum_{\s=\pm}
{\rm meas}(\gotK^{\s}_{\nn}) \le
\sum_{\substack{\nn\in\ZZZ^{D+1} \\ |\nn| \ge w_{0}/\eps_{0}^{w_{1}} }}
\frac{8\g}{c_{1}}
|\nn|^{-\tau_{2}-c_{0}} \left( \overline C|\nn|^{2\al} \right)
\le C \eps_{0}^{w_{1}(\tau_{2}+c_{0}-2\al-D-1)} , \nonumber
\end{equation}
for some positive constant $C$, provided $\tau_{2}+c_{0}-2\al-D=
\tau+c_{0}-2\al-D-\xi>1$, so that (\ref{eq:7.4}) implies that
${\rm meas}(\gotE_{0}\setminus\gotE_{1}) \le C \eps_{0}^{w_{1}
(\tau_{2}+ c_{0}-D-1)}$.

Therefore, the assertion follows provided
$\min\{\tau_{1},\tau_{2}-2\al\}> D-c_{0}+1$,
$y_{1}(\tau_{1}+c_{0}-D-1) > 1$ and
$w_{1}(\tau_{2}+c_{0}-2\al-D-1) > 1$.\EP
%%%%%%%%%%%%%%%%%%%%%%%%%%%%%%%%%%%%%%%%%%%%%%%%%%%%%%%%%%%%%%%%%%%%%%%%%%

\appendix

%%%%%%%%%%%%%%%%%%%%%%%%%%%%%%%%%%%%%%%%%%%%%%%%%%%%%%%%%%%%%%%%%%%%%%%%%
%%%%%%%%%%%%%%%%%%%%%%%%%%%%%%%%%%%%%%%%%%%%%%%%%%%%%%%%%%%%%%%%%%%%%%%%%
\zerarcounters
\section{Proof of Lemma \ref{lem:1}}
\label{app:A}
%%%%%%%%%%%%%%%%%%%%%%%%%%%%%%%%%%%%%%%%%%%%%%%%%%%%%%%%%%%%%%%%%%%%%%%%%

Lemma \ref{lem:1} is a consequence of the following elementary
proposition in Galois theory.

%%%%%%%%%%%%%%%%%%%%%%%%%%%%%%%%%%%%%%%%%%%%%%%%%%%%%%%%%%%%%%%%%%%%%%%%%
\begin{prop}
If $p_1,\dots,p_k$  are distinct primes  then the field $$F:=\mathbb
Q[\sqrt{p_1},\sqrt{p_2},\dots,\sqrt{p_k}]$$ obtained from the rational
numbers  $\mathbb Q$ by adding the  $k$  square roots $\sqrt{p_i}$ has
dimension $2^k$ over $\mathbb Q$ with basis  the elements $\prod_{i\in
I}\sqrt{p_i}$ as $I$ varies on the $2^k$ subsets of $\{1, 2,\dots,k\}$.

The group of automorphisms\footnote{i.e. the linear transformations
$\tau$ such that $\tau(uv)=\tau(u)\tau(v)$.} of $F$ which fix $\mathbb Q$
(i.e. the Galois group of $F/\mathbb Q$) is an Abelian group generated
by  the automorphisms $\tau_{i}$ defined by $\tau_{i}(\sqrt{p_j})=
(-1)^{\delta(i,j)}\sqrt{p_j}$.
\end{prop}
%%%%%%%%%%%%%%%%%%%%%%%%%%%%%%%%%%%%%%%%%%%%%%%%%%%%%%%%%%%%%%%%%%%%%%%%%

%%%%%%%%%%%%%%%%%%%%%%%%%%%%%%%%%%%%%%%%%%%%%%%%%%%%%%%%%%%%%%%%%%%%%%%%%
\prova We prove by induction both statements. Let us assume
the statements valid for $p_1,\dots,p_{k-1}$ and let $F':=\mathbb
Q[\sqrt{p_1},\sqrt{p_2},\dots,\sqrt{p_{k-1}}]$ so that
$F=F'[\sqrt{p_k}]$. We first prove that $\sqrt{p_k}\notin F'$.
Assume it to be false.
Since $(\sqrt{p_k})^2$ is integer, each element -- say $\tau$ -- of
the Galois group of $F'/\mathbb Q$ must either fix $\sqrt{p_k}$ or
transform it into $-\sqrt{p_k}$ (by definition $\tau(p_k)=
\tau(\sqrt p_k)^2=p_k $).

Now any element $b\in F'$ is by induction uniquely expressed as
$$ b= \sum_{I\subset \{1,2,\dots,k-1\}} a_{I} \prod_{i\in I}
\sqrt{p_i}, \qquad a_I\in \mathbb Q . $$ 
If $h$ of the numbers $a_I$ are non zero, it is easily seen
that $b$ has $2^h$ transforms (changing the signs of each of
the $a_I$) under the Galois group of $F'$.  Therefore $b=\sqrt p_k$
if and only if  $h=1$, that is one should have $\sqrt{p_k} = 
m/n \prod_{i\in I}\sqrt{p_i}, I\subset \{1,2,\dots,k-1\}$
for $m,n$ integers. This implies that  $p_kn^2=m^2\prod_{i\in I} p_i$
which is impossible by the unique factorisation of integers.
This proves the first statement.

To construct the Galois group of $F/\mathbb Q$ we extend the action
of $\tau_{i}$ for $i=1,\ldots,k-1$ by setting $\tau_{i}(\sqrt{p_{k}})
=\sqrt{p_{k}}$. Finally we define the automorphism $\tau_{k}$
as $\tau_{k}(\sqrt{p_j})= (-1)^{\delta(k,j)}\sqrt{p_j}$
for $j=1,\ldots,k$. \EP

\*

\noindent\textbf{Acknowledgements}. We thank Claudio Procesi
and Massimiliano Berti for useful discussions.

%%%%%%%%%%%%%%%%%%%%%%%%%%%%%%%%%%%%%%%%%%%%%%%%%%%%%%%%%%%%%%%%%%%%%%%%%%
%%%%%%%%%%%%%%%%%%%%%%%%%%%%%%%%%%%%%%%%%%%%%%%%%%%%%%%%%%%%%%%%%%%%%%%%%%
% References
%%%%%%%%%%%%%%%%%%%%%%%%%%%%%%%%%%%%%%%%%%%%%%%%%%%%%%%%%%%%%%%%%%%%%%%%%%
%%%%%%%%%%%%%%%%%%%%%%%%%%%%%%%%%%%%%%%%%%%%%%%%%%%%%%%%%%%%%%%%%%%%%%%%%%

\end{document}